\documentclass[10pt]{article}
\usepackage{amsmath, amssymb,amsthm, fullpage, color, relsize, enumerate, enumitem, hyperref, mathtools}
\usepackage{latexsym}
\usepackage{bm,bbm} 
\usepackage{graphicx}
\usepackage{chngcntr}
\usepackage{import}
\usepackage{comment}
\usepackage{amsfonts}
\usepackage[T2A, T1]{fontenc}
\usepackage[utf8]{inputenc}
\usepackage{mathtools}  
\usepackage{multirow,multicol}
\usepackage[noadjust]{cite}
\usepackage{mcite}
\usepackage{subfig}
\usepackage{caption}
\usepackage[russian, english]{babel}
\numberwithin{equation}{section}
\usepackage{authblk}
\usepackage[toc,page]{appendix}
\usepackage{pgfplots}
\pgfplotsset{compat=newest} 
\pgfplotsset{plot coordinates/math parser=false}

\title{Optimal change point detection in Gaussian processes}

\author[1]{Hossein Keshavarz}
\author[2,1]{Clayton Scott}
\author[1,2]{XuanLong Nguyen}

\affil[1]{Department of Statistics, University of Michigan}
\affil[2]{Department of Electrical Engineering and Computer Science, University of Michigan}
\affil[ ]{Email: \textit {\{hksh,clayscot,xuanlong\}@umich.edu}}

\theoremstyle{definition}
\newtheorem{defn}{Definition}[section]

\newtheorem{prop}{Proposition}[section]

\newtheorem{lem}{Lemma}[section]
\newtheorem{thm}{Theorem}[section]
\newtheorem{cor}{Corollary}[section]
\newtheorem{rem}{Remark}[section]
\newtheorem{assu}{Assumption}[section]

\makeatletter
\@addtoreset{clawithinpf}{proof}
\makeatother

\long\def\comment#1{}

\DeclareSymbolFont{cyrillic}{T2A}{cmr}{m}{n}
\def\makecyrsymbol#1#2{%
  \begingroup\edef\temp{\endgroup
    \noexpand\DeclareMathSymbol{\noexpand#1}
    {\noexpand\mathalpha}{cyrillic}%
    {\expandafter\expandafter\expandafter
     \calccyr\expandafter\meaning\csname T2A\string#2\endcsname\end}}%
  \temp}
\expandafter\def\expandafter\calccyr\string\char#1\end{#1}

\makecyrsymbol\russz\cyrz
\makecyrsymbol\RussL\CYRL
\makecyrsymbol\RussD\CYRD

\makeatletter
\def\amsbb{\use@mathgroup \M@U \symAMSb}
\makeatother

\makeatletter
\g@addto@macro\@floatboxreset\centering
\makeatother

\newcommand{\bb}[1]{\mathbb{#1}} 
\newcommand{\bbM}[1]{\mathbbm{#1}} 
\newcommand{\cc}[1]{\mathcal{#1}} 
\newcommand{\ff}[1]{\mathfrak{#1}} 
\newcommand{\zero}{\boldsymbol{0}}

\newcommand{\sign}{\mathop{\mathrm{sign}}}

\newcommand{\supp}{\mathop{\mathrm{supp}}}

\newcommand{\cov}{\mathop{\mathrm{cov}}} 
\newcommand{\var}{\mathop{\mathrm{var}}} 
\newcommand{\sinc}{\mathop{\mathrm{sinc}}} 
 
\newcommand{\tap}{\mathop{\mathrm{tap}}}
\newcommand{\dist}{\mathop{\mathrm{dist}}} 
 
\newcommand{\essinf}{\mathop{\mathrm{essinf}}}

\newcommand{\brac}[1]{\left[#1\right]}
\newcommand{\set}[1]{\left\{#1\right\}}
\newcommand{\abs}[1]{\left\lvert #1 \right\rvert}
\newcommand{\paren}[1]{\left(#1\right)}
\newcommand{\brapar}[1]{\left[#1\right)}
\newcommand{\parbra}[1]{\left(#1\right]}

\newcommand{\Bigbrac}[1]{\Big[#1\Big]}

\newcommand{\InnerProd}[2]{\langle #1,#2 \rangle}

\newcommand{\cp}[1]{\overset{#1}{\rightarrow}}    
  
\newcommand{\RelNum}[2]{\overset{#1}{#2}}   

\newcommand{\floor}[1]{\left\lfloor#1\right\rfloor}

\newcommand{\LpNorm}[2]{\left\| #1\right\|_{\ell_{#2}}}

\newcommand{\norm}[2]{\left\|#1\right\|_{#2}}
\newcommand{\OpNorm}[3]{\left\| #1\right\|_{#2\rightarrow#3}}
\newcommand{\RNum}[1]{\uppercase\expandafter{\romannumeral #1\relax}}

\pgfplotsset{every axis/.append style={
        scaled y ticks = false, 
        scaled x ticks = false, 
        y tick label style={/pgf/number format/.cd, fixed, fixed zerofill,
                            ,precision=1},
        x tick label style={/pgf/number format/.cd, fixed, fixed zerofill,
                            precision=2}
    }
}

\newlength\figureheight 
\newlength\figurewidth

\begin{document}
\maketitle
\begin{abstract}

We study the problem of detecting a change in the mean of one-dimensional Gaussian process data. This problem is investigated in the setting of increasing domain (customarily employed in time series analysis) and in the setting of fixed domain (typically arising in spatial data analysis). We propose a detection method based on the generalized likelihood ratio test (GLRT), and show that our method achieves nearly asymptotically optimal rate in the minimax sense, in both settings. The salient feature of the proposed method is that it exploits in an efficient way the data dependence captured by the Gaussian process covariance structure. When the covariance is not known, we propose the plug-in GLRT method and derive conditions under which the method remains asymptotically near optimal. By contrast, the standard CUSUM method, which does not account for the covariance structure, is shown to be asymptotically optimal only in the increasing domain. Our algorithms and accompanying theory are applicable to a wide variety of covariance structures, including the Matern class, the powered exponential class, and others. The plug-in GLRT method is shown to perform well for maximum likelihood estimators with a dense covariance matrix.
\end{abstract}

\section{Introduction}

Change point detection is the problem of detecting an abrupt change or changes arising in a sequence of observed samples. A common problem of this type involves detecting shifts in the mean of a temporal process. This problem has found a variety of applications in many fields, including audio analysis~\cite{O.Gillet}, EEG segmentation~\cite{M.Lavielle}, structural health monitoring~\cite{HY.Noh,X.Hu} and environment sciences~\cite{M.Last, J.Verbesselt}. Despite advances in the development of algorithms~\cite{Y.Kawahara, M.Lavielle, J.Liu, G.Rigaill} and asymptotic theory~\cite{P.Bertrand, A.Tartakovsky, X.Shao, C.Levy-leduc} for a number of contexts, such studies are mainly confined to the setting of (conditionally) independently distributed data. Existing works on optimal detection of shifts in the mean in temporal data with statistically dependent observations are far less common.

Incorporating dependence structures into the modelling of random processes is a natural approach. In fact, this has been considered in detecting changes of remotely collected data~\cite{V.Chandola, E.Gabriel}. For instance, Chandola et al.~\cite{V.Chandola} proposed a Gaussian process based algorithm to identify changes in Normalized Difference Vegetation Index (NDVI) time series for a particular location in California. Despite such statistical modelling considerations, to our knowledge, most researchers have not exploited the dependence structures of the underlying temporal process, e.g., its covariance function and spectral density, in designing a (minimax) optimal detection method. 

In this paper, we focus on the detection of a single change in the mean of a Gaussian process data sequence. Our main contribution is to show that neglecting the dependence structures in the data samples leads to suboptimal detection procedures, particularly in the presence of strong correlation among samples. Moreover, it is possible to exploit the underlying dependence structures to design asymptotically optimal detection algorithms.

Consider a simplified setting in which we let $G$ be a Gaussian process in a domain $\cc{D}\subseteq\bb{R}$ and $\cc{D}_n \coloneqq \set{t_k}^n_{k=1}\subset\cc{D}$ denote the finite index set of sampling points. Denote the observed samples by $\boldsymbol{X} = \set{X_k}^n_{k=1}$ in which $X_k = G\paren{t_k}$ for $k=1,\ldots,n$. Moreover, let $t\in\cc{C}_{n,\alpha}\subseteq\set{1,\ldots,n}$ (the parameter $\alpha$ is a positive scalar which will be introduced in Section \ref{DtctProcedure}) and $b>0$ represents the point of sudden change and the jump/shift value, respectively. Namely, there is $\mu\in\bb{R}$ (which will be assumed to be 0 for now) such that
\begin{equation}\label{ProbForm}
\bb{E}X_k = \paren{\mu-\frac{b}{2}}\bbM{1}\paren{k\leq t} + \paren{\mu+\frac{b}{2}}\bbM{1}\paren{k>t},\quad k\in\set{1,\ldots,n}.
\end{equation}
To analyze the performance of a detection procedure as sample size $n$ grows to infinity, one is confronted with two fundamentally different theoretical frameworks: \emph{increasing domain asymptotics} and \emph{fixed domain (infill) asymptotics} (\cite{TS.Rao}, Chapter 13). The former arises naturally in time series analysis, which is distinguished by the constraint that the distance between consecutive sampling \emph{time} points are bounded away from zero. The simplest instance of the sampling scenario in this regime arises when the diameter of $\cc{D}_n$ is of order $n$ and $\min\abs{t_{i+1}-t_i} > \delta$ for some strictly positive, fixed scalar $\delta$. In our notation the index set for the Gaussian process represents the sampling time points. Typically we set $\cc{D}=\bb{R}$ and $\bigcup^{\infty}_{n=1}\cc{D}_n=\bb{N}$ or $\bb{Z}$. See, \cite{J.Antoch, L.Horvath, L.Horvath2, P.Kokoszka, M.Rencova, Y.Yao} for examples of works studying change point detection via increasing domain asymptotics. 

Fixed domain asymptotics, one the other hand, is a more suitable setting when the index set of sampling points $\cc{D}$ is bounded, so that the observations get denser in $\cc{D}$ as $n$ increases. Particularly for $\cc{D}\subset\bb{R}$, we have that $\min \abs{t_{i+k}-t_i} = \cc{O}\paren{k/n}$ for positive integers $i,k$ with $i,\paren{i+k}\in\set{1,\ldots,n}$, and it can be extended to multidimensional domains in a straightforward way. This is the case for spatially distributed data \cite{ML.Stein}, where the domain of the index set is typically of one, two or three dimensions. But there are also other examples: a particularly useful approach to change point detection in non-stationary processes widely adopted in speech processing and finance can also be cast in this framework. In fact, this approach involves \emph{piecewise locally stationary (PLS) processes}, which can be interpreted as processes that are approximately piecewise stationary, due to the gradual and smooth change of the spectrum \cite{S.Adak, R.Dahlhaus}. Abrupt change detection in PLS processes has been considered in, for instance, \cite{S.Adak, M.Last}. Since the future samples of a PLS process may not carry any information about the current state of the process, the increasing domain setting cannot suitably capture the dependence structure of the observed data sequence. As a result, the theory of change detection and inference in PLS processes is also based upon fixed domain asymptotics.

Our goal in this work is to derive a computationally efficient detection procedure that effectively accounts for the underlying dependence structures of the observed sequential data. For the analysis of such a procedure we shall adopt both aforementioned asymptotic frameworks. Obviously, specific applications ultimately determine which one among the two is more suitable. There are also several motivations for this dual treatment. First, we are not aware of any established detection algorithm and accompanying asymptotic theory under the fixed-domain setting. Second, although there is a vast literature on the change point detection in the increasing domain scenario, existing work focus mostly on independent (or conditionally independent) data sequences, so our analysis in this setting still carries some notable novelty. Third, there is a fundamental difference in the behavior of the same detection algorithm when applied to the two asymptotic settings. This point is worth highlighting: we will show that in the fixed domain setting, ignoring the underlying dependence structure may result in suboptimal detection performance, but this is not the case for the increasing domain setting. Finally, we note that the fixed domain minimax detection problem considered in this paper also serves as a useful starting point in the study of optimal detection of discontinuities in Gaussian spatial processes, as considered recently by \cite{E.Gabriel, X.Shen}. 

The amount of correlation between the nearby samples is one of the differences among the two asymptotic frameworks that will plays a crucial role in our results. In the fixed domain regime, regardless of how large the sample size is, if $\abs{j-i}$ is of order $n^{\beta}$ for some $\beta\in\paren{0,1}$, the correlation among $X_i$ and $X_{j}$ is still so close to one. Roughly speaking, the effective sample size is much smaller than $n$. However, in the increasing domain setting, even for long range dependent processes the correlation among samples at time points $i,j$ is small when $\abs{j-i}$ is large. Accordingly, the asymptotic behavior of a statistical method applied to the fixed-domain regime is expected to be distinct from that of the increasing domain regime. Moreover, new techniques are needed to address the statistical dependence intrinsic in the former.

From here on, for fixed domain asymptotic results, we assume that $G$ is a one dimensional Gaussian process in $\cc{D} = \brac{0,1}$ with $n$ regularly spaced samples, i.e. $\cc{D}_n = \set{k/n}^n_{k=1}$. Under the increasing domain asymptotics framework, we let $\cc{D} = \brapar{0,\infty}$ and $\cc{D}_n = \set{1,\ldots,n}$. The remaining notations are the
same for both.

\paragraph{Previous works.}
One of the earliest attempts to study shift in mean detection was perhaps that of 
Chernoff et al.~\cite{H.Chernoff}. More general settings of this problem have been studied in subsequent works, e.g., \cite{I.B.MacNeill, J.Deshayes, Y.Yao}. For instance it is assumed in \cite{Y.Yao} that the sequence of $X_k$'s are independent Gaussian variables. They proposed a detection method based on generalized likelihood ratio test (GLRT), also known as the \emph{cumulative sum (CUSUM)} test, and given by
\begin{equation}\label{CUSUM}
T_{CUSUM} = \bbM{1}\set{\max_{t\in\cc{C}_{n,\alpha}} \biggr \{
	\sqrt{\frac{t\left(n-t\right)}{n}}\abs{\frac{1}{n-t}\sum\limits_{k=t+1}^{n} X_k - \frac{1}{t}\sum\limits_{k=1}^{t} X_k}  \biggr\} \geq R_n}.
\end{equation}
CUSUM compares the maximum of a test statistic over $C_{n,\alpha}$ with a critical values $R_n$. Non-asymptotic upper bounds on the error probabilities of this simple test were obtained by the authors under the Gaussian and i.i.d. assumptions. Interestingly, due to test's simplicity, when such assumptions do not hold, one can still apply the same test statistic to the data sequence. 
Most subsequent works appeared to follow this direction, in addition to adhering to the increasing domain asymptotic framework, e.g. \cite{J.Antoch,L.Horvath,L.Horvath2,M.Rencova}.
We wish to mention Rencova (\cite{M.Rencova}, chapter 4), who studied the same CUSUM test as \cite{Y.Yao}, but working with the assumption that $\boldsymbol{X}$ be a strong mixing time series. Kokoszka~\cite{P.Kokoszka} also analyzed the CUSUM test, but working with a different dependent observation model with sub-squared growth of the variance of partial sum, i.e. there is $\delta\in\paren{0,2}$ such that for any $k<m$, $\var\sum_{j=k}^{m} X_j\lesssim \paren{m-k+1}^{\delta}$. Horv{\'a}th et al. \cite{L.Horvath,L.Horvath2} and Antoch \cite{J.Antoch} studied the performance of the CUSUM test for the detection of a sudden change in the mean in linear processes, i.e. $X_t = \sum_{j=0}^{\infty} w_j\epsilon_{t-j}$, in which $\set{\epsilon_t}^{\infty}_{t=-\infty}$ are i.i.d. and zero mean random variables and the weights $\set{w_j}^{\infty}_{j=0}$ satisfy some properties such as absolute or square summability.

At first sight, it may seem puzzling how the CUSUM test continues to admit nearly optimal detection performance even as its test statistics apparently ignore the dependence among data samples. A high-level explanation can be made
regarding this phenomenon: when $G$ is a Gaussian process with $\bigcup^{\infty}_{n=1}\cc{D}_n = \bb{N}$ and the covariance function $\cov\paren{X_s,X_t}\rightarrow 0$ as $\abs{t-s}$ grows to infinity, the percentage of pairs $\paren{X_s,X_t}^n_{s,t=1}$ whose covariance is non-negligible tends to zero as $n\rightarrow\infty$. As a result, there is very little gain in accounting for the dependence structures underlying the sequence, and so the CUSUM statistic provides a good approximation of the likelihood ratio test for large $n$, leading to the asymptotic optimality of $T_{CUSUM}$ in the increasing domain setting. This is of course not the case for the fixed domain setting. Indeed, one of the contributions of this paper is to show that the CUSUM test is suboptimal when applied to the fixed domain setting of the detection problem. Moreover, to achieve mimimax optimal detection performance we shall develop a new test statistic that account for the underlying dependence structures in the data sequence. We also note in passing that in comparison to the increasing domain analysis typically encountered in the literature, the theoretical analysis for the fixed domain setting is considerably more challenging, as one needs to take into account the statistical dependence in the data sequence in a more fundamental way.

CUSUM test also applied to one dimensional processes with highly correlated samples, after a proper standardization. For instance, Horv{\'a}th et al. \cite{L.Horvath2} used a different normalizing factor for applying CUSUM to one dimensional Gaussian time series with long range dependence. However apart from standardizing factor, they do not directly incorporate the correlation structures of the data in the formulation of the test statistic. Note that, the proposed test in this paper achieve consistency under weaker condition on minimal detectable jump. Furthermore, Lai \cite{TL.Lai} adjusts CUSUM test for detecting abrupt changes in the mean of a stochastic process. However we are not aware of the analysis of change point problem, for the case of unknown before and after distribution parameters, in the infill regimes and its comparison with the increasing domain setup. There are also some notable works (see e.g., \cite{V.Spokoiny}) on estimating the volatility parameter of non-stationary time series using the change point models.

\paragraph{Overview of main results.}\label{Contributions} 
In this paper we study the change point detection problem that arises in Gaussian processes in both settings of increasing and fixed domains. We show when it is important to account for the dependence structure in the data sequence, and analyze a number of detection algorithms based upon a generalized likelihood ratio test. More specifically, our contributions are as follows.

\begin{enumerate}
\item Given an $n$-sample drawn from a one dimensional Gaussian process data with a known covariance structure, we propose a generalized likelihood ratio test for detecting a sudden shift in the mean. This method requires the knowledge of the dependence structure (via the covariance matrix), and will be shown to achieve asymptotically optimal detection performance
in both fixed and increasing domain settings. Our theory holds for a variety of covariance structures, such as the Matern class, powered exponential class, and several others specified in terms of the covariance kernel's spectral density.
The smoothness parameter for the Gaussian process (which determines how fast the corresponding spectral density decays) plays a central role in characterizing the minimax optimal detection performance --- but this is shown to be the case only in the fixed domain setting, \emph{not} in the increasing domain setting.

\item We provide an upper bound guarantee for the CUSUM detection method. This result confirms that the CUSUM is asymptotically optimal in the increasing domain setting, but it also suggests that the CUSUM is suboptimal in the fixed domain setting. The suboptimality is confirmed in our simulation study, which demonstrates a wide gap between the CUSUM and our method. This result makes sense, in light of the minimax result described in the preceding paragraph.

\item In practice, the covariance structure is not known a priori, and often has to be estimated from the data. To address this scenario, we propose a plug-in GLRT method, and analyze its performance. In particular, we derive detection performance bounds which also account for the quality of a particular covariance estimation method (such as MLE based methods with dense or tapered covariances) used in the plug-in GLRT. Most surprisingly, we show that as long as a consistent covariance estimate is employed (the notion of consistency will be defined in Section \ref{DtctRateUnKnwnSpDen}), regardless of its estimation rate, the plug-in GLRT achieves asymptotically optimal detection performance. Moreover, in some situations a plug-in GLRT with an inconsistent covariance estimate 
is shown to perform almost as well as the case of known covariance.
\end{enumerate}

In addition to the above contributions, our proof methods contain several useful techniques. Our proofs integrate four major technical tools: 
we exploit properties of the mutually absolutely Gaussian measures, 
the decorrelation of samples drawn from Gaussian processes in fixed domain, 
the non-asymptotic analysis of the inverse of large Toeplitz matrices, in addition to the classical theory of minimax detection. We also develop novel techniques for analyzing different norms of the decorrelation matrix of $\boldsymbol{X}$ in either of the asymptotic regimes. The Appendices provide several beneficial and easy-to-reference technicalities for theoretical problems in the area of Gaussian random fields and time series that may be of independent interest.

\paragraph{Structure of the paper.}
Section \ref{DtctAlg} precisely formulates detection of a shift in the mean, focusing on one dimensional Gaussian process data, and then introduces our proposed detection algorithm as hypothesis testing in the case of fully known spectral density of the underlying process. We also adapt our proposed detection technique, which will be referred as the plug-in test, to the much more realistic case of unknown spectral density. Note that Sections \ref{DtctRateKnwnSpDen}, \ref{DtctCUSUM} and \ref{LowBndDtctRate} are divided into two subsections focusing on theoretical results in the fixed domain and the increasing domain setting, respectively. Section \ref{DtctRateKnwnSpDen} presents sufficient conditions on shift value $b$ and spectral density to detect the existence of shift in mean with high probability. Section \ref{DtctRateUnKnwnSpDen} is devoted to analyze the  performance of plug-in test by imposing some sufficient conditions on the estimated spectral density. Section \ref{DtctCUSUM} serves as a comprehensive study of the CUSUM test. The minimax optimality of the proposed algorithms will be discussed in section \ref{LowBndDtctRate}. Section \ref{NumRes} is devoted to the numerical experiments and assessing the proposed algorithms using simulation studies. Section \ref{Discussion} contain concluding remarks with a concise discussion of future directions. Appendix \ref{Proofs} contains the proofs of the main results and Appendix \ref{AuxRes} presents and proves some auxiliary results used in Appendix \ref{Proofs}. Lastly, Appendix \ref{NonAsympInvLargToep} develops non-asymptotic results on the inverse of large Toeplitz matrices which are useful in the study of CUSUM test.

\paragraph{Notation.}
$\wedge$ and $\vee$ stand for minimum and maximum operators and the indicator function is represented by $\bbM{1}\paren{\cdot}$. For any $m\in\bb{N}$, $I_m$, $\zero_m$ and $\bbM{1}_m$ respectively denote the $m$ by $m$ identity matrix, all zeros column vector of length $m$, and all ones column vector of length $m$. For two matrices of the same size $M_1$ and $M_2$, $\InnerProd{M_1}{M_2}{}\coloneqq \sum_{i,j} \paren{M_1}_{ij}\paren{M_2}_{ij}$ denotes their usual inner product. For any symmetric matrix $M$, $\lambda_{\min}\paren{M}$ represents the smallest eigenvalue of $M$.  We will use the following matrix norms on $M\in\bb{R}^{m\times n}$. For any $1\leq p<\infty$, $\LpNorm{M}{p}\; \coloneqq \paren{\sum_{i,j} \abs{M_{ij}}^p}^{1/p}$ stands for the element-wise $\ell_p$ norm of $M$, while $\LpNorm{M}{\infty} \coloneqq \max_{i,j}\abs{M_{ij}}$ represents the sup norm of $M$. For a function $f:\cc{D}\mapsto\bb{R}$ and $p>0$, $\norm{f}{p}^p\coloneqq \int_{\cc{D}}\abs{f\paren{u}}^p du$. The special case of $p=\infty$ is defined by $\norm{f}{\infty} \coloneqq \sup_{u\in\cc{D}}\abs{f\paren{u}}$. For any  $f\in\bb{L}^{1}\paren{\bb{R}}$, $\hat{f}$ represents its Fourier transform defined by 
\begin{equation*}
\hat{f}\paren{\omega} = \int\limits_{-\infty}^{\infty} f\paren{t} e^{-j\omega t} dt,\quad\forall\;\omega\in\bb{R},
\end{equation*} 
where $j^2 = -1$ denotes the imaginary unit. Moreover, for a symmetric function $f\in\bb{L}^{\infty}\paren{\brac{-\pi,\pi}}$, $\set{f_m}_{m\in\bb{Z}}$ denotes the set of the its Fourier coefficients. Assuming that $\set{f_m}_{m\in\bb{Z}}$ is absolutely summable, $\cc{T}_{\bb{N}}\paren{f}$ represents the infinite Toeplitz matrix generated by the Fourier coefficient of $f$, i.e.
\begin{equation}
\cc{T}_{\bb{N}}\paren{f} \coloneqq \paren{\begin{array}{cccc}
f_0 & f_1 & f_2 & \cdots \\
f_1 & f_0 & f_1 & \ddots\\
f_2 & f_1 & f_0 & \ddots \\
\vdots & \ddots & \ddots & \ddots \\
\end{array}}.
\end{equation}
Moreover, $\cc{T}_n\paren{f} = \brac{\paren{\cc{T}_{\bb{N}}\paren{f}}_{ij}}^n_{i,j=1}$ denotes a $n\times n$ truncated Toeplitz matrix generated by $f$. For two functions $f$ and $g$ on $\bb{R}$, we write $f\paren{t}\asymp g\paren{t}$ as $t\rightarrow t_0$, if $C_1\leq \lim\limits_{t\rightarrow t_0} \abs{\frac{f\paren{t}}{g\paren{t}}} \leq C_2$ for some strictly positive bounded scalars $C_1\leq C_2$. In particular, we write $f\paren{t}\sim g\paren{t}$ as $t\rightarrow t_0$ to indicate the case that $C_1=C_2=1$. Furthermore, for sequences $a_n$ and $b_n$, we write $b_n=\Omega\paren{a_n}$ when $b_n$ is bounded below by $a_n$ asymptotically, i.e. $\lim\limits_{n\rightarrow\infty} \abs{b_n/a_n}\geq C$ for some positive $C$. In case $a_n$ and $b_n$ are random, $b_n= o_{\bb{P}}\paren{a_n}$ means that $b_n/a_n$ converges in probability to zero as $n\rightarrow\infty$. For two sets  $\Omega_1,\Omega_2\subset\bb{R}^m$, $\dist\paren{\Omega_1,\Omega_2} \coloneqq \inf_{\omega_i\in\Omega_i,\;i=1,2} \LpNorm{\omega_1-\omega_2}{2}$ stands for their mutual distance with respect to Euclidean distance. Lastly, $\Gamma\paren{\cdot}$ denotes the gamma function. 

\section{Problem formulation and detection algorithms}\label{DtctAlg}

In this section we present a formulation of the shift-in-mean detection problem associated with Gaussian process data, and then describe detection algorithms that account for the underlying process modelling assumptions. Recall that $G$ is a Gaussian process in $\cc{D}\subset\bb{R}$ and $\set{X_k = G\paren{t_k}}^n_{k=1}$ represents the set of $n$ samples of $G$ at times $\cc{D}_n = \set{t_k}^n_{k=1}\subset\cc{D}$. At this point we proceed to split the problem formulation into two subsections according to the two different asymptotic settings. 

\subsection{Gaussian process model in fixed domain setting}

In the fixed domain setting, $G-\bb{E} G$ is assumed to be a mean-zero stationary Gaussian process on a bounded domain $\cc{D} = \brac{0,1}$ and regularly sampled at $t_k = k/n$ for $k=1,\ldots,n$. Let the symmetric real functions $K:\bb{R}\mapsto\bb{R}$ and $\hat{K}:\bb{R}\mapsto\bb{R}$ respectively denote the covariance function and spectral density of $G$. Accordingly, $\Sigma_n \coloneqq \cov\paren{\set{X_k}^n_{k=1}}$ is a symmetric Toeplitz matrix given by
\begin{equation}\label{SigmanFxdDom}
\Sigma_n = \set{\cov X_rX_s}^n_{r,s=1} = \brac{K\paren{\frac{r-s}{n}}}^{n}_{r,s=1}.
\end{equation}
In the next section, we impose some regularity conditions on $K$. 

\subsection{Gaussian process model in increasing domain setting}\label{IncDomFormulation}

In the increasing domain setting, $G-\bb{E} G$ is a mean-zero Gaussian process in $\cc{D} = \brapar{0,\infty}$ and so $\set{X_k}_{k}$ endowed with a Toeplitz covariance function. This is a common setting for time series data, where 
the observed samples are indexed by time points in $\cc{D}$. It is customary to assume that the covariance of the observed samples decreases as the temporal distance increases. Define $\cov\paren{X_1, X_k} = f_k$ for any $k$, in which $\set{f_m}^{\infty}_{m=0}$ be a absolutely summable sequence with $f_0 = 1$. Due to the stationarity assumption, $\Sigma_{\bb{N}}\coloneqq\cov\paren{\set{X_k}^{\infty}_{k=1}}$ is an infinite symmetric Toeplitz matrix. We view $\set{X_k}^n_{k=1}$ as the observed part of an infinite stationary time series, $\set{X_k}^{\infty}_{k=1}$. Accordingly, the covariance matrix of $\set{X_k}^n_{k=1}$, denoted by $\Sigma_n$, is a symmetric (truncated) Toeplitz matrix. 

It is a known fact that (Chapter 4, \cite{RM.Gray}) there is a symmetric and almost surely (with respect to \emph{Lebesgue} measure) positive function, $f:\brac{-\pi,\pi}\mapsto \bb{R}$ such that $\Sigma_{\bb{N}} = \cc{T}_{\bb{N}}\paren{f}$. Thus $\Sigma_n = \cc{T}_n\paren{f}$. For studying the asymptotic properties of the change detection algorithm, certain regularity conditions are required on $f$.

\subsection{Detection procedure based on generalized likelihood ratio test}\label{DtctProcedure}

Now we proceed to formulate the detection of the existence of a sudden change in the mean of a one dimensional Gaussian process as a composite hypothesis testing problem. As noted above, we are dealing with two different settings of the domains. That is, we assume that $G$ satisfies either of the two conditions:
\begin{enumerate}[label=(\alph*),leftmargin=*]
\item Fixed domain setting: $G$ is restricted to $\cc{D}=\brac{0,1}$ where its spectral density admits Assumption \ref{FxdDomAssu} and $\cc{D}_n = \set{k/n}^n_{k=1}$.
\item Increasing domain setting: The domain of $G$ is $\cc{D} = \bb{R}$ and the samples are taken at $\cc{D}_n = \set{1,\ldots,n}$. Moreover, $\Sigma_n = \cc{T}_n\paren{f}$ for some $f$ fulfilling Assumption \ref{AssumpToepCov}.
\end{enumerate}

Although the domain settings are different, the detection procedure that we propose based on a generalized likelihood ratio test will be the same. The composite hypothesis testing problem is set out as follows. Under the null hypothesis, all the random variables have zero mean, i.e. $\bb{E} \boldsymbol{X} = \zero_n$. To specify the alternative hypothesis $\bb{H}_1$ we first introduce a few additional notations. Let $t\in\cc{C}_{n,\alpha}$ denote the occurrence time of the single change point. The set $\cc{C}_{n,\alpha}\subseteq \set{1,\ldots,n}$ contains plausible occurrence time of the change, and we assume there is $\alpha\in\paren{0,1/2}$ such that $\cc{C}_{n,\alpha} = \set{t:\;t\wedge\paren{n-t} > \alpha n}$. Another important parameter $b$ denotes the amount of shift in the mean before and after the change point. Thus, for a fixed $t\in \cc{C}_{n}$, the associated alternative hypothesis to $t$ can be stated as,
\begin{equation}\label{H1t}
H_{1,t}:\;\exists\;b\ne 0,\;\bb{E} \boldsymbol{X} = \frac{b}{2}\zeta_t,
\end{equation}
where $\zeta_t\in\bb{R}^n$ is given by $\zeta_t\paren{k}\coloneqq \sign\paren{k-t}$ for any $t\in\cc{C}_{n,\alpha}$. Since $t$ is not known a priori, the alternative hypothesis is specified by taking the union of $\bb{H}_{1,t}$. 
Thus, the composite hypothesis testing problem is given by
\begin{equation}\label{DtctProb}
\bb{H}_0: \bb{E} \boldsymbol{X} = \zero_n,\quad v.s.\quad
\bb{H}_1 = \bigcup_{t\in\cc{C}_{n,\alpha}} \bb{H}_{1,t},\;\mbox{i.e.}\;\exists\;t\in\cc{C}_{n,\alpha},\; b\ne 0,\;s.t.\;\;\bb{E} \boldsymbol{X} = \frac{b}{2}\zeta_t.
\end{equation}
Next, we propose a test statistic which is constructed by the generalized likelihood ratio (GLR). Note that the GLR is an explicit function of the joint density of samples and so the Gaussian process assumption is essential to its calculation.

\begin{prop}\label{GLRTProp}
Assuming that $\Sigma_n$ is known, there exists $R_{n,\delta}>0$ for which the GLRT is given by
\begin{equation}\label{GLRT}
T_{GLRT} = \bb{I}\paren{\max_{t\in\cc{C}_{n,\alpha}} \abs{\frac{\zeta^\top_t\paren{\Sigma_n}^{-1}\boldsymbol{X}}{\sqrt{\zeta^\top_t\paren{\Sigma_n}^{-1}\zeta_t}}}^2\geq R_{n,\delta}}.
\end{equation}
\end{prop}

The threshold value $R_{n,\delta}$ depends only on $n$ and some parameter $\delta$ determining the false alarm rate. The precise form of $R_{n,\delta}$ will be presented in subsequent sections. We also note that setting $\mu = 0$ in \eqref{ProbForm} results in a substantially simplified expression of the GLR, which eases the exposition of our analysis of the computational and theoretical properties of the proposed test. The general form of the GLRT, when $\mu$ is unknown, is presented as Proposition \ref{GLRTGenForm} in Appendix \ref{Proofs}.

Unlike the CUSUM test, cf. Eq. \eqref{CUSUM}, the covariance function of $G$ is explicitly taken into account in the GLRT. As a result, it will be shown in the sequel that the proposed detection method is optimal, while the same cannot be said for the CUSUM test, specifically in the setting of fixed domain asymptotics. In practice, however, the covariance is not known and needs to be estimated. To address such scenarios, we propose to approximate the likelihood ratio by plugging in a positive definite estimate of the covariance matrix, which will be indicated by $\tilde{\Sigma}_n$, in Eq. \eqref{GLRT}. The plug-in detection technique will be called \emph{plug-in GLRT}. Here is the formulation of the plug-in GLRT, while the choice of $\tilde{\Sigma}_n$ and the accompanying theory will be given later in Section \ref{DtctRateUnKnwnSpDen}.

\begin{defn}\label{AGLRT}
Let $\tilde{\Sigma}_n$ be a positive definite estimate of $\Sigma_n$. The plug-in GLRT is given by 
\begin{equation}\label{PGLRT}
\tilde{T}_{GLRT} = \bb{I}\paren{\max_{t\in\cc{C}_{n,\alpha}} \abs{\frac{\zeta^\top_t\paren{\tilde{\Sigma}_n}^{-1}\boldsymbol{X}}{\sqrt{\zeta^\top_t\paren{\tilde{\Sigma}_n}^{-1}\zeta_t}}}^2\geq \tilde{R}_{n,\delta}},
\end{equation}
for some strictly positive threshold value $\tilde{R}_{n,\delta}$.
\end{defn}

\section{Detection rate of GLRT: known $\Sigma_n$}\label{DtctRateKnwnSpDen}

This section is devoted to a treatment of the detection rate of the GLRT, given that $\Sigma_n$ is known. Section \ref{DtctRateKnwnSpDenFxdDom} addresses this problem in the fixed domain, while Section \ref{DtctRateKnwnSpDenIncDom} considers the increasing domain. We first define the risk measure that we will consider in the subsequent sections. 

\begin{defn}\label{CDEP}
For any change detection algorithm $T\in\set{0,1}$, the conditional detection error probability (CDEP) of $T$ which is denoted by $\varphi_n\paren{T}$, is defined as
\begin{equation*}
\varphi_n\paren{T} = \bb{P}\paren{T=1\mid \bb{H}_0} + \max_{t\in\cc{C}_{n,\alpha}} \bb{P}\paren{T=0\mid \bb{H}_{1,t}}.
\end{equation*}
\end{defn}

\begin{rem}\label{RemCDEP}
In words, $\varphi_n$ is the sum of the false positive rate and the maximal miss detection rate taken over the set of possible change point locations $\cc{C}_{n,\alpha}$. Clearly, CDEP hinges on the choices of $\cc{C}_{n,\alpha}$ -- the value of $\varphi_n$ increases as $\cc{C}_{n,\alpha}$ becomes a larger proper subset of $\set{1,\ldots,n}$. CDEP as a risk measure has been considered in the literature for detecting abnormal clusters in a network (see e.g., \cite{E.Arias-Castro, C.Butucea}). It also provides an upper bound on the Bayesian risk measure. We refer the reader to \cite{L.Addario-Berry} for a prudent comparison of CDEP and the Bayesian risk measure. 
\end{rem}

In the following results, we seek for sufficient conditions on the shift value $b$ under which the CDEP of GLRT is bounded above by some 
small $\delta \in (0,1)$. 

\subsection{Detection rate in fixed domain setting}\label{DtctRateKnwnSpDenFxdDom}
The results in this sub-section are guaranteed for two common classes of covariance function, one of which admits polynomially decaying spectral density, and the other is Gaussian covariance function. Throughout this section we assume that $G$ is a one dimensional Gaussian process restricted to $\cc{D} = \brac{0,1}$ whose covariance function and spectral density are denoted by $K$ and $\hat{K}$, respectively. We first focus on the case of polynomially decaying $\hat{K}$.

\begin{assu}\label{FxdDomAssu}
$K$ is an integrable positive definite covariance function. Moreover, there exist $\nu\in\paren{0,\infty}$ and $C_K$ (depending on $K$) for which $\hat{K}$ satisfies the following condition:
\begin{equation}\label{BslSpcCond}
C_K\coloneqq\sup_{\omega\in\bb{R}}\abs{\hat{K}\paren{\omega}\paren{1+\omega^2}^{\nu+\frac{1}{2}}} < \infty.
\end{equation}
\end{assu}

We shall always choose the largest possible $\nu$ that satisfies \eqref{BslSpcCond}. It is simple to see that Assumption \ref{FxdDomAssu} holds if and only if $\hat{K}$ is bounded at the origin and $\hat{K}\paren{\omega}\asymp \omega^{-\paren{2\nu+1}}$ as $\omega$ tends to infinity. It is well-known that the tail behavior of $\hat{K}$ is closely linked to the smoothness of $K$ at the origin (e.g., Section $2.8$, \cite{ML.Stein}). The following are a few examples of common covariance functions that will be studied in this paper.
\begin{enumerate}[label=(\alph*),leftmargin=*]
\item \emph{Matern}: This class is widely used in geostatistics, and has a fairly simple explicit form of spectral density.
\begin{equation}\label{MaternSpDen}
\hat{K}\paren{\omega} = \frac{\sqrt{4\pi}\Gamma\paren{\nu+1/2}}{\Gamma\paren{\nu}}\sigma^2\rho^{-2\nu}\paren{\frac{1}{\rho^2}+\omega^2}^{-\paren{\nu+1/2}},
\end{equation}
in which $\rho,\nu,\sigma\in\paren{0,\infty}$. Regardless of the choice of $\rho$ and $\sigma$, condition \eqref{BslSpcCond} holds for Matern spectral density with parameter $\nu$.
\item \emph{Powered exponential}: Another versatile class of covariance functions is 
\begin{equation}\label{PowExpCov}
K\paren{r} =\sigma^2 \exp\paren{-\abs{\frac{r}{\rho}}^{\beta}}
\end{equation}
for some $\beta\in\paren{0,2}$ and $\rho,\sigma\in\paren{0,\infty}$. Although the spectral density does not have a closed form in terms of simple functions, Lemma \ref{PowExpSpDenDecay} shows that $\hat{K}$ admits Assumption \ref{FxdDomAssu} with $\nu = \beta/2$.
\item \emph{Rational spectral densities}: Rational spectral densities form a general class admitting Assumption \ref{BslSpcCond}. For any $\hat{K}$ in this class, there are two polynomials, $Q_n$ and $Q_d$, with real coefficients, unit leading coefficients and $p\coloneqq\deg\paren{Q_d}-\deg\paren{Q_n}\in\bb{N}$, such that
\begin{equation}\label{RationalSpDen}
\hat{K}\paren{\omega} = \lambda\frac{\abs{Q_n\paren{j\omega}}^2}{\abs{Q_d\paren{j\omega}}^2}.
\end{equation}
Moreover, we assume that $Q_d$ has no root on the imaginary axis and $\lambda$ is a strictly positive scalar. Since $K\paren{0}<\infty$ and $\hat{K}\paren{\omega} \asymp \omega^{-2p}$ as $\omega\rightarrow\infty$, 
Assumption \ref{FxdDomAssu} holds with $\nu = p-1/2$. 
\item \emph{Triangular}: For $T,\sigma\in\paren{0,\infty}$, the covariance function and spectral density are given by 
\begin{equation*}
K\paren{r} = \sigma^2\paren{1-\abs{\frac{r}{\rho}}}_{+},\quad \hat{K}\paren{\omega} = \frac{\rho\sigma^2}{2}\abs{\sinc\paren{\frac{\rho\omega}{2}}}^2.
\end{equation*}
Triangular covariance is less favorable than the aforestated cases due to the oscillatory behaviour of $\hat{K}$ (p. $31$, \cite{ML.Stein}). One can easily show that this covariance fulfils Assumption \ref{FxdDomAssu} with $\nu = 1/2$.
\end{enumerate}

\begin{thm}\label{GLRTRateFxdDom}
Let $\delta\in\paren{0,1}$. Suppose that $G$ is a one dimensional Gaussian process restricted to $\cc{D} = \brac{0,1}$ whose associated spectral density $\hat{K}$ admits Assumption \ref{FxdDomAssu} for some $\nu$ and $C_K$. $G$ is regularly sampled on $i/n,\;i=1,\ldots,n$. There exist $R_{n,\delta}>0$ (depending only on $n$ and $\delta$), $n_0\coloneqq n_0\paren{K}$ and a positive universal constant $C$ such that if $n\geq n_0$ and 
\begin{equation}\label{LowBndDtctbaleJmp}
\abs{b}\geq Cn^{-\nu} \sqrt{C_{K}\paren{1+\frac{1}{\nu}}\log\paren{\frac{n\paren{1-2\alpha}}{\delta}}},
\end{equation}
we have
\begin{equation*}
\varphi_n\paren{T_{GLRT}} \leq \delta.
\end{equation*}
\end{thm}

See Appendix \ref{ProofThm1} for the proof of Theorem \ref{GLRTRateFxdDom}. We now make several comments regarding the roles of various quantities embedded in Theorem \ref{GLRTRateFxdDom}.

\begin{enumerate}[label=(\alph*),leftmargin=*]
\item $R_{n,\delta}$ in Theorem \ref{GLRTRateFxdDom} can be chosen as
\begin{equation}\label{CritcValGLRTEq}
R_{n,\delta} = 1+2\brac{\log\paren{\frac{2n\paren{1-2\alpha}}{\delta}} + \sqrt{\log\paren{\frac{2n\paren{1-2\alpha}}{\delta}}}}.
\end{equation}
We guarantee that CDEP is less than or equal $\delta$ by controlling false alarm and miss detection probabilities below $\delta/2$. Our trick provides an evidence to choose $R_{n,\delta}$. Notice that under null hypothesis, the test statistic in Eq. \eqref{GLRT} has the same distribution as the supremum of a $\chi^2_1$ process over $\cc{C}_{n,\alpha}$, which is represented by $\set{\Psi\paren{t}:\;t\in\cc{C}_{n,\alpha} }$. Strictly speaking for controlling the false alarm probability below $\delta/2$, $R_{n,\delta}$ needs to be chosen in such a way that
\begin{equation*}
\bb{P}\paren{\max_{t\in\cc{C}_{n,\alpha}} \Psi\paren{t} \geq R_{n,\delta}}\leq \frac{\delta}{2}.
\end{equation*}
The standard $\chi^2_1$ tail inequality in \cite{L.Birge} implies that if  $R_{n,\delta}$ is chosen based upon Eq. \eqref{CritcValGLRTEq}, then $\Psi\paren{t}$ is below $\delta/\set{2n\paren{1-2\alpha}}$ for any $t\in\cc{C}_{n,\alpha}$. Thus, the union bound inequality yields
\begin{equation*}
\bb{P}\paren{\sup_{t\in\cc{C}_{n,\alpha}} \Psi\paren{t} \geq R_{n,\delta}} \leq \abs{\cc{C}_{n,\alpha}}\max_{t\in \cc{C}_{n,\alpha} } \bb{P}\paren{ \Psi\paren{t} \geq R_{n,\delta}} \leq \frac{\delta\abs{\cc{C}_{n,\alpha}}}{2n\paren{1-2\alpha}} = \frac{\delta}{2}.
\end{equation*}
\item The minimal detectable shift is proportional to $\sqrt{C_K}$, as defined in \eqref{BslSpcCond}. Note that $C_K$ is determined by both low frequency and tail behaviour of spectral density via $\nu$.  (E.g., for Matern covariance functions given by \eqref{MaternSpDen}, $C_K = \frac{\sqrt{4\pi}\Gamma\paren{\nu+1/2}}{\Gamma\paren{\nu}}\sigma^2\paren{1 \vee \rho^{-2\nu}}$). It is easily verifiable by \eqref{BslSpcCond} that $C_K$ is linearly proportional to $\sqrt{K\paren{0}}$, meaning that $C_K$ also captures the notion of the standard deviation of the observations. Thus, Theorem \ref{GLRTRateFxdDom} implicitly expresses that change detection is more challenging for Gaussian processes with larger variance.
\item Sample size $n$ has two opposing effects on the detection rate. On the one hand, $n\paren{1-2\alpha}$ appearing in the logarithmic function, is closely connected to the size of alternative hypothesis which is determined by $\abs{\cc{C}_{n,\alpha}} = n\paren{1-2\alpha}$. On the other hand, the term $n^{-\nu}$ indicates the possibility of small shift detection as more observations are available. 
\end{enumerate}

We will see in Section \ref{DtctRateKnwnSpDenIncDom} that parameters $\delta$, variance of observations and sample size have almost analogous roles in the increasing domain change detection. The main difference between the two asymptotic settings is the role of the decay rate of $\hat{K}$ in fixed domain, which is encapsulated as $\nu$. Note that $\nu$ is closely related to the smoothness of $G$ with larger values of $\nu$ corresponds to smoother Gaussian process in the mean squared sense (cf. \cite{ML.Stein}, Chapter $2$). For smooth Gaussian processes, $G\paren{t_0}$ can be interpolated using the observations in the vicinity of $t_0$ with small estimation error. This leads to a simpler shift-in-mean detection for smoother processes. More precisely, as $n\rightarrow\infty$ the lower bound on detectable $b$, \eqref{LowBndDtctbaleJmp}, vanishes more rapidly for larger $\nu$.

\begin{rem}\label{Rem3.2}
We describe the rate of the minimal detectable jump for some specific commonly used classes of spectral densities, all of which satisfy Assumption \ref{FxdDomAssu}.
\begin{enumerate}[label=(\alph*),leftmargin=*]
\item \emph{Matern}: For the Matern class with parameters $\paren{\sigma,\rho,\nu}$ satisfying Assumption \ref{FxdDomAssu},
the smallest detectable shift in mean is $\abs{b} \asymp n^{-\nu}\sqrt{\log \paren{n\paren{1-2\alpha}/\delta}}$. 
\item \emph{Powered exponential}: For the power exponential class, any jump size of magnitude at least
$\abs{b} = \Omega\paren{\sqrt{n^{-\beta}\log \paren{n\paren{1-2\alpha}/\delta}}}$ is detectable. On the contrary to Matern class, obtaining a closed form for $C_K$ is quite difficult for powered exponential class.
\item \emph{Rational spectral densities}: %
It has been discussed previously 
that $\hat{K}\paren{\omega}\asymp \abs{\omega}^{-2p}$ as $\omega\rightarrow\infty$ and $\nu = p-1/2$, revealing that each $\abs{b} = \Omega\paren{n^{-\paren{p-1/2}}\sqrt{\log \paren{n\paren{1-2\alpha}/\delta}}}$ can be detected with high probability.
\item \emph{Triangular}: Since $\hat{K}$ satisfies Assumption \ref{FxdDomAssu} with $\nu = 1/2$, so any $\abs{b} = \Omega\paren{\sqrt{n^{-1}\log \paren{n\paren{1-2\alpha}/\delta}}}$ is detectable. Although it needs a great algebraic effort to find $C_K$, it can be shown easily that $C_K\leq \sigma\paren{2/\rho+\rho/2}$.
\end{enumerate}
\end{rem}

We conclude this section by a comprehensive explanation the role of $\alpha$ in Theorem \ref{GLRTRateFxdDom}. The dependence on $\alpha$ in Eq. \eqref{LowBndDtctbaleJmp} is logarithmic, which encoding how the size of $\cc{C}_{n,\alpha}$ affecting the detection rate.

\begin{rem}\label{RemRoleofAlpha}
The minor role of $\alpha$ in Eq. \eqref{LowBndDtctbaleJmp} may seem a bit surprising. Strictly speaking the asymptotic behavior of the smallest detectable jump remain the same, regardless of how small $\alpha$ has been chosen (even if $\alpha$ tends ro zero). It is also notable to mention that we did not use the assumption that $\alpha$ is a fixed and strictly positive scalar in our proof. This puzzling aspect can be resolved by a deeper look at the formulation of the hypothesis testing problem \eqref{DtctProb}. For algebraic convenience, we assume that the mean of $G$ fluctuates around $\mu = 0$ in Eq. \eqref{DtctProb}. The fact that we assumed $\mu$ is known is the main reason that $\alpha$ parameter in the detraction rate of GLRT, as we do not need to estimate $\mu$ from the data. That is why in this particular case $\alpha$ even be chosen as small as $\cc{O}\paren{1/n}$. We want to emphasize that the generic form of GLRT test for unknown $\mu$ are presented in Proposition \ref{GLRTGenForm}. We believe that in the analysis of the extended version of GLRT, the constant $C$ in Eq. \eqref{LowBndDtctbaleJmp} depends on $\alpha$ (without changing the dependence on $n$ and other parameters). 
\end{rem}

\paragraph{Gaussian covariance function.}\label{DtctRateKnwnSpDenFxdDomGaussCovFunc}

The Gaussian covariance function is given by
\begin{equation}\label{GaussCovFunc}
K\paren{r} = \sigma^2\exp\brac{-\frac{1}{2}\paren{\frac{r}{\rho}}^2},\quad \hat{K}\paren{\omega} = \rho\sigma^2\sqrt{2\pi}\exp\brac{-\frac{\paren{\rho\omega}^2}{2}}.
\end{equation}
It is widely used in practice for modeling of smooth Gaussian processes, e.g. in \cite{WL.Loh}. Regarding this choice of covariance, we have the following result:

\begin{thm}\label{GLRTRateFxdDomGaussCov}
Let $G$ be a Gaussian process on $\brac{0,1}$ which is observed at $i/n,\;i=1,\ldots,n$, whose covariance function is given by Eq. \eqref{GaussCovFunc}. Choose $\delta\in\paren{0,1}$. There are $R_{n,\delta}>0$, $n_0\coloneqq n_0\paren{\rho}$, $C_0\coloneqq C_0\paren{\rho}>0$ and a universal constant $C>0$, such that if $n\geq n_0$ and
\begin{equation}\label{GaussCovFxdDtctRate}
\abs{b}\geq C\sqrt{\exp\Bigbrac{-n\log\paren{C_0n}}\log\paren{\frac{n\paren{1-2\alpha}}{\delta}}},
\end{equation}
then
\begin{equation*}
\varphi_n\paren{T_{GLRT}} \leq \delta.
\end{equation*}
\end{thm}

The details of the proof is given in Appendix \ref{ProofThm2}. Because of the super-exponential decay of Gaussian spectral density, Assumption \ref{FxdDomAssu} is actually satisfied for any $\nu>0$. This result shows that it is possible to detect exponentially small jump size $b$ as $n$ increases. Moreover, the result of Theorem \ref{GLRTRateFxdDomGaussCov} is also compatible with that of Theorem \ref{DtctRateKnwnSpDenFxdDom}, in quantifying precisely the assertion that the smoother the Gaussian process is, the easier it is to detect the presence of the shift in the mean.

\subsection{Detection rate in increasing domain setting}\label{DtctRateKnwnSpDenIncDom}

Turning now to the increasing domain setting, recall that $G$ is assumed to follow the setting described in Section \ref{IncDomFormulation}. We also assume that $\Sigma_n = \cc{T}_n\paren{f}$ for some function $f$ satisfying the following conditions.

\begin{assu}\label{AssumpToepCov}
$f:\brac{-\pi,\pi}\mapsto \bb{R}$ is a real symmetric function such that
\begin{enumerate}[label=(\alph*),leftmargin=*]
\item There are two positive universal scalars, $0 < m_f \leq M_f < \infty$ such that
\begin{equation*}
m_f \coloneqq \inf_{\omega\in\brac{-\pi,\pi}} f\paren{\omega} \leq M_f \coloneqq \sup_{\omega\in\brac{-\pi,\pi}} f\paren{\omega}.
\end{equation*}
\item There exist positive universal constants $c$ and $\lambda$ such that
\begin{equation}\label{UppBndfk}
\abs{f_k} \leq c \paren{1+k}^{-\paren{1+\lambda}}.
\end{equation}
\end{enumerate}
\end{assu}

Note that the first condition regarding the infimum of $f$ is necessary to have a positive definite infinite covariance matrix, i.e., $\nu^\top \Sigma_{\bb{N}}\nu > 0$ for any non-zero $\nu\in\bb{R}^{\bb{N}}$. Moreover, the polynomial decay of $f_k$'s as stated in \eqref{UppBndfk} is a sufficient condition to ensure that $f$ can be equivalently expressed by its Fourier series. Such condition is common in the non-asymptotic analysis of Toeplitz matrices (see, e.g., \cite{RM.Gray}).

\begin{thm}\label{GLRTRateIncDom}
Let $\delta\in\paren{0,1}$ and suppose that $\Sigma_n = \cc{T}_n\paren{f}$ in which $f$ admits Assumption \ref{AssumpToepCov} for some positive scalars $c$ and $\lambda$. There exist $n_0\in\bb{N}$, $C>0$ (depending only on $c$ and $\lambda$) and $R_{n,\delta}>0$ such that for any $n\geq n_0$, if
\begin{equation}\label{IncDomDtctRate}
\abs{b}\geq C\sqrt{f\paren{0}n^{-1}\log\paren{\frac{n\paren{1-2\alpha}}{\delta}}},
\end{equation}
then
\begin{equation*}
\varphi_n\paren{T_{GLRT}} \leq \delta.
\end{equation*}
\end{thm}

See Appendix \ref{ProofThm3} for the proof of Theorem \ref{GLRTRateIncDom}. Some comments are in order. First, the threshold $R_{n,\delta}$ in Theorem \ref{GLRTRateIncDom} is chosen in exactly the same way as in the fixed domain setting, as given by Eq. \eqref{CritcValGLRTEq}. Second, in contrast to the fixed domain setting, the dependence structure for $G$ no longer plays the central role in the characterization of the detection performance. In particular, $f\paren{0}$ is the only factor in \eqref{IncDomDtctRate} that captures the correlation in the samples, but this scalar quantity evidently has an insignificant effect: the asymptotic behaviour of GLRT remains the same (up to some constant factor) for different Gaussian processes satisfying Assumption \ref{AssumpToepCov}. A related observation that arises by comparing between \eqref{LowBndDtctbaleJmp} and \eqref{IncDomDtctRate} is that the correlation structure of observations, which is encapsulated into $\nu$ or $f\paren{0}$, and the quantities encoding the marginal density information such as $n$ have been completely decoupled in the rate of GLRT in increasing domain. An examination of the proof reveals that the decoupling effect in the increasing domain setting arises due to the short-range correlation assumption ($\cov\paren{X_r,X_s}\rightarrow 0$ polynomially in $\abs{r-s}$). It follows that as $n$ increases the correlation for most pairs of observed sample become negligible.

\begin{rem}
Although for algebraic convenience, throughout the paper we assume that $\alpha$ is a fixed and positive constant, it is not necessary for our proof technique of Theorem \ref{GLRTRateIncDom}. We basically can generalize Theorem \ref{GLRTRateIncDom} to the case that $\alpha$ tends to zero as $n\rightarrow\infty$. A deeper look at our proof (see \ref{ProofThm3}) and the auxiliary results in Appendix \ref{NonAsympInvLargToep} reveals that Eq. \eqref{IncDomDtctRate} can be replaced by a complicated from
\begin{equation*}
\abs{b}\geq C\sqrt{\log\paren{\frac{n\paren{1-2\alpha}}{\delta}}}\sqrt{f\paren{0}- \xi_n},
\end{equation*}
In which $\xi_n$ is a vanishing sequence of $n$, which can be easily eliminated from the detection rate by adjusting the universal constant $C$. We refer the reader to Remark \ref{RemRoleofAlpha} for a detailed discussion on the role of $\alpha$ in the case that $\mu$ is unknown in the formulation of the shift in mean problem (\eqref{ProbForm}).
\end{rem}

\section{Detection rate of plug-in GLRT}\label{DtctRateUnKnwnSpDen}

As we have shown in Proposition \ref{GLRTProp}, full knowledge of $\Sigma_n$ is central to computing the generalized likelihood ratio. In practice, the spectral density and covariance function of $G$ are not known a priori, and so we take a plug-in approach, by approximating the GLRT by estimating the covariance estimate $\Sigma_n$ (see Definition \ref{AGLRT}). This section serves to investigate various ways of constructing plug-in GLRT and assessing its detection performance. We focus only on the fixed domain setting, because the dependence structure underlying the hypothesis plays an important role in determining the detection error rate, as shown in Section \ref{DtctRateKnwnSpDenIncDom}.

We first assume that $G$ is a Matern Gaussian process on $\cc{D = \brac{0,1}}$ with unknown parameters $\eta = \paren{\sigma,\rho}\in\Omega$ (see Eq. \eqref{MaternSpDen}), and is regularly observed on $\set{i/n}^n_{i=1}$. We use $\tilde{\eta}_m = \paren{\tilde{\sigma}_m,\tilde{\rho}_m}$ to indicate the estimated parameters using $m$ regularly spaced samples in $\cc{D}$. We also assume that $\cc{C}_{n,\alpha}\subseteq\set{k:\; \alpha n\leq k\leq \paren{1-\alpha}n}$. Namely, the Gaussian process is under control for a certain number of observations. The controlled samples before the sudden change, $\boldsymbol{X}_{\cc{B}}\coloneqq\set{X_k:\; k\leq \alpha n}$, will be used to estimate $\eta$. The parameter estimation stage is typically called the \emph{burn-in} period in the literature.

It is known in the Gaussian processes literature that $\eta$ is not consistently estimable in the fixed domain setting when the number of the observations in $\cc{D}$ grows to infinity --- see 
\cite{Z.Ying,H.Zhang} for further details. Zhang \cite{H.Zhang} showed that neither $\sigma$ or $\rho$ are consistently estimable but the quantity $\sigma \rho^{-\nu}$ can be consistently estimated using MLE. The profound reason behind the inconsistency is the existence of a class of mutually absolutely continuous models for $G$ which are almost surely impossible to discern by observing one realization of $G$. Strictly speaking, the induced measures corresponding to two Matern Gaussian processes with parameters $\eta$ and $\eta'$ are absolutely continuous with respect to each other, whenever $\sigma\rho^{-\nu} = \sigma'\rho'^{-\nu}$. Furthermore Zhang \cite{H.Zhang}  showed that if one fixes $\rho$ at an arbitrary value, then the maximum likelihood estimator for $\sigma\rho^{-\nu}$ is consistent. We shall show that despite the inconsistency in estimating $\eta$, plug-in GLRT exhibits an analogous performance as GLRT with fully known covariance function whenever the estimate of $\eta$ is consistent up to the equivalence class.

It has been discussed in that fixing $\tilde{\rho}_m$ at large values has trifling impact on predictive performance. Note that due to the complicated dependence of Matern covariance function to $\rho$, estimating $\rho$ is a computationally challenging task, particularly for large data sets. So we can accelerate the whole detection procedure without estimating $\rho$. Plug-in GLRT change detector is a two stage algorithm as follows:
\begin{itemize}[leftmargin=*]
\item\emph{Estimation step:}
\begin{enumerate}
\item Fix $\tilde{\rho}_m$ at the largest possible element in $\Omega$. Namely, $\tilde{\rho}_m$ is a deterministic quantity given by $\tilde{\rho}_m = \sup\set{\rho:\paren{\sigma,\rho}\in\Omega}$.
\item Estimate $\sigma\rho^{-\nu}$ given the controlled samples $\boldsymbol{X}_{\cc{B}}$, using any consistent procedure such as maximum likelihood (MLE) \cite{H.Zhang}, weighted local Whittle likelihood \cite{WY.Wu}, averaging quadratic variation \cite{E.Anderes}. We use the term \emph{consistent} to refer to the cases that $\abs{ \sigma\rho^{-\nu} - \tilde{\sigma}_m\tilde{\rho}^{-\nu}_m} \cp{\bb{P}} 0$ as $m$ grows to infinity.
\item Construct the approximated covariance matrix of $\boldsymbol{X}$, as $\tilde{\Sigma}_{n} = \brac{K\paren{\frac{r-s}{n},\tilde{\eta}_m}}^n_{r,s=1}$ (here $m\coloneqq \floor{\alpha n}$).
\end{enumerate}
\item \emph{Detection step:}
\begin{enumerate}
\item Applying the GLRT by plugging $\tilde{\Sigma}_{n}$ in place of $\Sigma_n$ into \eqref{GLRT}, as described in Definition \ref{AGLRT}. 
\end{enumerate}
\end{itemize}

Now we turn to state the main result of this section regarding the rate of the plug-in GLRT. 

\begin{thm}\label{PlugInGLRTRateFxdDom}
Let $\delta\in\paren{0,1}$. Let $G$ be Gaussian process whose associated spectral density $\hat{K}$ has Matern form with unknown parameters $\paren{\sigma,\rho}\in\Omega$. Given regular samples of one realization of $G$, there are finite scalar $C,n_0\in\bb{N}$, a non-negative sequence $\lim\limits_{m\rightarrow\infty}\tau_m = 0$, and threshold level $R_{n,\delta}>0$ such that for any $n\geq n_0$,
\begin{equation*}
\varphi_n\paren{\tilde{T}_{GLRT}} \leq \delta+2\tau_m,
\end{equation*}
whenever

\begin{equation}\label{ALGRTRate}
\abs{b}\geq Cn^{-\nu}\sqrt{C_{K}\paren{1+\frac{1}{\nu}} \log\paren{\frac{n\paren{1-2\alpha}}{\delta}}}.
\end{equation}
\end{thm}

See Appendix \ref{ProofThm4} for the proof of Theorem \ref{PlugInGLRTRateFxdDom}.

\begin{rem}\label{CritcValUAGLRT}
The threshold value for the plug-in GLRT is chosen exactly same as in Theorem \ref{GLRTRateFxdDom}.
\begin{equation}\label{CriticValPGLRT}
R_{n,\delta} = \brac{1+2\paren{\log\paren{\frac{2n\paren{1-2\alpha}}{\delta}} + \sqrt{\log\paren{\frac{2n\paren{1-2\alpha}}{\delta}}}}}.
\end{equation}
\end{rem}

Since $\set{\tau_k}_{k\in\bb{N}}$ is a vanishing sequence and $m = \floor{\alpha n}$ is an increasing function of $n$, $\paren{\delta+2\tau_m}$ lies in the vicinity of $\delta$ for large $n$. It is worthwhile to mention that $C$ appearing in Theorem \ref{PlugInGLRTRateFxdDom} is larger than the previously introduced scalar in Theorem \ref{GLRTRateFxdDom}, which can be viewed as the cost of mis-specifying $\rho$. The most interesting aspect of
Theorem \ref{PlugInGLRTRateFxdDom} is perhaps that if some consistent estimate of $\sigma\rho^{-\nu}$ is available then regardless of its rate, the plug-in GLRT has asymptotically the same rate as the GLRT with fully known covariance function.

We conclude this section by studying the performance of plug-in GLRT when both variance and range parameter are consistently estimable. In this case the plug-in GLRT test is constructed by replacing the estimated parameters in the GLRT test statistic. Suppose that $G$ has a powered exponential covariance function, introduced in Eq. \eqref{PowExpCov}. Anderes \cite{E.Anderes} proposed a consistent estimate of covariance parameters using empirical average of the quadratic variation of $G$. According to Theorem $5$ of \cite{E.Anderes}, unlike the Matern class, both $\sigma_0$ and $\rho_0$ are consistently estimable when $\beta\in\paren{0,1/2}$. Namely, $\abs{\rho-\tilde{\rho}_m} \vee \abs{\sigma-\tilde{\sigma}_m}\cp{\bb{P}} 0$, for the method introduced in \cite{E.Anderes}. The following result, which has a similar flavor as Theorem \ref{PlugInGLRTRateFxdDom}, determines the detection rate of plug-in GLRT for  one dimensional powered exponential Gaussian processes.

\begin{thm}\label{PlugInGLRTRateFxdDom2}
Let $\delta\in\paren{0,1}$. Let $G$ be Gaussian with powered exponential covariance function with unknown parameters $\paren{\sigma,\rho}\in\Omega$ and known $\beta\in\paren{0,1/2}$. Given regular samples of one realization of $G$, there are finite scalar $n_0\in\bb{N}$ and $C$ (which depends on the covariance parameters $\beta,\sigma$ and $\rho$), a non-negative sequence $\lim\limits_{m\rightarrow\infty}\tau_m = 0$, such that for any $n\geq n_0$,
\begin{equation*}
\varphi_n\paren{\tilde{T}_{GLRT}} \leq \delta+2\tau_m,
\end{equation*}
whenever
	
\begin{equation*}
\abs{b}\geq C\sqrt{n^{-\beta} \log\paren{\frac{n\paren{1-2\alpha}}{\delta}}}.
\end{equation*}
and 
\begin{equation*}
R_{n,\delta} = \brac{1+2\paren{\log\paren{\frac{2n\paren{1-2\alpha}}{\delta}} + \sqrt{\log\paren{\frac{2n\paren{1-2\alpha}}{\delta}}}}}.
\end{equation*}
\end{thm}

Theorem \ref{PlugInGLRTRateFxdDom2} states that given a consistent estimate of $\eta=\paren{\sigma_0,\rho_0}$, the plug-in GLRT procedure has the same asymptotic behavior as GLRT with fully known parameters (see part $\paren{b}$ of Remark \ref{Rem3.2} for the detection rate of GLRT with known $\sigma_0$ and $\rho_0$).

\section{Detection rate of CUSUM}\label{DtctCUSUM}

In this section we revisit the classical CUSUM test and present several results regarding its detection rate in both asymptotic settings. These results should be contrasted with our earlier theorems on the performance of the proposed exact and plug-in GLRT tests, and highlight the need for accounting for the dependence structures underlying the data, especially in the fixed domain setting. In the following, Theorem \ref{CUSUMFxdDom} introduces sufficient condition on $\abs{b}$ under which CUSUM can distinguish null and alternative hypotheses with high probability. Theorem \ref{CUSUMIncDom} studies the performance of CUSUM in the increasing domain setting. 

\begin{thm}\label{CUSUMFxdDom}
Let $G$ be a Gaussian process in $\brac{0,1}$ satisfying $\norm{K}{1}<\infty$ and $\norm{\hat{K}'}{\infty}<\infty$. Moreover let $\delta\in\paren{0,1}$, $\alpha\in\paren{0,1/2}$ and $\cc{C}_{n,\alpha} = \brac{\alpha n,\paren{1-\alpha}n}$. Given $n$ samples of one realization of $G$ at $i/n,\;i=1,\ldots,n$, there are $R_{n,\delta}>0$, and $n_0\coloneqq n_0\paren{\delta,\alpha}$ such that if $n\geq n_0$ and
\begin{equation}\label{SuffCondCUSUMFxdDom}
\abs{b}\geq 4\sqrt{\frac{\log\paren{\frac{2n\paren{1-2\alpha}}{\delta}}}{\alpha\paren{1-\alpha}}},
\end{equation}
then,
\begin{equation*}
\varphi_n\paren{T_{CUSUM}} \leq \delta.
\end{equation*}
\end{thm}

We refer the reader to Appendix \ref{PlugInGLRTRateFxdDom} for the proof of above result. The risk of fixed domain-CUSUM has been controlled from above under mild conditions on $K$, which holds true for all considered examples of covariance functions in this paper. It is indeed obvious form the following inequality that $K$ satisfies the assumptions in Theorem \ref{CUSUMFxdDom} if $a\paren{r}\coloneqq rK\paren{r}$ is absolutely integrable:
\begin{equation*}
\norm{\hat{K}'}{\infty} = \sup_{\omega\in\bb{R}} \abs{\int\limits_{-\infty}^{\infty} a\paren{r}e^{-j\omega r}dr}\leq \int\limits_{-\infty}^{\infty}\abs{rK\paren{r}}dr.
\end{equation*}
The main feature of the above theorem is the sufficient condition that the jump size increases (at the order of $\log n$ at least) in order to have an upper bound guarantee on the detection error. Although we do not have a proof that this sufficient condition is also necessary, our result suggests that the CUSUM test is \emph{inconsistent} in the fixed domain setting: the detection error may not vanish as data sample size increases, when the jump size is a constant. This statement is 
in fact verified by simulations. By contrast, we have shown earlier that using the GLRT based approach, we can guarantee the detection error to vanish as long as the jump size is either constant or (better yet) 
bounded from below by a suitable vanishing term.

\begin{rem}\label{Rem5.1}
Let us give a qualitative argument for the inconsistency of the CUSUM test in the fixed domain setting. Suppose that $b$ tends to zero as $n\rightarrow\infty$. Define
\begin{equation*}
U_t \coloneqq \sqrt{\frac{t\paren{n-t}}{n}}\brac{\frac{1}{n-t}\sum\limits_{k=t+1}^{n} X_k - \frac{1}{t}\sum\limits_{k=1}^{t} X_k}.
\end{equation*}
The expected value of $U_t$ is zero, under the null hypothesis and for any $t$. Regardless of the existence of a shift in the mean, the standard deviation of $U_t$ remains the same. A careful look at the proof of Theorem \ref{CUSUMFxdDom} reveals that the smallest value of the standard deviation of $U_t$ over $t\in\cc{C}_{n,\alpha}$ is order $\sqrt{n}$. Moreover, if there is a shift in the mean occuring at the change point $\bar{t}\in\cc{C}_{n,\alpha}$, then the expected value of $U_{\bar{t}}$ is given by $ b\sqrt{\bar{t}\paren{n-\bar{t}}/n} = \cc{O}\paren{b\sqrt{n}}$ (Recall that $\alpha n\leq \bar{t}\leq \paren{1-\alpha}n$). Generally speaking, as the mean of $U_t$ under the null hypothesis, denoted by $\bb{E}\paren{U_t\mid \bb{H}_0}$, is zero for any $t\in\cc{C}_{n,\alpha}$, the CUSUM test cannot distinguish between the null and the alternative (even for large sample size), since
\begin{equation*}
\abs{\frac{\bb{E}\paren{U_{\bar{t}}\mid \bb{H}_1}}{\sqrt{\var\paren{U_{\bar{t}}}}}} = \cc{O}\paren{\frac{b\sqrt{n}}{\sqrt{n}}} = \cc{O}\paren{b}\rightarrow 0,\quad\mbox{as}\; n\nearrow\infty.
\end{equation*}
Here, $\bb{E}\paren{U_{\bar{t}}\mid \bb{H}_1}$ represents the expected value of $U_t$ under the alternative. This suggests that regardless of the sample size, the CUSUM test cannot detect the existence of a small shift in the mean in the fixed domain setting.
\end{rem}

\begin{rem}
The threshold value of the CUSUM test in Theorem \ref{CUSUMFxdDom} is given by 
\begin{equation*}
R_{n,\delta} = \sqrt{n\paren{1+2\log\paren{\frac{2n\paren{1-2\alpha}}{\delta}}+ 2\sqrt{\log\paren{\frac{2n\paren{1-2\alpha}}{\delta}}} }}.
\end{equation*}
This threshold has different form of dependence to $n$ than that of the threshold of GLRT in Eq. \eqref{CritcValGLRTEq}, since unlike GLRT the CUSUM test do not reduce the correlation among the samples. In order to remove the gap between the threshold of GLRT and CUSUM in the fixed domain setting, we further normalize  $U_t$ by considering $U^{\star}_n = U_n/\sqrt{n}$. So equivalently CUSUM test in this regime can be written as  
\begin{equation*}
T_{CUSUM} = \bbM{1}\paren{ \max_{t\in\cc{C}_{n,\alpha}} \abs{U^{\star}_n}^2> R^{\star}_{n,\delta} \coloneqq \frac{R^2_{n,\delta}}{n} }.
\end{equation*}
Here $R^{\star}_{n,\delta}$ is exactly same as the critical value of GLRT test.
\end{rem}

Now we aim to study CUSUM test in the increasing domain asymptotic regime. Recall that for this scenario, $G$ is a Gaussian process with covariance matrix $\Sigma_n = \cc{T}_n\paren{f}$

\begin{thm}\label{CUSUMIncDom}
Let $\delta\in\paren{0,1}$, $\vartheta > 0$ and $\cc{C}_{n,\alpha} = \brac{\alpha n,\paren{1-\alpha}n}$. Assume that $f$ satisfies Assumption \ref{AssumpToepCov} for some $c$ and $\lambda$. There are $n_0 = n_0\paren{f,\vartheta}$ and universal constant $C\paren{\lambda,c}>0$, such that if $n\geq n_0$ and
\begin{equation}\label{DtctRateCUSUMIncDom}
\abs{b}\geq C\sqrt{\frac{\paren{1+\vartheta}f\paren{0}}{n\alpha\paren{1-\alpha}}\log\paren{\frac{n\paren{1-2\alpha}}{\delta}}},
\end{equation}
then 
\begin{equation*}
\varphi_n\paren{T_{CUSUM}} \leq \delta.
\end{equation*}
\end{thm}

Th reader can find the detailed proof of Theorem \ref{CUSUMIncDom} in Appendix \ref{ProofThm6}. By comparing between Theorems \ref{GLRTRateIncDom} and \ref{CUSUMIncDom}, it is clear that the CUSUM test exhibits a similar detection performance as the GLRT test in the increasing domain setting. In fact, we will show in the next section that both tests achieve minimax optimality in that setting. However, according to the numerical studies, using GLRT slightly improve the detection performance comparing to CUSUM, especially in the presence of strong long range dependence. 
Thus, in practice, one can afford to ignore the covariance structure of $G$ in the increasing domain setting, and due to its simplicity, CUSUM is to be preferred.

\begin{rem}
The parameter $\alpha$, which characterizes the prior knowledge on the location of possible abrupt change, plays a similar role in both Theorems \ref{CUSUMFxdDom} and \ref{CUSUMIncDom}. The fact that $\alpha$ is a fixed, strictly positive quantity, means that the detectable jump for CUSUSM in the increasing domain setting is of order $\sqrt{n^{-1}\log n}$. However, such restriction on $\alpha$ is not critical for our proof approach and the CDEP of CUSUM is still less than $\delta$, if condition \eqref{DtctRateCUSUMIncDom} holds for the case of $\alpha\rightarrow 0$ with the sample size growing. For instance even if $\alpha \asymp n^{\beta-1}$ for some $\beta\in\paren{0,1}$, CUSUSM is still consistent with the detection rate $\abs{b}\gtrsim \sqrt{n^{-\beta}\log n}$. 
\end{rem}

\section{Minimax lower bound on detection rate}\label{LowBndDtctRate}

In this section, we establish minimax lower bounds on the detectable jump in the mean of $G$. Theorem \ref{MinMaxThmFxdDom} shows that the obtained rate for GLRT and plug-in GLRT (Theorems \ref{GLRTRateFxdDom} and \ref{PlugInGLRTRateFxdDom}) are nearly optimal (up to some logarithmic term in $n$) for rational spectral densities. Section \ref{MinMaxLowBndIncDom} demonstrates the near minimax optimality of the GLRT and CUSUM algorithms in increasing domain setting. Before jumping to the main result of this section, let us rigorously introduce the notion of near minimax optimality.

\begin{defn}
Given $n$ samples, let $T\in\set{0,1}$ be a shift in mean detection algorithm whose CDEP is denoted by $\varphi_n\paren{T}$. $T$ is said to be \emph{near minimax optimal (up to some logarithmic term in $n$) in the asymptotic sense}, if for any $\delta\in\paren{0,2}$, there are two vanishing sequences $\set{h_{1,n}}^{\infty}_{n=1}$ and $\set{h_{2,n}}^{\infty}_{n=1}$ relying on $n,\delta$ and spectral density, such that
\begin{enumerate}
\item $\varphi_n\paren{T}\leq \delta$ whenever $\abs{b}\geq h_{1,n}$.
\item When $\abs{b}\leq h_{2,n}$, there is no algorithm whose CDEP is strictly less than $\delta$.
\item There exists a positive, bounded scalar $\beta$, for which $h_{1,n}/h_{2,n} = \cc{O}\paren{\log^{\beta} n}$ as $n\rightarrow\infty$.
\end{enumerate}
\end{defn}

Put simply, when the sample size is large enough, no algorithm is considerably superior to a near optimal $T$, regardless of how complicated its formulation might be.

\subsection{Lower bound in the fixed domain}\label{MinmaxLowBndFxdDom}

We begin this section by recalling that in the fixed domain regime, $G$ is a Gaussian process in $\brac{0,1}$ which is observed at $\set{i/n}^n_{i=1}$. We formally introduce the class of spectral densities that we consider in this section. While the following conditions on $\hat{K}$ are more restrictive than Assumption \ref{FxdDomAssu}, it is still provides a rich class of commonly used spectral densities.

\begin{assu}\label{MinMaxAssu}
There are constants $p\in\bb{N}$ and $\beta\in\paren{1/2,\infty}$ such that
\begin{enumerate}
\item $\lim\limits_{\omega\rightarrow\infty} \hat{K}\paren{\omega}\abs{\omega}^{2p}$ exists and $C'_K \coloneqq \lim\limits_{\omega\rightarrow\infty} \hat{K}\paren{\omega}\abs{\omega}^{2p}\in\paren{0,\infty}$.
\item $\limsup\limits_{\omega\rightarrow\infty}\abs{\paren{\frac{\hat{K}\paren{\omega}\abs{\omega}^{2p}}{C'_K}-1} \omega^{\beta}} <\infty$.
\end{enumerate}
\end{assu}

Generally speaking, Assumption \ref{MinMaxAssu} contains the class of spectral densities $\hat{K}\paren{\omega}$ for which there is some $p\in\bb{N}$ such that $\hat{K}\paren{\omega}\asymp \abs{\omega}^{-2p}$ as $\omega$ tends to infinity. Note that the second condition in Assumption \ref{MinMaxAssu} is of theoretical purposes and does not have a simple qualitative interpretation. It can be observed that Assumption \ref{MinMaxAssu} excludes any $\hat{K}\paren{\omega}$ satisfying Assumption \ref{FxdDomAssu} with $\paren{\nu+1/2}\notin \bb{N}$. For instance, Assumption\ref{MinMaxAssu} does not hold for Matern covariance functions with $\paren{\nu+1/2}\notin\bb{N}$. 

\begin{rem}\label{MinMaxOptClass}
Here, we name a salient class of spectral densities satisfying Assumption \ref{MinMaxAssu}.

\begin{itemize}
\item Simple calculations show that any rational spectral density $\hat{K}$ (See \eqref{RationalSpDen}) admits Assumption \ref{MinMaxAssu} with $C'_K = \lambda$, $\beta=1$ and $p=\deg\paren{Q_d}-\deg\paren{Q_n}\in\bb{N}$. Moreover, $\hat{K}$ satisfies Assumption \ref{FxdDomAssu} with $\nu=p-1/2$.
We discussed in Remark that Matern covariance function with $p\coloneqq\paren{\nu+1/2}\in\bb{N}$ has indeed a rational spectral density. These particular instances of Matern covariance, which are commonly used in machine learning and geostatistics, are of the form $K\paren{r} = Q\paren{\abs{r}}e^{-d\abs{r}}$, where $Q\paren{\cdot}$ is a polynomial of degree $p-1$. 
\end{itemize}

\end{rem}

\begin{thm}\label{MinMaxThmFxdDom}
Let $\delta\in\paren{0,2}$ and assume that Assumption \ref{MinMaxAssu} holds for $K$. Consider the problem \eqref{DtctProb} in which $\cov\paren{\boldsymbol{X}} = \brac{K\paren{\frac{r-s}{n}}}^n_{r,s=1}$. There are positive scalars $\bar{C}_K$ and $n_0\coloneqq n_0\paren{K}$ such that if $n\geq n_0$ and
\begin{equation}\label{LowBndRate}
\abs{b}\leq \bar{C}_Kn^{-p+1/2}\sqrt{\log\paren{\frac{1}{\delta\paren{2-\delta}}}},
\end{equation}
then for any test $T$,
\begin{equation*}
\varphi_n\paren{T}\geq\delta.
\end{equation*}
\end{thm}

See Appendix \ref{ProofThm7} for the proof of Theorem \ref{MinMaxThmFxdDom}.

\begin{rem}
Comparing the detection rate of GLRT (see Theorem \ref{GLRTRateFxdDom}) and plug-in GLRT (recall from Theorem \ref{PlugInGLRTRateFxdDom}), with the rate described in Eq. \eqref{LowBndRate} establishes near minimax optimality of the GLRT with known covariance structure and plug-in GLRT for the class of spectral densities considered in Remark \ref{MinMaxOptClass}, in the asymptotic sense. Strictly speaking, under the fixed domain setting, there is a gap of order $\sqrt{\log n}$ between \eqref{LowBndRate} and the detection rate of formerly studied GLRT based algorithms. Although we do not have a proof to establish the near minimax optimality of the GLRT and plug-in GLRT for spectral densities satisfying Assumption \ref{FxdDomAssu}, our conjecture is that Theorem \ref{MinMaxThmFxdDom} can be extended to this broader class.
\end{rem}

\subsection{Lower bound in the increasing domain setting}\label{MinMaxLowBndIncDom}

Turning to the increasing domain setting, we give a condition on jump size $\abs{b}$ according to which no algorithm in the increasing domain can properly detect the existence of a shift in the mean. Unlike Section \ref{MinmaxLowBndFxdDom}, there is no distinction between the assumptions used to obtain the minimax lower bound and Assumption \ref{AssumpToepCov}. 
 
\begin{thm}\label{MinMaxThmIncDom}
Let $\delta\in\paren{0,2}$, $\vartheta > 0$ and $\cc{C}_{n,\alpha} = \brac{\alpha n,\paren{1-\alpha}n}$. Suppose that $\Sigma_n = \cc{T}_n\paren{f}$ in which $f$ satisfies Assumption \ref{AssumpToepCov}. There are $n_0\coloneqq n_0\paren{f,\vartheta}$ and a universal constant $C>0$ such that if $n\geq n_0$ and
\begin{equation*}
\abs{b}\leq C\sqrt{\frac{\paren{1+\vartheta}f\paren{0}\log\paren{\frac{1}{\delta\paren{2-\delta}}} }{\alpha n}},
\end{equation*}
then for any test $T$,
\begin{equation*}
\varphi_n\paren{T}\geq\delta.
\end{equation*}
\end{thm}

The direct comparison between the detection rate of both CUSUM (in Theorem \ref{CUSUMIncDom}) and GLRT (see Theorem \ref{GLRTRateIncDom}) test with the above result indicates the minimax optimality (up to some order $\log n$ term) of both os these procedures in the increasing domain setting.

\section{Simulation study}\label{NumRes}

To illustrate the performance of the proposed shift-in-mean detection algorithms, we conduct a set of controlled simulation studies for verifying the results in Sections \ref{DtctRateKnwnSpDen}, \ref{DtctRateUnKnwnSpDen} and \ref{DtctCUSUM}. Our goals are two-fold:

\begin{enumerate}[label=(\alph*),leftmargin=*]
\item comparing the performance of the GLRT based algorithms with the standard CUSUM test. 
\item assessing the sensitivity of algorithm \eqref{PGLRT} to the parameters of the covariance function and tapering of $\Sigma_n$.
\end{enumerate}
In all the numerical studies in this section we fix $n=500$ and $\alpha=0.1$. 

The area under the receiver operating characteristic (ROC) curve, which will be referred as AUC, is a standard way for assessing the performance of a test. The ROC curve plots the power against the false alarm probability. Since the ROC curve is confined in the unit square, the AUC ranges in [0,1]. The ROC curve of a test based on pure random guessing is the diagonal line between origin and $(1,1)$ and so the AUC of any realistic test is at least $0.5$.

The subsequent figures in this section exhibit empirical AUC versus $b$. For a fixed value of $b$, covariance function $K$ and a detection algorithm $T$, we apply the following method to compute the AUC of $T$:

\begin{enumerate}[leftmargin=*]
\item Set $T_1=500$ and $T_2=50$.
\item For $k=1$ to $T_2$ repeat independently
\begin{enumerate}
\item For $\ell=1$ to $T_1$ repeat independently
\begin{enumerate}
\item Choose $p\in\set{0,1}$ with equal probability which denotes null or alternative hypotheses. Thus, approximately $T_1/2=250$ experiments correspond to both null and alternative.
\item If $p=0$, generate zero mean $\boldsymbol{X}\in\bb{R}^n$ according to covariance function $K$. That is, $\boldsymbol{X}$ are sampled from a Gaussian process with no abrupt shift in mean. Otherwise, choose $t\in\brac{\alpha n,\paren{1-\alpha}n}=\set{50,51,\cdots,450}$ uniformly at random (recall that $t$ represents the location of the mean shift) and generate $\boldsymbol{X}\in\bb{R}^n$ according to $\bb{H}_{1,t}$.
\item Compute $T$ score.
\end{enumerate}
\item Numerically obtain the ROC curve of $T$ based upon $T_1$ experiments in part $i$. 
\item Given the ROC curve, compute $AUC_k$ using trapezoidal integration method. 
\end{enumerate}
\item Compute the average AUC by $\overline{AUC} = \frac{1}{T_2}\sum\limits_{k=1}^{T_2} AUC_k$. 
\end{enumerate} 

\begin{figure}
\centering
\setlength\figureheight{5cm} 
\setlength\figurewidth{6.5cm} 
%
%
%
\definecolor{mycolor1}{rgb}{0.00000,0.49804,0.00000}%
\begin{tikzpicture}

\pgfplotsset{every axis/.append style={
        scaled y ticks = false, 
        scaled x ticks = false, 
        y tick label style={/pgf/number format/.cd, fixed, fixed zerofill,
                            ,precision=1},
        x tick label style={/pgf/number format/.cd, fixed, fixed zerofill,
                            precision=1}
    }
}

\begin{axis}[%
width=\figurewidth,
height=\figureheight,
scale only axis,
xmin=0.01,
xmax=0.6,
xlabel={$b$},
xmajorgrids,
ymin=0.49,
ymax=1.01,
ylabel={$\overline{AUC}$},
ymajorgrids,
name=plot5,
legend style={at={(0.6,0.5)},anchor=south west,draw=black,fill=white,legend cell align=left}
]
\addplot [color=mycolor1,dashed,line width=1.0pt,mark size=1.0pt,mark=triangle,mark options={solid}]
  table[row sep=crcr]{0.01	0.5\\
0.02	0.4955\\
0.03	0.501\\
0.04	0.5068\\
0.05	0.5031\\
0.06	0.4999\\
0.07	0.4989\\
0.08	0.5036\\
0.09	0.5096\\
0.1	0.5058\\
0.11	0.4984\\
0.12	0.5015\\
0.13	0.5093\\
0.14	0.5046\\
0.15	0.5133\\
0.16	0.5144\\
0.17	0.5098\\
0.18	0.5124\\
0.19	0.5097\\
0.2	0.5197\\
0.21	0.5146\\
0.22	0.5219\\
0.23	0.5232\\
0.24	0.5124\\
0.25	0.5176\\
0.26	0.5215\\
0.27	0.519\\
0.28	0.532\\
0.29	0.5286\\
0.3	0.5308\\
0.31	0.5369\\
0.32	0.5425\\
0.33	0.5399\\
0.34	0.5478\\
0.35	0.5474\\
0.36	0.5466\\
0.37	0.5534\\
0.38	0.5488\\
0.39	0.5516\\
0.4	0.5578\\
0.41	0.5622\\
0.42	0.5679\\
0.43	0.5587\\
0.44	0.5664\\
0.45	0.5746\\
0.46	0.5663\\
0.47	0.5687\\
0.48	0.5868\\
0.49	0.5849\\
0.5	0.5886\\
0.51	0.5816\\
0.52	0.5877\\
0.53	0.5938\\
0.54	0.5971\\
0.55	0.6023\\
0.56	0.6035\\
0.57	0.6101\\
0.58	0.6016\\
0.59	0.6142\\
0.6	0.6184\\
};
\addlegendentry{CUSUM};

\addplot [color=black,dotted,line width=1.0pt,mark size=1.0,mark=*,mark options={solid}]
  table[row sep=crcr]{0.01	0.5051\\
0.02	0.5101\\
0.03	0.4949\\
0.04	0.5087\\
0.05	0.5066\\
0.06	0.5043\\
0.07	0.5025\\
0.08	0.5172\\
0.09	0.5168\\
0.1	0.5261\\
0.11	0.5372\\
0.12	0.5532\\
0.13	0.5723\\
0.14	0.5967\\
0.15	0.6062\\
0.16	0.6415\\
0.17	0.6676\\
0.18	0.6938\\
0.19	0.7213\\
0.2	0.7481\\
0.21	0.7706\\
0.22	0.8026\\
0.23	0.8229\\
0.24	0.8526\\
0.25	0.8736\\
0.26	0.8907\\
0.27	0.9097\\
0.28	0.9263\\
0.29	0.9371\\
0.3	0.9492\\
0.31	0.9563\\
0.32	0.9644\\
0.33	0.9711\\
0.34	0.9798\\
0.35	0.9816\\
0.36	0.9872\\
0.37	0.9899\\
0.38	0.9915\\
0.39	0.9939\\
0.4	0.9949\\
0.41	0.9967\\
0.42	0.9972\\
0.43	0.9981\\
0.44	0.9988\\
0.45	0.9989\\
0.46	0.9992\\
0.47	0.9995\\
0.48	0.9997\\
0.49	0.9997\\
0.5	0.9998\\
0.51	0.9999\\
0.52	0.9999\\
0.53	0.9999\\
0.54	1\\
0.55	1\\
0.56	1\\
0.57	1\\
0.58	1\\
0.59	1\\
0.6	1\\
};
\addlegendentry{GLRT};

\addplot [color=blue,solid,line width=1.0pt,mark size=1.0pt,mark=asterisk,mark options={solid}]
  table[row sep=crcr]{0.01	0.5062\\
0.02	0.5049\\
0.03	0.4956\\
0.04	0.5056\\
0.05	0.5113\\
0.06	0.5031\\
0.07	0.5033\\
0.08	0.5188\\
0.09	0.5149\\
0.1	0.5254\\
0.11	0.5313\\
0.12	0.5492\\
0.13	0.5618\\
0.14	0.5882\\
0.15	0.5974\\
0.16	0.6288\\
0.17	0.6543\\
0.18	0.6808\\
0.19	0.7036\\
0.2	0.7327\\
0.21	0.7519\\
0.22	0.7833\\
0.23	0.7993\\
0.24	0.8294\\
0.25	0.8527\\
0.26	0.8706\\
0.27	0.8895\\
0.28	0.9049\\
0.29	0.9175\\
0.3	0.9323\\
0.31	0.9403\\
0.32	0.9505\\
0.33	0.9579\\
0.34	0.9685\\
0.35	0.9711\\
0.36	0.9783\\
0.37	0.9825\\
0.38	0.985\\
0.39	0.9876\\
0.4	0.9899\\
0.41	0.9927\\
0.42	0.9938\\
0.43	0.9957\\
0.44	0.9968\\
0.45	0.9968\\
0.46	0.9976\\
0.47	0.9983\\
0.48	0.9989\\
0.49	0.9989\\
0.5	0.9992\\
0.51	0.9995\\
0.52	0.9996\\
0.53	0.9996\\
0.54	0.9998\\
0.55	0.9999\\
0.56	0.9999\\
0.57	0.9999\\
0.58	0.9999\\
0.59	1\\
0.6	1\\
};
\addlegendentry{PGLRT};

\end{axis}

\pgfplotsset{every axis/.append style={
        scaled y ticks = false, 
        scaled x ticks = false, 
        y tick label style={/pgf/number format/.cd, fixed, fixed zerofill,
                            ,precision=1},
        x tick label style={/pgf/number format/.cd, fixed, fixed zerofill,
                            precision=2}
    }
}

\begin{axis}[%
width=\figurewidth,
height=\figureheight,
scale only axis,
xmin=0.0025,
xmax=0.05,
xlabel={$b$},
xmajorgrids,
ymin=0.49,
ymax=1.01,
ylabel={$\overline{AUC}$},
ymajorgrids,
at=(plot5.right of south east),
anchor=left of south west,
legend style={at={(0.6,0.5)},anchor=south west,draw=black,fill=white,legend cell align=left}
]
\addplot [color=mycolor1,dashed,line width=1.0pt,mark size=1.0pt,mark=triangle,mark options={solid}]
  table[row sep=crcr]{0.0025	0.4972\\
0.005	0.5048\\
0.0075	0.4988\\
0.01	0.4993\\
0.0125	0.4985\\
0.015	0.4979\\
0.0175	0.5027\\
0.02	0.4982\\
0.0225	0.4952\\
0.025	0.4974\\
0.0275	0.5011\\
0.03	0.5034\\
0.0325	0.4975\\
0.035	0.5051\\
0.0375	0.507\\
0.04	0.5055\\
0.0425	0.5033\\
0.045	0.5017\\
0.0475	0.4992\\
0.05	0.504\\
};
\addlegendentry{CUSUM};

\addplot [color=black,dotted,line width=1.0pt,mark size=1.0,mark=*,mark options={solid}]
  table[row sep=crcr]{0.0025	0.496\\
0.005	0.4961\\
0.0075	0.5164\\
0.01	0.5356\\
0.0125	0.5913\\
0.015	0.6506\\
0.0175	0.7395\\
0.02	0.8099\\
0.0225	0.8699\\
0.025	0.92\\
0.0275	0.9507\\
0.03	0.9726\\
0.0325	0.9859\\
0.035	0.9937\\
0.0375	0.997\\
0.04	0.9986\\
0.0425	0.9995\\
0.045	0.9998\\
0.0475	0.9999\\
0.05	1\\
};
\addlegendentry{GLRT};

\addplot [color=blue,solid,line width=1.0pt,mark size=1.0pt,mark=asterisk,mark options={solid}]
  table[row sep=crcr]{0.0025	0.497\\
0.005	0.4926\\
0.0075	0.5141\\
0.01	0.5288\\
0.0125	0.5831\\
0.015	0.6369\\
0.0175	0.723\\
0.02	0.7878\\
0.0225	0.8491\\
0.025	0.9007\\
0.0275	0.9345\\
0.03	0.9593\\
0.0325	0.9763\\
0.035	0.9877\\
0.0375	0.993\\
0.04	0.9963\\
0.0425	0.9983\\
0.045	0.9991\\
0.0475	0.9996\\
0.05	0.9998\\
};
\addlegendentry{PGLRT};

\end{axis}

\pgfplotsset{every axis/.append style={
        scaled y ticks = false, 
        scaled x ticks = false, 
        y tick label style={/pgf/number format/.cd, fixed, fixed zerofill,
                            ,precision=1},
        x tick label style={/pgf/number format/.cd, fixed, fixed zerofill,
                            precision=3}
    }
}

\begin{axis}[%
width=\figurewidth,
height=\figureheight,
scale only axis,
xmin=0.0015,
xmax=0.03,
xlabel={$b$},
xmajorgrids,
ymin=0.49,
ymax=1.01,
ylabel={$\overline{AUC}$},
ymajorgrids,
at=(plot5.below south west),
anchor=above north west,
legend style={at={(0.6,0.5)},anchor=south west,draw=black,fill=white,legend cell align=left}
]
\addplot [color=mycolor1,dashed,line width=1.0pt,mark size=1.0pt,mark=triangle,mark options={solid}]
  table[row sep=crcr]{0.0015	0.4937\\
0.003	0.4958\\
0.0045	0.4948\\
0.006	0.5036\\
0.0075	0.5004\\
0.009	0.5016\\
0.0105	0.5099\\
0.012	0.5012\\
0.0135	0.4961\\
0.015	0.5023\\
0.0165	0.4995\\
0.018	0.4993\\
0.0195	0.4988\\
0.021	0.4953\\
0.0225	0.5047\\
0.024	0.5006\\
0.0255	0.5002\\
0.027	0.4977\\
0.0285	0.4981\\
0.03	0.4994\\
};
\addlegendentry{CUSUM};

\addplot [color=black,dotted,line width=1.0pt,mark size=1.0,mark=*,mark options={solid}]
  table[row sep=crcr]{0.0015	0.5058\\
0.003	0.5184\\
0.0045	0.5892\\
0.006	0.714\\
0.0075	0.8346\\
0.009	0.9243\\
0.0105	0.9694\\
0.012	0.9904\\
0.0135	0.9972\\
0.015	0.9995\\
0.0165	0.9999\\
0.018	1\\
0.0195	1\\
0.021	1\\
0.0225	1\\
0.024	1\\
0.0255	1\\
0.027	1\\
0.0285	1\\
0.03	1\\
};
\addlegendentry{GLRT};

\addplot [color=blue,solid,line width=1.0pt,mark size=1.0pt,mark=asterisk,mark options={solid}]
  table[row sep=crcr]{0.0015	0.5092\\
0.003	0.518\\
0.0045	0.5818\\
0.006	0.6948\\
0.0075	0.8117\\
0.009	0.9055\\
0.0105	0.9561\\
0.012	0.9823\\
0.0135	0.9938\\
0.015	0.9981\\
0.0165	0.9994\\
0.018	0.9999\\
0.0195	1\\
0.021	1\\
0.0225	1\\
0.024	1\\
0.0255	1\\
0.027	1\\
0.0285	1\\
0.03	1\\
};
\addlegendentry{PGLRT};

\end{axis}
\end{tikzpicture}%
\caption{The above figures assess the performance of different detection algorithms when $G$ is one dimensional Matern Gaussian process, with parameters $\paren{\nu,\sigma_0,\rho_0}$, and regularly sampled in $\brac{0,1}$. From left to right then from top to bottom, $\paren{\nu,\sigma_0,\rho_0}=\paren{0.5,1,0.5},\paren{1,1,0.5},\paren{1.5,1,0.5}$. In each panel horizontal axis displays jump value $b$ and the three curves (dashed black, solid blue and green) respectively exhibit AUC of GLRT with known covariance structure, plug in GLRT (PGLRT) using full MLE and CUSUM.}
\label{Fig:Fig1}
\end{figure}
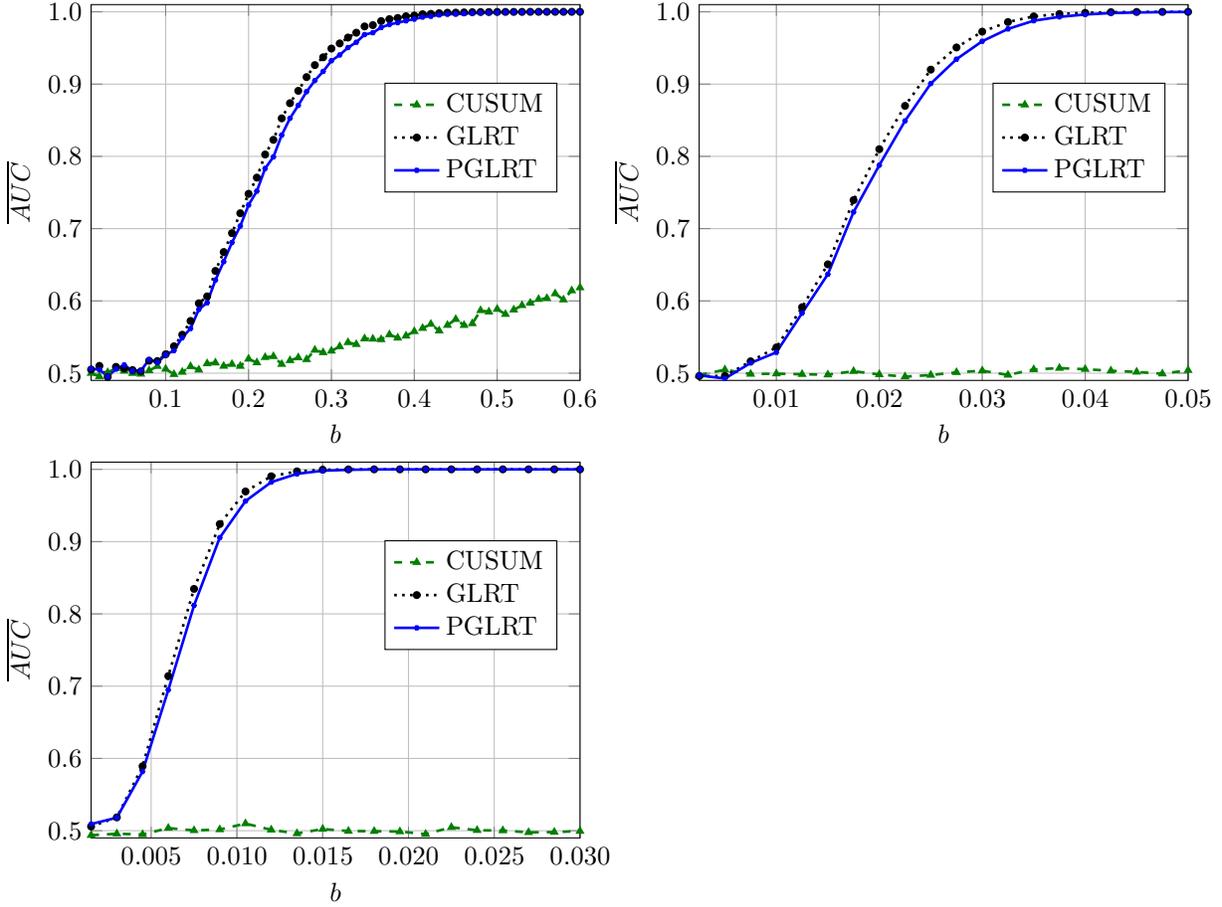

The first simulation study aims to compare CUSUM and GLRT based algorithms in the fixed domain regime and assess the role of smoothness and other parameters of $K$ in the performance of the GLRT. For this experiment $G$ is a Gaussian process in $\brac{0,1}$ which is observed at regularly spaced samples, $\cc{D}_n = \set{k/n}^n_{k=1}$, i.e., $X_k = G\paren{k/n},\;k=1,\ldots,n$ . The covariance function of $G$ is assumed to has Matern form with parameters $\paren{\sigma_0,\rho_0,\nu}$. Strictly speaking,
\begin{align*}
&\cov\paren{X_i,X_l} = \sigma^2_0 K_{\nu}\paren{\abs{\frac{i-l}{n\rho_0}}}, \quad i,l = 1,\ldots,n,\\
&K_{\nu}\paren{x} = \frac{\sqrt{4\pi}\Gamma\paren{\nu+1/2}}{\Gamma\paren{\nu}} \int_{-\infty}^{\infty} e^{-j\omega x} \paren{1+\omega^2}^{-\paren{\nu+1/2}} d\omega,\quad \forall x\geq 0.
\end{align*}
We consider three different scenarios on $\nu$, $0.5,1,$ and $1.5$. We also set $\sigma_0 = 1$ and $\rho_0 = 1/2$. As customary in the literature, we assume that $\nu$ is known and so $\nu$ will not be estimated. For conducting the plug-in GLRT procedure, both parameters $\paren{\sigma_0,\rho_0}$ are estimated using full MLE. Due to the low dimensionality of unknown parameters, the most effective way to estimate $\paren{\sigma_0,\rho_0}$ is to apply brute force grid search over a pre-specified set $\cc{P}$. Here, we choose $\cc{P} = \set{0.2,0.4,\cdots,2}\times\set{1/4,1/3.9,\ldots,1/0.1}$. The final results of this numerical study is exhibited in Figure \ref{Fig:Fig1}. We observe the following:
\begin{itemize}
\item GLRT and plug-in GLRT have a significantly better detection performance than CUSUM. This performance improvement is more pronounced for smoother covariance function (larger $\nu$). In particular, the CUSUM is completely impractical for detection of an small change when $\nu=1$ or $1.5$.
\item In each panel of Figure \ref{Fig:Fig1}, the GLRT has a slightly larger AUC than that of plug-in GLRT. Thus, Figure \ref{Fig:Fig1} verifies the existence of a small gap between the smallest detectable jump of GLRT and plug-in GLRT. Note that this fact can also be observed by comparing \eqref{LowBndDtctbaleJmp} and \eqref{ALGRTRate}. In short, having full knowledge of covariance parameters slightly improves the detection performance and so our proposed algorithm is robust to the estimation error of the unknown parameters of $K$. 
\item Comparing the range of $b$ in each panel of Figure \ref{Fig:Fig1} discloses that more rapid decay of the spectral density can decrease the smallest detectable jump. This observation substantiates the role of $\nu$ in the theory established in Sections \ref{DtctRateKnwnSpDen} and \ref{DtctRateUnKnwnSpDen}.
\end{itemize}

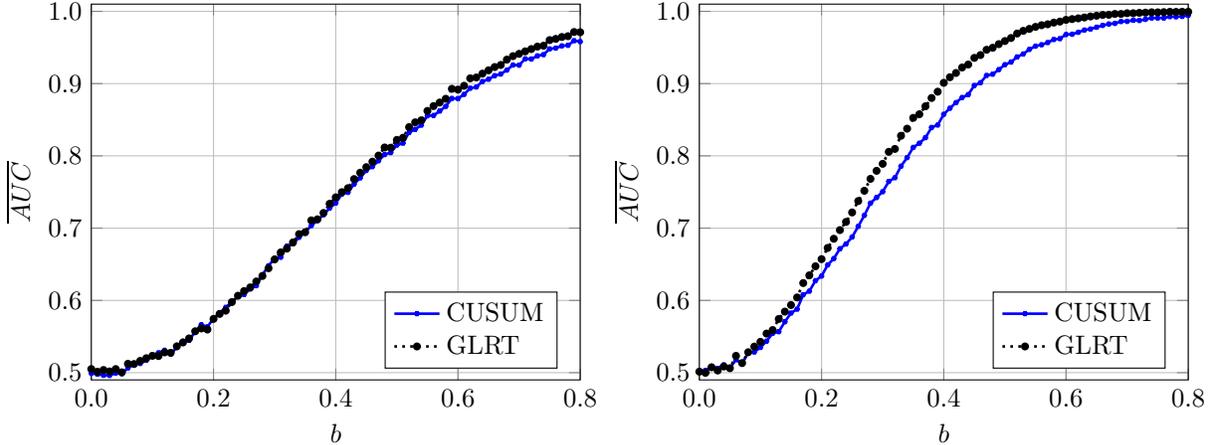
\begin{figure}
\centering
\setlength\figureheight{5cm} 
\setlength\figurewidth{6.5cm} 
%
%

\pgfplotsset{every axis/.append style={
        scaled y ticks = false, 
        scaled x ticks = false, 
        y tick label style={/pgf/number format/.cd, fixed, fixed zerofill,
                            ,precision=1},
        x tick label style={/pgf/number format/.cd, fixed, fixed zerofill,
                            precision=1}
    }
}

\begin{tikzpicture}

\begin{axis}[%
width=\figurewidth,
height=\figureheight,
scale only axis,
xmin=0,
xmax=0.8,
xlabel={$b$},
xmajorgrids,
ymin=0.49,
ymax=1.01,
ylabel={$\overline{AUC}$},
ymajorgrids,
name=plot3,
legend style={at={(0.6,0.03)},anchor=south west,draw=black,fill=white,legend cell align=left}
]
\addplot [color=blue,solid,line width=1.0pt,mark size=1.0pt,mark=asterisk,mark options={solid}]
  table[row sep=crcr]{0	0.4991\\
0.01	0.4979\\
0.02	0.4964\\
0.03	0.4962\\
0.04	0.4994\\
0.05	0.4986\\
0.06	0.5065\\
0.07	0.5117\\
0.08	0.5127\\
0.09	0.5171\\
0.1	0.5243\\
0.11	0.5275\\
0.12	0.531\\
0.13	0.5259\\
0.14	0.5349\\
0.15	0.5434\\
0.16	0.5448\\
0.17	0.5566\\
0.18	0.5664\\
0.19	0.5634\\
0.2	0.5751\\
0.21	0.5827\\
0.22	0.5904\\
0.23	0.5978\\
0.24	0.6054\\
0.25	0.6082\\
0.26	0.6158\\
0.27	0.6204\\
0.28	0.6354\\
0.29	0.6479\\
0.3	0.6581\\
0.31	0.6597\\
0.32	0.6751\\
0.33	0.6783\\
0.34	0.6873\\
0.35	0.6943\\
0.36	0.7035\\
0.37	0.7105\\
0.38	0.7186\\
0.39	0.7277\\
0.4	0.7346\\
0.41	0.747\\
0.42	0.7495\\
0.43	0.7609\\
0.44	0.7695\\
0.45	0.7798\\
0.46	0.7851\\
0.47	0.7932\\
0.48	0.8018\\
0.49	0.8045\\
0.5	0.8149\\
0.51	0.8176\\
0.52	0.8324\\
0.53	0.8367\\
0.54	0.8425\\
0.55	0.8555\\
0.56	0.8559\\
0.57	0.862\\
0.58	0.8685\\
0.59	0.8795\\
0.6	0.8793\\
0.61	0.8851\\
0.62	0.8936\\
0.63	0.8953\\
0.64	0.9029\\
0.65	0.906\\
0.66	0.911\\
0.67	0.9132\\
0.68	0.9187\\
0.69	0.9257\\
0.7	0.9255\\
0.71	0.9343\\
0.72	0.934\\
0.73	0.9384\\
0.74	0.9401\\
0.75	0.948\\
0.76	0.9491\\
0.77	0.9522\\
0.78	0.953\\
0.79	0.9594\\
0.8	0.9581\\
};
\addlegendentry{CUSUM};

\addplot [color=black,dotted,line width=1.0pt,mark size=1.0,mark=*,mark options={solid}]
  table[row sep=crcr]{0	0.505\\
0.01	0.5013\\
0.02	0.5037\\
0.03	0.5017\\
0.04	0.505\\
0.05	0.5005\\
0.06	0.5123\\
0.07	0.512\\
0.08	0.5163\\
0.09	0.52\\
0.1	0.5233\\
0.11	0.523\\
0.12	0.5281\\
0.13	0.5275\\
0.14	0.5365\\
0.15	0.5418\\
0.16	0.5474\\
0.17	0.5573\\
0.18	0.5612\\
0.19	0.5597\\
0.2	0.5744\\
0.21	0.5813\\
0.22	0.5858\\
0.23	0.5977\\
0.24	0.6065\\
0.25	0.6129\\
0.26	0.618\\
0.27	0.6262\\
0.28	0.6338\\
0.29	0.6445\\
0.3	0.6567\\
0.31	0.6664\\
0.32	0.6717\\
0.33	0.6802\\
0.34	0.6917\\
0.35	0.6943\\
0.36	0.7107\\
0.37	0.7122\\
0.38	0.721\\
0.39	0.7337\\
0.4	0.7426\\
0.41	0.7498\\
0.42	0.7553\\
0.43	0.7679\\
0.44	0.7767\\
0.45	0.7844\\
0.46	0.792\\
0.47	0.8004\\
0.48	0.8116\\
0.49	0.8113\\
0.5	0.8221\\
0.51	0.825\\
0.52	0.8398\\
0.53	0.8464\\
0.54	0.8492\\
0.55	0.8622\\
0.56	0.8689\\
0.57	0.8738\\
0.58	0.879\\
0.59	0.8928\\
0.6	0.8917\\
0.61	0.897\\
0.62	0.9073\\
0.63	0.9085\\
0.64	0.9138\\
0.65	0.9184\\
0.66	0.9229\\
0.67	0.9258\\
0.68	0.9332\\
0.69	0.938\\
0.7	0.9412\\
0.71	0.9446\\
0.72	0.9478\\
0.73	0.9512\\
0.74	0.9524\\
0.75	0.9602\\
0.76	0.9619\\
0.77	0.9645\\
0.78	0.9656\\
0.79	0.9714\\
0.8	0.9711\\
};
\addlegendentry{GLRT};

\end{axis}

\begin{axis}[%
width=\figurewidth,
height=\figureheight,
scale only axis,
xmin=0,
xmax=0.8,
xlabel={$b$},
xmajorgrids,
ymin=0.49,
ymax=1.01,
ylabel={$\overline{AUC}$},
ymajorgrids,
at=(plot3.right of south east),
anchor=left of south west,
legend style={at={(0.6,0.03)},anchor=south west,draw=black,fill=white,legend cell align=left}
]
\addplot [color=blue,solid,line width=1.0pt,mark size=1.0pt,mark=asterisk,mark options={solid}]
  table[row sep=crcr]{0	0.5024\\
0.01	0.503\\
0.02	0.5087\\
0.03	0.5059\\
0.04	0.5108\\
0.05	0.5064\\
0.06	0.518\\
0.07	0.5124\\
0.08	0.5279\\
0.09	0.5286\\
0.1	0.5349\\
0.11	0.5436\\
0.12	0.5544\\
0.13	0.5567\\
0.14	0.5705\\
0.15	0.5825\\
0.16	0.5879\\
0.17	0.6079\\
0.18	0.6128\\
0.19	0.6273\\
0.2	0.6337\\
0.21	0.6491\\
0.22	0.6579\\
0.23	0.6715\\
0.24	0.6781\\
0.25	0.6877\\
0.26	0.7025\\
0.27	0.7177\\
0.28	0.7346\\
0.29	0.7424\\
0.3	0.7504\\
0.31	0.7646\\
0.32	0.7698\\
0.33	0.7857\\
0.34	0.7976\\
0.35	0.8117\\
0.36	0.8175\\
0.37	0.8252\\
0.38	0.8393\\
0.39	0.8427\\
0.4	0.8576\\
0.41	0.8659\\
0.42	0.8737\\
0.43	0.8808\\
0.44	0.8848\\
0.45	0.8974\\
0.46	0.9012\\
0.47	0.9113\\
0.48	0.913\\
0.49	0.9194\\
0.5	0.9262\\
0.51	0.9298\\
0.52	0.9369\\
0.53	0.941\\
0.54	0.9473\\
0.55	0.9518\\
0.56	0.9538\\
0.57	0.9568\\
0.58	0.961\\
0.59	0.9622\\
0.6	0.9679\\
0.61	0.9683\\
0.62	0.971\\
0.63	0.9741\\
0.64	0.9755\\
0.65	0.978\\
0.66	0.9805\\
0.67	0.9825\\
0.68	0.9833\\
0.69	0.9859\\
0.7	0.9861\\
0.71	0.9877\\
0.72	0.9873\\
0.73	0.9889\\
0.74	0.9908\\
0.75	0.9908\\
0.76	0.9908\\
0.77	0.9927\\
0.78	0.9924\\
0.79	0.993\\
0.8	0.9944\\
};
\addlegendentry{CUSUM};

\addplot [color=black,dotted,line width=1.0pt,mark size=1.0,mark=*,mark options={solid}]
  table[row sep=crcr]{0	0.5014\\
0.01	0.4998\\
0.02	0.5075\\
0.03	0.503\\
0.04	0.5085\\
0.05	0.5062\\
0.06	0.5233\\
0.07	0.5133\\
0.08	0.5283\\
0.09	0.5361\\
0.1	0.5427\\
0.11	0.5541\\
0.12	0.559\\
0.13	0.5743\\
0.14	0.5845\\
0.15	0.5936\\
0.16	0.6041\\
0.17	0.624\\
0.18	0.6347\\
0.19	0.6472\\
0.2	0.6571\\
0.21	0.6726\\
0.22	0.6853\\
0.23	0.6973\\
0.24	0.7086\\
0.25	0.7217\\
0.26	0.7375\\
0.27	0.7517\\
0.28	0.7682\\
0.29	0.7793\\
0.3	0.789\\
0.31	0.8056\\
0.32	0.8096\\
0.33	0.8278\\
0.34	0.8376\\
0.35	0.8526\\
0.36	0.8576\\
0.37	0.8686\\
0.38	0.8801\\
0.39	0.8888\\
0.4	0.9011\\
0.41	0.9088\\
0.42	0.9147\\
0.43	0.9223\\
0.44	0.9264\\
0.45	0.9356\\
0.46	0.9395\\
0.47	0.9468\\
0.48	0.9496\\
0.49	0.9546\\
0.5	0.9589\\
0.51	0.9631\\
0.52	0.9693\\
0.53	0.973\\
0.54	0.9757\\
0.55	0.9785\\
0.56	0.9808\\
0.57	0.9821\\
0.58	0.9844\\
0.59	0.9858\\
0.6	0.9882\\
0.61	0.9893\\
0.62	0.9903\\
0.63	0.9915\\
0.64	0.9929\\
0.65	0.9939\\
0.66	0.995\\
0.67	0.9952\\
0.68	0.9964\\
0.69	0.9964\\
0.7	0.9974\\
0.71	0.9974\\
0.72	0.9979\\
0.73	0.9982\\
0.74	0.9987\\
0.75	0.9986\\
0.76	0.9989\\
0.77	0.9993\\
0.78	0.9992\\
0.79	0.9994\\
0.8	0.9994\\
};
\addlegendentry{GLRT};

\end{axis}
\end{tikzpicture}%
\caption{The above figure assesses the performance of increasing domain detection algorithms. In each panel the horizontal axis displays the jump value $b$ and the two curves (dashed black and solid blue) respectively exhibit the AUC of the GLRT with known covariance structure and CUSUM. In the right panel, we choose $\cov\paren{X_i,X_l} = \sigma^2_0\paren{1+\abs{i-l}/\rho_0}^{-\paren{1+\lambda}}$ in which $\paren{\sigma_0,\rho_0}=\paren{1,2}$ and $\lambda=0.5$. For the left panel, the covariance function is given by $\cov\paren{X_i,X_l} = \sigma^2_0\exp\paren{-\abs{i-l}/\rho_0}$ where $\paren{\sigma_0,\rho_0}=\paren{1,2}$.}
\label{Fig:Fig4}
\end{figure}

Next, we compare the performance of the GLRT with known parameters and the CUSUM in the increasing domain setting. Recalling from Theorems \ref{GLRTRateIncDom} and \ref{CUSUMIncDom}, these two methods have analogous asymptotic rates. In the left panel of Figure \ref{Fig:Fig4}, we choose an exponentially decaying covariance function 
\begin{equation*}
\cov\paren{X_i,X_l} = f_{\abs{i-l}} = \sigma^2_0\exp\paren{-\frac{\abs{i-l}}{\rho_0}},\quad\;i,l = 1,\ldots,n,
\end{equation*}
in which $\sigma_0 = 1$ and $\rho_0 = 2$. That is $\Sigma_n$ has exponentially decaying off-diagonal entries. However, in the right panel, the chosen covariance function has a polynomially decaying tail given by 
\begin{equation*}
\cov\paren{X_i,X_j} = f_{\abs{i-j}} = \sigma^2_0 \paren{1+\frac{\abs{i-l}}{\rho_0}}^{-\paren{1+\lambda}},
\end{equation*}
with $\sigma_0 = 1$, $\rho_0 = 2$ and $\lambda = 0.5$. In this case, $\Sigma_n$ has heavier off-diagonal terms. Note that Assumption \ref{AssumpToepCov} is satisfied in either of the two cases. It is evident from Figure \ref{Fig:Fig4} that the GLRT exhibits a slightly better performance than the CUSUM, and the gap between the two AUC curves is more visible in the case of polynomially decaying covariance function. Thus, we still recommend the use of GLRT in the presence of strong correlation among samples in applications described by the increasing domain regime.

\section{Discussion}\label{Discussion}

As indicated in the Introduction, the comprehensive analysis of the detection of shift-in-mean of a Gaussian process in the fixed domain regime has remained relatively unexplored. However, the considered model in \eqref{ProbForm} is only one of several plausible scenarios which should be subject to more thorough investigation. We note that the probabilistic model of $G$ can be extended in some possible ways for future research. 

\begin{enumerate}[label = (\alph*),leftmargin=*]
\item Here we deal with a single abrupt change in $\bb{E}G$. However, we believe that our techniques can be extended to rigorously formulate the minimax optimal rate of the GLRT and plug-in GLRT for 
detecting multiple shifts in $\bb{E}G$.
\item Recently, Ivanoff et al. \cite{BG.Ivanoff} studied the problem of \emph{change-set} detection in two dimensional Poisson processes. Specifically, $G$ is a Poisson process in $\bb{R}^2$ and there are scalars $\mu_0\ne \mu_1$ and $\Omega\subset\bb{R}^2$ such that the intensity of $G$ can be formulated by
\begin{equation}\label{ChnageSet}
\bb{E}G\paren{s} = \mu_0\bbM{1}_{s\in\Omega}+\mu_1\bbM{1}_{s\notin\Omega},\quad\in\bb{R}^2.
\end{equation}
The objective is to detect $\Omega$ as well as possible based on observation of $G$. However, a comprehensive study of minimax optimal change-set detection methods for multi-dimensional Gaussian processes remains unavailable. Note that the fixed domain setting is the natural way to study asymptotic behaviour of algorithms regarding spatial processes. So, this paper can provide valuable intuition about the minimax rate of change-set detection in Gaussian spatial processes. We expect that aside from the sample size and smoothness of the covariance functions, the geometric properties of change-sets will have a crucial role in the design and analysis of detection algorithms.
\end{enumerate}

\appendix
\appendixpage
\makeatletter
\def\@seccntformat#1{\csname Pref@#1\endcsname \csname the#1\endcsname\quad}
\def\Pref@section{~}
\makeatother

\counterwithin{thm}{section}
\counterwithin{assu}{section}
\counterwithin{lem}{section}
\counterwithin{cor}{section}
\counterwithin{prop}{section}

Appendix \ref{Proofs} contains the proofs of the main results in Sections \ref{DtctAlg}-\ref{LowBndDtctRate}. Appendix \ref{AuxRes} states and proves the technical results required in the proofs of the main results. In addition, Appendix \ref{NonAsympInvLargToep} establishes properties of the 
inverses of large Toeplitz matrices. Such results are useful for the proofs of 
Theorems \ref{GLRTRateIncDom}, \ref{CUSUMIncDom} and \ref{MinMaxThmIncDom}, and may also be of independent interest.

\section{Proofs}\label{Proofs}

\subsection{Proofs for Section \ref{DtctAlg}}

\begin{proof}[Proof of Proposition \ref{GLRTProp}]
In the following $\ff{L}$ stands for the generalized negative log-likelihood ratio.
\begin{equation}\label{GLR}
2\ff{L} = \boldsymbol{X}^\top\paren{\Sigma_n}^{-1}\boldsymbol{X}
- \min_{t\in\cc{C}_{n,\alpha}}\min_{b\ne 0}\brac{\paren{\boldsymbol{X}-\frac{b}{2}\zeta_{t}}^\top\paren{\Sigma_n}^{-1}\paren{\boldsymbol{X}-\frac{b}{2}\zeta_{t}}}.
\end{equation}
Note that the objective function in \eqref{GLR} is quadratic in terms of $b$. The explicit form of $2\ff{L}$ can be obtained with a bit of algebraic derivations. The algebra has been skipped to save space; we arrive at
\begin{equation*}
2\ff{L} = \max_{t\in\cc{C}_{n,\alpha}}\max_{b\ne 0} \paren{-\frac{\zeta^\top_t\paren{\Sigma_n}^{-1}\zeta_t}{4}b^2+ b \zeta^\top_t\paren{\Sigma_n}^{-1}\boldsymbol{X}} = \max_{t\in\cc{C}_{n,\alpha}} \abs{\frac{\zeta^\top_t\paren{\Sigma_n}^{-1}\boldsymbol{X}}{\sqrt{\zeta^\top_t\paren{\Sigma_n}^{-1}\zeta_t}}}^2.
\end{equation*}
So, there is a threshold value, $R_{n,\delta}>0$, for which the GLRT is given by \eqref{GLRT}.
\end{proof}

The following result expressing the form of GLRT in the generic case of unknown $\mu$ can be proved in an analogous way as Proposition \ref{GLRTProp}.

\begin{prop}\label{GLRTGenForm}
There is $R_{n,\delta}>0$ for which the GLRT is given by
\begin{equation}
T_{GLRT} = \bb{I}\paren{\max_{t\in\cc{C}_{n,\alpha}} \abs{\frac{\InnerProd{\boldsymbol{Y}}{\zeta_{t}- B_1\paren{t}\bbM{1}_n}{}}{\sqrt{B_2\paren{t}}}}^2\geq R_{n,\delta}},
\end{equation}
where $\boldsymbol{Y} = \paren{\Sigma_n}^{-1}\boldsymbol{X}$ and
\begin{equation*}
B_1\paren{t} = \frac{\zeta_{t}^\top\paren{\Sigma_n}^{-1}\bbM{1}_n}{\bbM{1}_n^\top\paren{\Sigma_n}^{-1}\bbM{1}_n},\quad
B_2\paren{t} = \zeta_{t}^\top\paren{\Sigma_n}^{-1}\zeta_{t} - \frac{\paren{\zeta_{t}^\top\paren{\Sigma_n}^{-1}\bbM{1}_n}^2}{\bbM{1}_n^\top\paren{\Sigma_n}^{-1}\bbM{1}_n}.
\end{equation*}
\end{prop}

\subsection{Proofs for Section \ref{DtctRateKnwnSpDen}}

\begin{proof}[Proof of Theorem \ref{GLRTRateFxdDom}]\label{ProofThm1}
Let $p = \lceil\nu+1/2\rceil$, $P=\set{1,\ldots,p}$ and $\theta_n = \exp\paren{-1/n}$. Construct a banded triangular matrix $A_n\in\bb{R}^{n\times n}$ by the following procedure.
\begin{align*}
&A_n\brac{k,k-j}= {p \choose j}\paren{-\theta_n}^j,\; j\in\set{0,\dots,p},\;k\in\set{p+1,\ldots,n},\\
&\paren{A_n}_{P,P} = n^{-2\nu}I_p.
\end{align*} 
It is relatively simple to verify that $A_n$ is invertible. In addition, for brevity let $Z_t = \frac{\zeta^\top_t\paren{\Sigma_n}^{-1}\boldsymbol{X}}{\sqrt{\zeta^\top_t\paren{\Sigma_n}^{-1}\zeta_t}}$ for any $t\in\cc{C}_{n,\alpha}$, in which $\zeta_t$ has been defined in \eqref{H1t}. Lastly, define $U_{n,t}\coloneqq A_n\zeta_t\in\bb{R}^n$, $W\coloneqq A_n\boldsymbol{X}$ and $D_n \coloneqq \cov\paren{W}$. 

Easy calculations show that under the null hypothesis $\set{Z_t}_{t\in\cc{C}_{n,\alpha}}$ is a set of standard Gaussian random variables and so, by Lemma \ref{ConcIneqNonCentChi2}, we have
$\bb{P}\paren{\max_{t\in\cc{C}_{n,\alpha}} Z^2_t\geq R_{n,\delta}}\leq \delta/2$. That is, the false alarm probability is less than $\delta/2$. Moreover if the alternative hypothesis $\bb{H}_{1,\tilde{t}}$ (for some $\tilde{t}\in\cc{C}_{n,\alpha}$) holds then $\set{Z^2_t}_{t\in\cc{C}_{n,\alpha}}$ are non-central $\chi^2_1$ random variables and the non-centrality parameter of $Z^2_{\tilde{t}}$ is given by 
\begin{equation*}
\bb{E}\paren{Z_{\tilde{t}}\mid \bb{H}_{1,\tilde{t}}} = \frac{\abs{b}}{2} \sqrt{\zeta^\top_{\tilde{t}}\paren{\Sigma_n}^{-1}\zeta_{\tilde{t}}}.
\end{equation*}
Applying Lemma \ref{ConcIneqNonCentChi2} ($\sigma_0 = \sigma_k=1$ for any $k$) demonstrates that $\varphi_n\paren{T_2}\leq \delta$, whenever 
\begin{equation}\label{LowBndExpValH1}
\abs{b} \sqrt{\zeta^\top_{\tilde{t}}\paren{\Sigma_n}^{-1}\zeta_{\tilde{t}}}\geq \abs{b} \min_{t\in\cc{C}_{n,\alpha}}\sqrt{\zeta^\top_t\paren{\Sigma_n}^{-1}\zeta_t} \geq 8\sqrt{\log\paren{\frac{4n}{\delta}}}.
\end{equation}
Thus, in order to get a sufficient condition on detectable $b$, it suffices to find a tight uniform lower bound on $\zeta^\top_t\paren{\Sigma_n}^{-1}\zeta_t$ for $t\in\cc{C}_{n,\alpha}$. 

The identity $\Sigma^{-1}_n=A^\top_n\paren{D_n}^{-1}A_n$ can be shown using the linearity of covariance operator and non-singularity of $A$. Choose $t\in\cc{C}_{n,\alpha}$ in an arbitrary way. As a result of this alternative representation of $\Sigma^{-1}_n$, we have $\zeta^\top_t\paren{\Sigma_n}^{-1}\zeta_t = U^\top_{n,t}\paren{D_n}^{-1} U_{n,t}$. Applying \emph{Kantorovich} inequality (cf. Appendix B) and the triangle inequality yields
\begin{equation}\label{LowBndQuadTerm}
\zeta^\top_t\paren{\Sigma_n}^{-1}\zeta_t = U^\top_{n,t}\paren{D_n}^{-1} U_{n,t} \geq \frac{\LpNorm{U_{n,t}}{2}^4}{U^\top_{n,t}D_n U_{n,t}} \geq \brac{\frac{\LpNorm{U_{n,t}}{2}^2}{\LpNorm{U_{n,t}}{1}}}^2 \frac{1}{\LpNorm{D_n}{\infty}}.
\end{equation}
Now, we show that $\frac{\LpNorm{U_{n,t}}{2}^2}{\LpNorm{U_{n,t}}{1}}\geq \frac{1}{3}$, for large enough $n$. 
Indeed, after some algebra, we can get
\begin{eqnarray}\label{LowBndL2Norm}
\LpNorm{U_{n,t}}{2}^2 &\geq& \sum\limits_{k=t+1}^{t+p} U^2_{n,t}\paren{k} = \sum\limits_{k=1}^{p} \brac{-\paren{1-\theta_n}^{p}+2\sum\limits_{j=0}^{k-1}{p \choose j}\paren{-\theta_n}^j}^2\nonumber\\
&\RelNum{\paren{a}}{\geq}& 2\sum\limits_{k=1}^{p} \brac{\sum\limits_{j=0}^{k-1}{p \choose j}\paren{-1}^j}^2 = 2\sum\limits_{k=1}^{p} \brac{{p-1\choose k-1}\paren{-1}^{k-1}}^2 = 2{2\paren{p-1}\choose p-1}\geq 2^p,
\end{eqnarray}
where inequality $\paren{a}$ follows from the fact that for large enough $n$, $\theta_n$ is arbitrarily close to $1$. To get an upper bound on $\LpNorm{U_{n,t}}{1}$,
\begin{eqnarray}\label{UppBndL1Norm}
\LpNorm{U_{n,t}}{1} &=& \sum\limits_{k=1}^{p} \abs{U_{n,t}\paren{k}} + \sum\limits_{k=p+1}^{t} \abs{U_{n,t}\paren{k}} + \sum\limits_{k=t+p+1}^{n} \abs{U_{n,t}\paren{k}} + \sum\limits_{k=t+1}^{t+p} \abs{U_{n,t}\paren{k}}\nonumber\\
&=&\sum\limits_{k=1}^{p} n^{-2\nu} + \sum\limits_{k=p+1}^{t} \paren{1-\theta_n}^{p} + \sum\limits_{k=t+p+1}^{n} \paren{1-\theta_n}^{p} + \sum\limits_{k=1}^{p} \abs{-\paren{1-\theta_n}^{p}+2\sum\limits_{j=0}^{k-1}{p \choose j}\paren{-\theta_n}^j}\nonumber\\
&\leq& pn^{-2\nu} + n^{1-p} + 2\sum\limits_{k=1}^{p} \abs{\sum\limits_{j=0}^{k-1}{p \choose j}\paren{-\theta_n}^j}\RelNum{\paren{b}}{\leq} 2+2\sum\limits_{k=1}^{p} \abs{\sum\limits_{j=0}^{k-1}{p \choose j}\paren{-\theta_n}^j}\nonumber\\
&\leq& 2 + 4\sum\limits_{k=1}^{p} \abs{\sum\limits_{j=0}^{k-1}{p \choose j}\paren{-1}^j} = 2+4\sum\limits_{k=1}^{p}\abs{{p-1\choose k-1}\paren{-1}^{k-1}} = 2+2^{p+1} \leq 3\;2^p.
\end{eqnarray}
Note that inequality $\paren{b}$ is valid when $pn^{-2\nu} + n^{1-p}\leq 2$, which obviously holds for sufficiently large $n=\cc{O}\paren{1}$. The remaining inequalities and identities in \eqref{UppBndL1Norm} can be easily verified via
basic properties of the binomial coefficients. Combining \eqref{LowBndL2Norm} and \eqref{UppBndL1Norm} yields the desired goal. Now, inequality \eqref{LowBndQuadTerm} can be rewritten as
\begin{equation}\label{LowBndQuadTerm2}
\zeta^\top_t\paren{\Sigma_n}^{-1}\zeta_t \geq \frac{1}{9\LpNorm{D_n}{\infty}} = 
\brac{9\max_{1\leq k\leq n}\var\paren{W_k}}^{-1}.
\end{equation}
In the final phase of the proof, we achieve a tight upper bound on $\max_{1\leq k\leq n}\var\paren{W_k}$. It is obvious from the formulation of $A_n$ and the stationarity of $\boldsymbol{X}-\bb{E}\boldsymbol{X}$ that $\max_{1\leq k\leq n}\var\paren{W_k} = n^{-2\nu}\vee \var\paren{W_{p+1}}$. So, the goal is reduced to give an upper bound on the variance of $W_{p+1}$.
\begin{eqnarray}\label{DiagEntriesofDn}
\var\paren{W_{p+1}} &=& \var\paren{\sum\limits_{r=0}^{p} {p \choose r}\paren{-\theta_n}^rX_{p+1-r}}\RelNum{\paren{c}}{=} \frac{1}{2\pi}\int\limits_{\bb{R}}\hat{K}\paren{\omega}\abs{\sum\limits_{r=0}^{p}{p \choose r}\paren{-\theta_n}^r\exp\paren{\frac{-jr\omega}{n}}}^2d\omega\nonumber\\
&=& \frac{1}{2\pi}\int\limits_{\bb{R}}\hat{K}\paren{\omega}\abs{\sum\limits_{r=0}^{p}{p \choose r}\paren{-\exp\paren{\frac{-\paren{1+j\omega}}{n}}}^r}^2d\omega = \frac{1}{2\pi}\int\limits_{\bb{R}}\hat{K}\paren{\omega}\abs{1-e^{\frac{-\paren{1+j\omega}}{n}}}^{2p}d\omega\nonumber\\
&=& \frac{1}{2\pi}\int\limits_{\bb{R}}\hat{K}\paren{\omega}\brac{1+\theta^2_n-2\theta_n\cos\paren{\omega/n}}^{p}d\omega\RelNum{\paren{d}}{\leq} \frac{C_K}{2\pi}\int\limits_{\bb{R}} \frac{\brac{1+\theta^2_n-2\theta_n\cos\paren{\frac{\omega}{n}}}^{p}}{\paren{1+\omega^2}^{\nu+1/2}}d\omega,
\end{eqnarray}
where, identity $\paren{c}$ is implied by \emph{Bochner theorem} (cf. \cite{ML.Stein}, Chapter $2$) and $\paren{d}$ is immediate consequence of Assumption \ref{FxdDomAssu}. Notice that
\begin{equation*}
1+\theta^2_n-2\theta_n\cos\paren{\frac{\omega}{n}}\leq \paren{1-\theta_n}^2+2\theta_n\paren{1-\cos\paren{\frac{\omega}{n}}}\leq \frac{1}{n^2}+2\paren{1-\cos\paren{\frac{\omega}{n}}}= \frac{1}{n^2} + \brac{\frac{\omega}{n}\sinc\paren{\frac{\omega}{2n}}}^2.
\end{equation*}
Let $\xi=p-\paren{\nu+1/2} < 1$. Henceforth, for any $R>0$,
\begin{eqnarray}\label{UppBndInt}
\frac{2\pi n^{2\nu}}{C_K}\var\paren{W_{p+1}}&\leq& n^{2\nu}\int\limits_{\bb{R}} \frac{\set{\frac{1}{n^2} + \brac{\frac{\omega}{n}\sinc\paren{\frac{\omega}{2n}}}^2}^{p}}{\paren{1+\omega^2}^{\nu+1/2}}d\omega = \int\limits_{\bb{R}} \frac{\set{1/n^2 + \brac{\omega\sinc\paren{\frac{\omega}{2}}}^2}^{p}}{\paren{1/n^2+\omega^2}^{\nu+1/2}}d\omega\nonumber\\
&=& \int\limits_{-R}^{R} \frac{\set{1/n^2 + \brac{\omega\sinc\paren{\frac{\omega}{2}}}^2}^{p}}{\paren{1/n^2+\omega^2}^{\nu+1/2}}d\omega + \int\limits_{\abs{\omega}\geq R} \frac{\set{1/n^2 + \brac{\omega\sinc\paren{\frac{\omega}{2}}}^2}^{p}}{\paren{1/n^2+\omega^2}^{\nu+1/2}}d\omega\nonumber\\
&\RelNum{\paren{e}}{\leq}& \int\limits_{-R}^{R} \paren{1/n^2+\omega^2}^{\xi}d\omega + \int\limits_{\abs{\omega}\geq R} \frac{\set{1/n^2 + \brac{\omega\sinc\paren{\frac{\omega}{2}}}^2}^{p}}{\paren{1/n^2+\omega^2}^{\nu+1/2}}d\omega\nonumber\\
&\RelNum{\paren{f}}{\leq}& \int\limits_{-R}^{R} \paren{1/n^2+\omega^2}^{\xi}d\omega + 5^p\int\limits_{\abs{\omega}\geq R} \abs{\omega}^{-\paren{2\nu+1}} d\omega\RelNum{\paren{g}}{\leq} 3R^3 + 5^p \frac{R^{-2\nu}}{\nu}.
\end{eqnarray}
Inequality $\paren{e}$ follows form the fact that $\sup_{\omega\in\bb{R}}\abs{\sinc\paren{\omega/2}}\leq 1$. In order to justify $\paren{f}$, observe that $\abs{\omega\sinc\paren{\omega/2}}\leq 2$ for any $\omega\in\bb{R}$. Thus, for large enough $n$ and $\abs{\omega}\geq R$, we get
\begin{equation*}
\frac{\set{1/n^2 + \brac{\omega\sinc\paren{\frac{\omega}{2}}}^2}^{p}}{\paren{1/n^2+\omega^2}^{\nu+1/2}}\leq \abs{\omega}^{-\paren{2\nu+1}}\paren{1/n^2 + 4}^{p}\leq 5^p\abs{\omega}^{-\paren{2\nu+1}}.
\end{equation*}
Note that there is some $n_0\coloneqq n_0\paren{R,\nu}$ such that $\sup_{\omega\in\bb{R}}\paren{1/n^2+\omega^2}^{\xi}\leq 3/2R^2$ for all $n>n_0$. This immediately entails inequality $\paren{g}$. 

Finally, minimizing the obtained upper bound in \eqref{UppBndInt} over $R>0$, we get
\begin{equation}\label{LowBndVar}
\var\paren{W_{p+1}}\leq CC_Kn^{-2\nu}\paren{1+\frac{1}{\nu}}
\end{equation}
for some universal constant $C>0$. Thus, there is another strictly positive universal constant, $C'$, for which $\max_{1\leq k\leq n}\var\paren{W_k} = n^{-2\nu}\vee \var\paren{W_{p+1}} \leq C'C_Kn^{-2\nu}\paren{1+\frac{1}{\nu}}$. So, \eqref{LowBndQuadTerm2} implies that 
\begin{equation}\label{KeyLowBnd}
\zeta^\top_t\paren{\Sigma_n}^{-1}\zeta_t \gtrsim \frac{n^{2\nu}}{C_K\paren{1+\frac{1}{\nu}}}.
\end{equation}
The combination of \eqref{LowBndExpValH1} and \eqref{KeyLowBnd} completes our proof.
\end{proof}

\begin{proof}[Proof of Theorem \ref{GLRTRateFxdDomGaussCov}]\label{ProofThm2}
The proof proceeds in a similar manner as that of the preceding theorem, in the sense that it is required to show that inequality \eqref{LowBndExpValH1} holds. 
Let $\theta_n = \exp\paren{-\frac{\rho^2}{n^2}}$. $A_n$ represents the inverse of the Cholesky factorization of $\Sigma_n$. For any $k\leq j$ and $q\in\brac{0,1}$, $G\paren{k,j;q}$ denotes the following rational function.
\begin{equation*}
G\paren{k,j;q} = \prod\limits_{\ell=j-k+1}^{j}\paren{1-q^{\ell}}\brac{\prod\limits_{\ell=1}^{k}\paren{1-q^{\ell}}}^{-1},
\end{equation*}
and $G\paren{k,j;1}={j \choose k}$. $G\paren{k,j;q}$ is usually referred to Gaussian binomial coefficients in the combinatorics literature. Finally, let $U_{n,t}\coloneqq A_n\zeta_t$. Similar to \eqref{LowBndExpValH1}, the aim is to obtain a universal lower bound on $\zeta^\top_t\paren{\Sigma_n}^{-1}\zeta_t$ for $t\in\cc{C}_{n,\alpha}$. Observe that, $\zeta^\top_t\paren{\Sigma_n}^{-1}\zeta_t = \LpNorm{U_{n,t}}{2}^2$. 

In order to achieve a tight lower bound on $\LpNorm{U_{n,t}}{2}$, it is pivotal to study the non-asymptotic behaviour of the entries of $A_n$. According to Proposition $1$ of \cite{WL.Loh}, the entries of $A_n$ are given by
\begin{equation*}
\paren{A_n}_{jk} = \paren{-\sqrt{\theta_n}}^{\paren{j-k}}\frac{G\paren{k-1,j-1;\theta_n}}{\sqrt{\prod\limits_{\ell=1}^{j-1}\paren{1-\theta^{\ell}_n}}} \bbM{1}_{\set{j\geq k}}.
\end{equation*}
Since $\frac{\ell\rho^2}{n^2}$ tends to 0 as $n$ gets large for any $\ell\in\set{0,\ldots,n}$ and $\lim\limits_{x\searrow 0} \frac{1-e^{-x}}{x} = 1$, we get
\begin{equation}\label{DenumAsympRel}
\brac{\prod\limits_{\ell=1}^{j-1}\paren{1-\theta^{\ell}_n}}^{-1} = \brac{\prod\limits_{\ell=1}^{j-1}\paren{1-\exp\paren{-\frac{\ell\rho^2}{n^2}}}}^{-1}\asymp \frac{1}{\paren{j-1}!}\paren{\frac{n}{\rho}}^{2\paren{j-1}}.
\end{equation}
Direct calculations show that $G\paren{k-1,j-1;\theta_n}\asymp {j-1 \choose k-1}$ for any $\theta_n$ in a small neighborhood of $1$ and $j,k\in\set{1,\ldots,n}$. Thus,
\begin{equation}\label{NumAsympRel}
\paren{-\sqrt{\theta_n}}^{\paren{j-k}}G\paren{k-1,j-1;\theta_n}\asymp \paren{-1}^{\paren{j-k}}{j-1 \choose k-1}.
\end{equation}
The asymptotic identities \eqref{DenumAsympRel} and \eqref{NumAsympRel} come in handy to analyze $\LpNorm{U_n}{2}^2$:
\begin{eqnarray*}
\LpNorm{U_n}{2}^2 &\geq& \sum\limits_{j=t+1}^{n} \paren{U_n}^2_{j} =  \sum\limits_{j=t+1}^{n}\brac{\sum\limits_{k=t+1}^{j}\paren{A_n}_{jk}-\sum\limits_{k=1}^{t}\paren{A_n}_{jk}}^2\\ &\asymp& \sum\limits_{j=t+1}^{n}\frac{1}{\paren{j-1}!}\paren{\frac{n}{\rho}}^{2\paren{j-1}}\brac{\sum\limits_{k=t+1}^{j}\paren{-1}^{\paren{j-k}}{j-1 \choose k-1}-\sum\limits_{k=1}^{t}\paren{-1}^{\paren{j-k}}{j-1 \choose k-1}}^2\\
&=&\sum\limits_{j=t+1}^{n}\frac{1}{\paren{j-1}!}\paren{\frac{n}{\rho}}^{2\paren{j-1}}\brac{\sum\limits_{k=1}^{j}\paren{-1}^{\paren{j-k}}{j-1 \choose k-1}-2\sum\limits_{k=1}^{t}\paren{-1}^{\paren{j-k}}{j-1 \choose k-1}}^2\\
&=&\sum\limits_{j=t+1}^{n}\frac{1}{\paren{j-1}!}\paren{\frac{n}{\rho}}^{2\paren{j-1}}\brac{0-2\paren{-1}^j{j-1 \choose t}}^2 \asymp \sum\limits_{j=t+1}^{n}\frac{{j-1 \choose t}^2}{\paren{j-1}!}\paren{\frac{n}{\rho}}^{2\paren{j-1}}.
\end{eqnarray*}
Thus, there are universal constants $C,C'>0$ and $C_0$ depending on $\alpha$ and $\rho$ such that
\begin{equation*}
\LpNorm{U_n}{2}^2\geq C\sum\limits_{j=t+1}^{n}\frac{{j-1 \choose t}^2}{\paren{j-1}!}\paren{\frac{n}{\rho}}^{2\paren{j-1}}\geq C\frac{{n-1 \choose t}^2}{\paren{n-1}!}\paren{\frac{n^2}{\rho^2}}^{\paren{n-1}}\RelNum{\paren{a}}{\geq} C'\paren{\frac{n}{t}}^t \frac{1}{\sqrt{n}}\paren{\frac{en^2}{n\rho^2}}^{\paren{n-1}}\RelNum{\paren{b}}{\geq}\paren{C_0n}^n.
\end{equation*}
Note that inequality $\paren{a}$ can be shown using \emph{Stirling's formula} and $\paren{b}$ is obvious implication of the fact that $t\leq \paren{1-\alpha}n$ (Recall $\cc{C}_{n,\alpha}$ from Section \ref{DtctProcedure}). In summary, we have that
\begin{equation*}
\abs{b} \sqrt{\zeta^\top_n\paren{\Sigma_n}^{-1}\zeta_n}\gtrsim \abs{b}\paren{C_0n}^{n/2}.
\end{equation*}
We conclude the proof by appealing to Lemma \ref{ConcIneqNonCentChi2}.
\end{proof}

\begin{proof}[Proof of Theorem \ref{GLRTRateIncDom}]\label{ProofThm3}
The proof is also similar to that of Theorem \ref{GLRTRateFxdDom}. By
applying Lemma \ref{ConcIneqNonCentChi2}, we need to show \eqref{LowBndExpValH1} holds. That is,
\begin{equation}\label{SuffCondIncDom}
\abs{b} \min_{t\in\cc{C}_{n,\alpha}}\sqrt{\zeta^\top_t\paren{\Sigma_n}^{-1}\zeta_t}\geq 8\sqrt{\log\paren{\frac{4n\paren{1-2\alpha}}{\delta}}}.
\end{equation} 
So, the sufficient detectability condition will be obtained by finding a tight uniform lower bound on $\zeta^\top_t\paren{\Sigma_n}^{-1}\zeta_t$ on $\cc{C}_{n,\alpha}$. For any $t\in\cc{C}_{n,\alpha}$, define $a_t,a'_t\in\set{0,1}^n$ by $a_t\paren{i} = \bbM{1}_{\set{i\leq t}}$ and $a'_t=\bbM{1}_n-a_t$. Observe that $a_t,a'_t$ have non-overlapping support and $\zeta_t = a_t-a'_t$. Thus, 
\begin{equation*}
\zeta^\top_t\paren{\Sigma_n}^{-1}\zeta_t= a^\top_t\paren{\Sigma_n}^{-1}a_t + a'^\top_t\paren{\Sigma_n}^{-1}a'_t -2a^\top_t\paren{\Sigma_n}^{-1}a'_t\geq a^\top_t\paren{\Sigma_n}^{-1}a_t-2a^\top_t\paren{\Sigma_n}^{-1}a'_t
\end{equation*}
The last inequality leads to the following key result
\begin{eqnarray}\label{LowBndIncDom}
\min_{t\in\cc{C}_{n,\alpha}}\frac{\zeta^\top_t\paren{\Sigma_n}^{-1}\zeta_t}{n}& \geq & 
\min_{t\in\cc{C}_{n,\alpha}}\frac{a^\top_t\paren{\Sigma_n}^{-1}a_t}{n}-2\max_{t\in\cc{C}_{n,\alpha}}\abs{\frac{a^\top_t\paren{\Sigma_n}^{-1}a'_t}{n}} \RelNum{\paren{a}}{\geq} \min_{t\in\cc{C}_{n,\alpha}}\frac{a^\top_t\paren{\Sigma_n}^{-1}a_t}{n}-\xi_n\nonumber\\
&\RelNum{\paren{b}}{\geq}& \frac{1}{f\paren{0}}-\xi'_n-\xi_n\RelNum{\paren{c}}{\geq} \frac{1}{2f\paren{0}},
\end{eqnarray}
in which $\set{\xi_n}^{\infty}_{n=1}$ and $\set{\xi'_n}^{\infty}_{n=1}$ are appropriately chosen non-negative vanishing sequences, based on the developed results in Appendix \ref{NonAsympInvLargToep}. In \eqref{LowBndIncDom}, inequality $\paren{a}$ and the explicit form $\xi_n$ are obtained from Lemma \ref{UppBndTauSnScn}. One can find a closed form of $\xi'_n$ and verifies inequality $\paren{b}$ using Corollary \ref{TauSigmaInvCor}. Furthermore, $\paren{c}$ holds whenever $n$ is greater than some $n_0$, which depends on $c$ and $\lambda$. The proof of Theorem \ref{GLRTRateIncDom} is completed by combining \eqref{SuffCondIncDom} and \eqref{LowBndIncDom}.
\end{proof}

\subsection{Proofs for Section \ref{DtctRateUnKnwnSpDen}}

\begin{proof}[Proof of Theorem \ref{PlugInGLRTRateFxdDom}]\label{ProofThm4}
Recall $d_{K,\infty}$ from \eqref{DkInfConvto0} and define the event $\cc{A}_m$ by $\cc{A}_m \coloneqq \brac{d_{K,\infty}\paren{\tilde{\eta}_m,\Omega_{\eta}}\leq \rho_m}$. Moreover, we use $\tilde{Z}_t$ to denote $\frac{\zeta^\top_t\paren{\tilde{\Sigma}_n}^{-1}\boldsymbol{X}}{\sqrt{\zeta^\top_t\paren{\tilde{\Sigma}_n}^{-1}\zeta_t}}$. Notice that $\tilde{\eta}_m$, $\tilde{\Sigma}_n$ and are measurable functions on the sample space. Furthermore for any $u\in\cc{A}_m$, $\tilde{\eta}_m\paren{u}$ and $\tilde{\Sigma}_n\paren{u}$ represent the value of these measurable functions at $u$. Lastly, for $u\in\cc{A}_m$ define the random variable $\tilde{Z}_t\paren{u}$ by
\begin{equation*}
\tilde{Z}_t\paren{u} = \frac{\zeta^\top_t\paren{\tilde{\Sigma}_n\paren{u}}^{-1}\boldsymbol{X}}{\sqrt{\zeta^\top_t\paren{\tilde{\Sigma}_n\paren{u}}^{-1}\zeta_t}}.
\end{equation*}
As the range of $\varphi_n$ is $\brac{0,2}$ (See Definition \ref{CDEP}), we have by Assumption \ref{PlgCPDAssu1} that

\begin{equation}\label{barPhiUppBnd}
\varphi_n\paren{\tilde{T}_{GLRT}}\leq 2\tau_m + \bb{E}\paren{\varphi_n\paren{\tilde{T}_{GLRT}}\mid \cc{A}_m}.
\end{equation}
As $\limsup_{m\rightarrow\infty}\rho_m < 1$, there are $\gamma\in\paren{0,1}$ and $m_0\in\bb{N}$ such that $\rho_m\leq \gamma$ for any $m\geq m_0$. Choose $\gamma<\gamma_0<1$ in an arbitrary fashion. We aim to obtain a sufficient condition on $b$ to control the second term in the right hand side of \eqref{barPhiUppBnd} below $\delta$. Notice that throughout the proof we assume that $m\geq m_0$. Conditioning on the occurrence of $\cc{A}_m$, the following statement trivially holds for $m$.

\begin{equation}\label{ModalIneq}
\exists\; \bar{\eta}_m\in\Omega_{\eta}\;\; s.t. \;\; \norm{\frac{\hat{K}\paren{\cdot,\tilde{\eta}_m}}{\hat{K}\paren{\cdot,\bar{\eta}_m}}-1}{\infty}\leq \frac{\rho_m}{\gamma_0}<1.
\end{equation}
Now choose $u\in\cc{A}_m$ arbitrarily. Define the covariance matrix of observations associated to $\bar{\eta}_m\paren{u}$ by $\bar{\Sigma}_n\paren{u} = \brac{K\paren{\frac{r-s}{n},\bar{\eta}_m\paren{u}}}^n_{r,s=1}$. Similar to the proof of Theorem \ref{GLRTRateFxdDom}, it is necessary to study the two following quantities: 1. variance of $\tilde{Z}\paren{u}$ and 2. expected value of $\tilde{Z}\paren{u}$ under the alternative hypothesis, to control the false alarm and miss detection probabilities. Notice that
\begin{eqnarray}\label{SigmatUppBnd}
\sigma_t\paren{u}&\coloneqq& \var \tilde{Z}\paren{u} =  \frac{\zeta^\top_t\paren{\tilde{\Sigma}_n\paren{u}}^{-1}\Sigma_n\paren{\tilde{\Sigma}_n\paren{u}}^{-1}\zeta_t}{\zeta^\top_t\paren{\tilde{\Sigma}_n\paren{u}}^{-1}\zeta_t} \RelNum{\paren{a}}{\leq} B\; \frac{\zeta^\top_t\paren{\tilde{\Sigma}_n\paren{u}}^{-1}\bar{\Sigma}_n\paren{u}\paren{\tilde{\Sigma}_n\paren{u}}^{-1}\zeta_t}{\zeta^\top_t\paren{\tilde{\Sigma}_n\paren{u}}^{-1}\zeta_t}\nonumber\\
&\RelNum{\paren{b}}{\leq}& B\paren{1-\frac{\rho_m}{\gamma_0}}^{-1}.
\end{eqnarray}
Lemma \ref{EqGaussMeasLemma} ensures the existence of some scalar $B\in\paren{1,\infty}$ for which inequality $\paren{a}$ in \eqref{SigmatUppBnd} holds (since $\eta$ and $\bar{\eta}_m\paren{u}$ belong to the same equivalence class). Moreover, $\paren{b}$ can be easily deduced from the combination of \eqref{ModalIneq} and the second inequality in Lemma \ref{VarContIneq}. Namely, $\sigma_t\paren{u}\leq \sigma_0\coloneqq B\paren{1-\frac{\rho_m}{\gamma_0}}^{-1}$. Now, using Lemma \ref{ConcIneqNonCentChi2}, we get

\begin{equation}\label{UppBndType1Err}
\bb{P}\paren{\max_{1\leq t\leq n}\tilde{Z}^2_t\geq R_{n,m} \mid \cc{A}_m}\leq \frac{\delta}{2}.
\end{equation} 
Note that Lemma \ref{VarContIneq} suggests to take $R_{n,\delta,\rho_m} = \sigma_0\brac{1+2\paren{\log\paren{\frac{2n}{\delta}} + \sqrt{\log\paren{\frac{2n}{\delta}}}}}$. So, we have controlled type $1$ error from above in \eqref{UppBndType1Err}. Now we turn to control the type $\RNum{2}$ error from above. Assume that there is a sudden change in the mean of $G$ at $\bar{t}\in\cc{C}_{n,\alpha}$. According to Lemma \ref{ConcIneqNonCentChi2}, type $\RNum{2}$ error is less than $\delta/2$ whenever for any $u\in\cc{A}_m$ 

\begin{equation}\label{eqModal1}
 \frac{\abs{b}}{2}\sqrt{\zeta^\top_{\bar{t}}\paren{\tilde{\Sigma}_n\paren{u}}^{-1}\zeta_{\bar{t}}}\geq 4\sigma_0\sqrt{\log\paren{\frac{2n}{\delta}}} = \frac{4B}{\paren{1-\frac{\rho_m}{\gamma_0}}}\sqrt{\log\paren{\frac{2n}{\delta}}}.
\end{equation}
Applying Lemma \ref{VarContIneq} and then Lemma \ref{EqGaussMeasLemma}, one can easily show that
\begin{equation}\label{KeyLowBnd2}
\zeta^\top_{\bar{t}}\paren{\tilde{\Sigma}_n\paren{u}}^{-1}\zeta_{\bar{t}}\geq \paren{1-\frac{\rho_m}{\gamma_0}}\zeta^\top_{\bar{t}}\paren{\bar{\Sigma}_n\paren{u}}^{-1}\zeta_{\bar{t}}\geq \frac{1-\frac{\rho_m}{\gamma_0}}{B}\zeta^\top_{\bar{t}}\paren{\bar{\Sigma}_n\paren{u}}^{-1}\zeta_{\bar{t}}.
\end{equation}
The combination of last inequality and \eqref{eqModal1} along one line of algebra leads to the following sufficient condition to control type $\RNum{2}$ error holds if 
\begin{equation}\label{KeyLowBnd3}
\abs{b}\min_{t}\sqrt{\zeta^\top_{t}\paren{\Sigma_n}^{-1}\zeta_{t}}\geq 8\paren{\frac{B}{1-\frac{\rho_m}{\gamma_0}}}^{3/2}\sqrt{\log\paren{\frac{2n}{\delta}}},
\end{equation}
We conclude the proof by invoking \eqref{KeyLowBnd} to obtain \eqref{KeyLowBnd3}.
\end{proof}

\subsection{Proofs for Section \ref{DtctCUSUM}}

\begin{proof}[Proof of Theorem \ref{CUSUMFxdDom}]\label{ProofThm5}
Choose $t\in C_{n,\alpha}$ and define
\begin{equation*}
U^{\star}_t\coloneqq \sqrt{\frac{t\paren{n-t}}{n^2}}\paren{\frac{1}{n-t}\sum\limits_{k=t+1}^{n} X_k-\frac{1}{t}\sum\limits_{k=1}^{n} X_k}.
\end{equation*}
Moreover set
\begin{equation*}
R_{n,\delta} = \sqrt{n\paren{1+2\log\paren{\frac{2n\paren{1-2\alpha}}{\delta}}+ 2\sqrt{\log\paren{\frac{2n\paren{1-2\alpha}}{\delta}}} }}.
\end{equation*}
Note that under the null hypothesis, $U^{\star}_t$ is a zero mean random variable and 
\begin{eqnarray*}
\lim\limits_{n\rightarrow\infty}\var\paren{U^{\star}_t}&\RelNum{\paren{a}}{=}&\lim\limits_{n\rightarrow\infty}\frac{t\paren{n-t}}{n^2}\int\limits_{-\infty}^{\infty}\frac{\hat{K}\paren{\omega}}{2\pi} \abs{\frac{1}{n-t}\sum\limits_{k=t+1}^{n}\exp\paren{-jk\omega/n}-\frac{1}{t}\sum\limits_{k=1}^{n}\exp\paren{-jk\omega/n}}^2d\omega\\
&=&\lim\limits_{n\rightarrow\infty}\int\limits_{-\infty}^{\infty}\frac{\hat{K}\paren{\omega}}{2\pi} \abs{\sqrt{\frac{\beta}{1-\beta}}\sum\limits_{k=t+1}^{n}\frac{\exp\paren{-jk\omega/n}}{n}-\sqrt{\frac{1-\beta}{\beta}}\sum\limits_{k=1}^{n}\frac{\exp\paren{-jk\omega/n}}{n}}^2d\omega\\
&\RelNum{\paren{a}}{=}&\int\limits_{-\infty}^{\infty}\frac{\hat{K}\paren{\omega}}{2\pi}\abs{\sqrt{\frac{\beta}{1-\beta}} \int\limits_{\beta}^{1}e^{-j\omega u}du-\sqrt{\frac{1-\beta}{\beta}}\int\limits_{0}^{\beta}e^{-j\omega u}du}^2d\omega = \int\limits_{-\infty}^{\infty}\frac{\hat{K}\paren{\omega}G_{\beta}\paren{\omega}}{2\pi}d\omega,
\end{eqnarray*}
where
\begin{equation}\label{GBetaFunc}
G_{\beta}\paren{\omega}\coloneqq\brac{\paren{1-\beta}\sinc\paren{\frac{\beta\omega}{2}}}^2 
+\brac{\beta\sinc\paren{\frac{\paren{1-\beta}\omega}{2}}}^2\\
+4\beta\paren{1-\beta}\sinc\paren{\frac{\beta\omega}{2}}\sinc\paren{\frac{\paren{1-\beta}\omega}{2}}\sin^2\paren{\frac{\omega}{2}}.
\end{equation}

The identity $\paren{a}$ is implied by Bochner Theorem and $\paren{b}$ follows from the dominated convergence theorem. It is easy to see that $\norm{G_{\beta}}{\infty}\leq 1$ and so $\lim\limits_{n\rightarrow\infty}\var\paren{U^{\star}_n}\leq 1$ by the triangle inequality. Moreover, Lemma \ref{LowBndGBeta} shows that the achieved upper bound on $\sigma^2_n=\var\paren{U^{\star}_n}$ is tight up to some constant whenever $\hat{K}$ has a uniformly bounded derivative. Namely, there is a universal constant $c\in\paren{0,1}$ such that $c\leq \lim\limits_{n\rightarrow\infty}\var\paren{U^{\star}_n} \leq 1$ for any $\beta\in\paren{0,1}$. Let $R^{\star}_{n,\delta} = R^2_{n,\delta}/n$. Thus
\begin{equation}\label{Type1ErrUppbnd}
\bb{P}\paren{T=1\mid \bb{H}_0}=  \bb{P}\paren{\max_{t\in\cc{C}_{n,\alpha}} \abs{U_t}\geq R_{n,\delta}\mid \bb{H}_0} = \bb{P}\paren{ \max_{t\in\cc{C}_{n,\alpha}} \abs{U^{\star}_t}^2\geq R^{\star}_{n,\delta}\mid \bb{H}_0}.
\end{equation}
For any $t\in\cc{C}_{n,\alpha}$, $\abs{U^{\star}_t}^2$ is a (non-normalized) $\chi^2_1$ random variable, as $\sigma^2_n\leq 1$. Moreover $\abs{\cc{C}_{n,\alpha}} = n\paren{1-2\alpha}$. So the part $\paren{a}$ of Lemma \ref{ConcIneqNonCentChi2} says that 
\begin{equation*}
\bb{P}\paren{ \max_{t\in\cc{C}_{n,\alpha}} \abs{U^{\star}_t}^2\geq R^{\star}_{n,\delta}\mid \bb{H}_0} \leq \frac{\delta}{2}.
\end{equation*}
Now we turn to control the miss detection probability. Without loss of generality assume that $b>0$. Choose an arbitrary $t\in\cc{C}_{n,\alpha}$. A line of algebra shows that 
\begin{equation}\label{UnivLowBndExpVal}
\bb{E}\paren{U^{\star}_t\mid \bb{H}_{1,t}}\geq b\sqrt{\alpha\paren{1-\alpha}}.
\end{equation}
Eq. \eqref{SuffCondCUSUMFxdDom} on $b$ implies that $\bb{E}\paren{U^{\star}_t\mid \bb{H}_{1,t}}\geq 4\sqrt{\log\paren{2n\paren{1-2\alpha}/\delta}}$. In other words, given a sudden jump at $t$, $\abs{U^{\star}_s}^2,\;s\in\cc{C}_{n,\alpha}$ are non-central $\chi^2_1$ random variables satisfying the conditions of the part $\paren{b}$ of Lemma \ref{ConcIneqNonCentChi2}. Hence
\begin{equation}\label{Type2ErrUppbnd}
\bb{P}\paren{T = 0\mid \bb{H}_{1,t}}
= \bb{P}\paren{\max_{s\in\cc{C}_{n,\alpha}} \abs{U^{\star}_s}^2 
\leq R^{\star}_{n,\delta} \mid \bb{H}_{1,t} } \leq \frac{\delta}{2}.
\end{equation}
\end{proof}

\begin{proof}[Proof of Theorem \ref{CUSUMIncDom}]\label{ProofThm6}
We continue to use the same notation as the proof of Theorem \ref{CUSUMFxdDom}. Note that there are three appropriately chosen vanishing sequences $\set{a^1_n}_{n\in\bb{N}}$, $\set{a^2_n}_{n\in\bb{N}}$ and $\set{a^3_n}_{n\in\bb{N}}$ such the

\begin{eqnarray}
\var\paren{U_t} &=& \frac{t}{n\paren{n-t}} \var\paren{\sum\limits_{k=t+1}^{n} X_k} + \frac{n-t}{nt} \var\paren{\sum\limits_{k=1}^{t} X_k} -\frac{2}{n}\cov\paren{\sum\limits_{k=1}^{t} X_k,\sum\limits_{k=t+1}^{n} X_k}\nonumber\\
&\RelNum{\paren{a}}{\leq}& \frac{t}{n\paren{n-t}} \var\paren{\sum\limits_{k=t+1}^{n} X_k} + \frac{n-t}{nt} \var\paren{\sum\limits_{k=1}^{t} X_k} + a^1_n\nonumber\\
&\RelNum{\paren{b}}{=}& \frac{t}{n}f\paren{0} + \frac{n-t}{n}f\paren{0} + a^1_n + a^2_n = f\paren{0} + a^3_n,
\end{eqnarray}
in which inequality $\paren{a}$ follows from Lemma \ref{UppBndTauSnScn2} and $\paren{b}$ is implied by identity ($2.1.3$) in \cite{L.Zhengyan}. Thus, there is $n_0\in\bb{N}$ (depending on $f$ and $\vartheta$) such that for any $n\geq n_0$, $\max_{t\in\cc{C}_{n,\alpha}} \var\paren{U_t}\leq \paren{1+\vartheta}f\paren{0}$. The rest of proof will be omitted because of the analogy to \eqref{Type1ErrUppbnd}-\eqref{Type2ErrUppbnd} in the proof of Theorem \ref{CUSUMFxdDom}.
\end{proof}

\subsection{Proofs for Section \ref{LowBndDtctRate}}

\begin{proof}[Proof of Theorem \ref{MinMaxThmFxdDom}]\label{ProofThm7}
We follow the standard method for bounding the Bayes risk from below. Observe that
\begin{eqnarray*}
\inf_{T} \varphi_n\paren{T} &=& 1-\sup_{T}\inf_{t\in\cc{C}_{n,\alpha}} \brac{\bb{P}\paren{T=0\mid \bb{H}_0}-\bb{P}\paren{T=0\mid \bb{H}_{1,t}} }\\
&\geq& 1-\inf_{t\in\cc{C}_{n,\alpha}}\sup_{T} \abs{\bb{P}\paren{T=0\mid \bb{H}_0}-\bb{P}\paren{T=0\mid \bb{H}_{1,t}} } \RelNum{\paren{a}}{\geq} 1-\inf_{t\in\cc{C}_{n,\alpha}} H\paren{\bb{P}_0,\bb{P}_{1,t}},
\end{eqnarray*}
where $\paren{a}$ follows from inequality $2.27$ in \cite{AB.Tsybakov}. So, it suffices to show that $\inf_{t\in\cc{C}_{n,\alpha}} H^2\paren{\bb{P}_0,\bb{P}_{1,t}}\leq\paren{1-\delta}^2$. A few lines of straightforward algebra on the explicit form of Hellinger distance of Gaussian measures indicates that $\inf_{T} \varphi_n\paren{T}\geq \delta$, whenever
\begin{equation}\label{HelDistnBound}
b^2\inf_{t\in\cc{C}_{n,\alpha}} \zeta^\top_t\paren{\Sigma_n}^{-1}\zeta_t \leq 32\log\paren{\frac{1}{\delta\paren{2-\delta}}}.
\end{equation}
Henceforth, it is enough to obtain a tight upper bound on $\inf_{t\in\cc{C}_{n,\alpha}}\zeta^\top_t\paren{\Sigma_n}^{-1}\zeta_t$. 

Let $\sigma = 1$ and choose $d>0$ by $d^{2p-1} = \frac{C'_K\Gamma\paren{p-1/2}}{\sqrt{4\pi}\Gamma\paren{p}}$. Furthermore, let $\hat{F}_{d,p,\sigma}:\bb{R}\mapsto\bb{R}$ denote the Matern spectral density parametrized by $p$ and $d$ and $\sigma$ as \eqref{MaternSpDen}. Note that $d$ is well defined due to the first condition in Assumption \ref{MinMaxAssu}. Define $\xi_t \in\bb{R}^n$ by $\xi_t\paren{k} = \bbM{1}_{k>t}$ and let $\xi'_t = \xi_t - \zeta_t$ for any $t\in\cc{C}_{n,\alpha}$. Moreover, let $\theta_n = \exp\paren{-d/n}$ and $S_t = \set{t+1,\ldots,n}$. Finally, define the covariance matrix $\Psi_n\in\bb{R}^{n\times n}$ by $\Psi_n = \brac{F_{d,p,\sigma}\paren{\paren{r-s}/n}}^n_{r,s=1}$. 
Observe that
\begin{eqnarray}\label{HelDistnBound2}
\zeta^\top_t\paren{\Sigma_n}^{-1}\zeta_t & = & 
2\paren{\xi^\top_t\paren{\Sigma_n}^{-1}\xi_t+\xi'^\top_t\paren{\Sigma_n}^{-1}\xi'_t} - \bbM{1}^\top_n\paren{\Sigma_n}^{-1}\bbM{1}_n \notag \\
& \leq & 4 \paren{\xi^\top_t\paren{\Sigma_n}^{-1}\xi_t\vee \xi'^\top_t\paren{\Sigma_n}^{-1}\xi'_t}.
\end{eqnarray}
We aim to prove that there is a constant $C\coloneqq C\paren{p}>0$ for which $\xi^\top_t\paren{\Sigma_n}^{-1}\xi_t\leq Cn^{2p-1}$. The same upper bound can be obtained for $\xi'^\top_t\paren{\Sigma_n}^{-1}\xi'_t$ in an analogous manner. 

We first show that  
\begin{equation}\label{EquivalenceCond}
\paren{\frac{\hat{K}}{\hat{F}_{d,p,\sigma}}-1}\in\bb{L}^2\paren{\bb{R}}.
\end{equation}
Let $M$ represent the finite $\limsup$ in the second condition of Assumption \ref{MinMaxAssu}. Without loss of generality, we can assume that $\beta<2$ in Assumption \ref{MinMaxAssu}. Using a few lines of algebra 
along with this condition, we get
\begin{eqnarray*}
\limsup\limits_{\omega\rightarrow\infty}\abs{\omega^{\beta}\paren{\frac{\hat{K}\paren{\omega}}{\hat{F}_{d,p}\paren{\omega}}-1}} &=& \limsup\limits_{\omega\rightarrow\infty}\abs{\omega^{\beta}\paren{\frac{\hat{K}\paren{\omega}\omega^{2p}}{C'_K}\paren{1+\frac{d^2}{\omega^2}}^{p}-1}}\leq \limsup\limits_{\omega\rightarrow\infty}\abs{\omega^{\beta}\paren{\frac{\hat{K}\paren{\omega}\omega^{2p}}{C'_K}-1}}\\
&+&\limsup\limits_{\omega\rightarrow\infty}\abs{\omega^{\beta}\brac{\frac{\hat{K}\paren{\omega}\omega^{2p}}{C'_K}\paren{\paren{1+\frac{d^2}{\omega^2}}^{p}-1}}}\\
&\RelNum{\paren{a}}{=}& M+\frac{2pd^2}{C'_K}\limsup\limits_{\omega\rightarrow\infty}\hat{K} \paren{\omega}\abs{\omega}^{2p-2+\beta} \RelNum{\paren{b}}{=} M.
\end{eqnarray*}
Notice that, identity $\paren{a}$ follows from Assumption \ref{MinMaxAssu} and first order Taylor expansion of $\paren{1+x}^p$ for infinitesimal $x>0$. Moreover, $\paren{b}$ follows from the combination of $\beta<2$ and the first condition in Assumption \ref{MinMaxAssu}. Namely, there is $R>0$ such that
\begin{equation*}
\abs{\frac{\hat{K}\paren{\omega}}{\hat{F}_{d,p}\paren{\omega}}-1}\leq \frac{2M}{\abs{\omega}^{\beta}},\quad\forall\;\abs{\omega}\geq R,
\end{equation*}
which substantiates \eqref{EquivalenceCond} as $\beta>1/2$. 

It is known ($4.31$, Chapter $\RNum{3}$, \cite{IA.Ibragimov}) that there is a function $\phi\in\bb{L}^2\paren{\bb{R}}$ with bounded support such that $\hat{F}_{d,p,\sigma}\paren{\omega}\asymp \abs{\hat{\phi}\paren{\omega}}^2$ as $\abs{\omega}\rightarrow\infty$. Theorem $4$ of Skorokhod \cite{AV.Skorokhod} implies that the associated zero mean Gaussian measures to spectral densities $\hat{K}$ and $\hat{F}_{d,p,\sigma}$ are equivalent. Based upon Lemma \ref{EqGaussMeasLemma}, there exists a constant $\ff{H}\in\paren{0,\infty}$ such that
\begin{equation*}
\frac{1}{\ff{H}}\leq\abs{\lim\limits_{n\rightarrow\infty}\frac{\xi^\top_t\paren{\Sigma_n}^{-1}\xi_t}{\xi^\top_t\paren{\Psi_n}^{-1}\xi_t}}\leq \ff{H}.
\end{equation*}
So, it suffices to show that $\xi^\top_t\paren{\Psi_n}^{-1}\xi_t\leq C' n^{2p-1}$ for some appropriately chosen $C'>0$ depending on $\ff{H}$ and $C$. 

Letting $\nu = p-1/2$ and recalling $A_n$, $W$ and $D_n$ form the proof of Theorem \ref{GLRTRateFxdDom}, we have  
\begin{equation}\label{KeyIneqLowBndThm}
\xi^\top_t\paren{\Psi_n}^{-1}\xi_t = \paren{A_n\xi_t}^\top D^{-1}_n\paren{A_n\xi_t}\RelNum{\paren{b}}{\leq} \frac{\LpNorm{A_n\xi_t}{2}^2}{\lambda_{\min}\paren{D_n\paren{S_t,S_t}}}.
\end{equation}
Note that inequality $\paren{b}$ is inferred from $\supp\paren{A_n\xi_t} = S_t$. Applying a similar technique as
\eqref{LowBndL2Norm}, we get
\begin{eqnarray}\label{UppBndNorm2}
\LpNorm{A_n\xi_t}{2}^2 &=& \paren{n-t-p}\paren{1-\theta_n}^p + \sum\limits_{k=1}^{p}\paren{\sum\limits_{j=0}^{k-1}{p \choose j}\paren{-\theta_n}^j}^2\leq n\paren{1-\theta_n}^p +2\sum\limits_{k=1}^{p}\paren{\sum\limits_{j=0}^{k-1}{p \choose j}\paren{-1}^j}^2\nonumber\\
&\leq& 2{2p-2 \choose p-1} + n\paren{1-\theta_n}^p \leq \alpha^pn^{-\paren{p-1}} + 2\paren{2e}^{p-1}\leq \paren{2e}^p.
\end{eqnarray}
So, $\xi^\top_t\paren{\Psi_n}^{-1}\xi_t\leq \paren{2e}^{2p}\brac{\lambda_{\min}\paren{D_n\paren{S_t,S_t}}}^{-1}$. 

Next, we control the smallest eigenvalue of $D_n\paren{S_t,S_t}$ from the below. We first control the diagonal entries from below. Note that all the diagonal entries of $D_n\paren{S_t,S_t}$ are the same and given by (cf. \eqref{DiagEntriesofDn})
\begin{eqnarray}\label{LowBndQ}
Q &=& \int\limits_{\bb{R}}\frac{\hat{F}_{d,p}\paren{\omega}}{2\pi}\brac{1+\theta^2_n-2\theta_n\cos\paren{\omega/n}}^{p}d\omega \RelNum{\paren{c}}{\propto}\int\limits_{\bb{R}}\frac{\paren{d^2+\omega^2}^{-p}}{d^{-2\nu}}\brac{1+\theta^2_n-2\theta_n\cos\paren{\omega/n}}^{p}d\omega\nonumber\\
&=&n^{-2\nu}\int\limits_{\bb{R}}\brac{\frac{\paren{1-\theta_n}^2+4\theta_n\sin^2\paren{d\omega/2}}{1/n^2+\omega^2}}^{p} d\omega\RelNum{\paren{d}}{\geq} \frac{n^{-2\nu}}{2}\int\limits_{\bb{R}} \brac{\sinc\paren{d\omega/2}}^{2p}d\omega = C'_{d}n^{-2\nu},
\end{eqnarray}
where $\paren{c}$ is obtained from \eqref{MaternSpDen} and the inequality $\paren{d}$ 
follows from the fact that for any $\gamma\in\paren{0,1}$ (here we put $\gamma = 2^{-\frac{1}{p}}$), there is $n_0\paren{\gamma}$ such that for any $n\geq n_0\paren{\gamma}$,
\begin{equation*}
\frac{\paren{1-\theta_n}^2+4\theta_n\sin^2\paren{d\omega/2}}{1/n^2+\omega^2}\geq \gamma d^2 \brac{\sinc\paren{d\omega/2}}^2.
\end{equation*} 
The proof of the last inequality will be skipped due to its simplicity. 

Now, let $\Xi \coloneqq D_n\paren{S_t,S_t}/Q$. The combination of \eqref{KeyIneqLowBndThm}, \eqref{UppBndNorm2} and
\eqref{LowBndQ} shows that
\begin{equation*}
\xi^\top_t\paren{\Psi_n}^{-1}\xi_t \leq \frac{C_0n^{-2\nu}}{\lambda_{\min}\paren{\Xi}} \;\;\Rightarrow\;\;\xi^\top_t\paren{\Sigma_n}^{-1}\xi_t \leq \frac{C'_0n^{2\nu}}{\lambda_{\min}\paren{\Xi}} = \frac{C'_0n^{2p-1}}{\lambda_{\min}\paren{\Xi}},
\end{equation*}
for some constants, $C_0\paren{p}$ and $C'_0$ depending on $C_0$, $\ff{H}$ and $K$. It can be shown using identity $1.2$ of \cite{D.Bolin} that there is some integrable function $g:\brac{-\pi,\pi}\mapsto\bb{R}$ with $m_g\coloneqq \essinf\paren{g}>0$ such that $\Xi$ is a $p-$banded correlation matrix, i.e. $\Xi\paren{r,s} = 0$ for $\abs{r-s}\geq p$, and $\Xi = \cc{T}_n\paren{f}$. In remains to note that Lemma $6$ of \cite{RM.Gray} implies that $\lambda_{\min}\paren{\Xi}>m_g$ for any $n$, which concludes the proof.
\end{proof}

\begin{proof}[Proof of Theorem \ref{MinMaxThmIncDom}]
The proof is similar to the proof of Theorem \ref{MinMaxThmFxdDom}, with some minor difference in 
the detail.
Applying the classical technique of bounding $\varphi_n\paren{T}$ from below in terms of Hellinger distance, we need to verify \eqref{HelDistnBound}. Recalling the formulation of $\xi_t$ and $\xi'_t$ from the proof of Theorem \ref{MinMaxThmFxdDom}, we get the following inequality for any $t\in\cc{C}_{n,\alpha}$,

\begin{equation}\label{KeyLowBndIneq}
\zeta^\top_t\paren{\Sigma_n}^{-1}\zeta_t \RelNum{\paren{a}}{\leq} 4 \paren{\xi^\top_t\paren{\Sigma_n}^{-1}\xi_t\vee \xi'^\top_t\paren{\Sigma_n}^{-1}\xi'_t}\RelNum{\paren{b}}{\leq} \frac{4n}{f\paren{0}} + o\paren{n}\RelNum{\paren{c}}{\leq} \frac{4n\paren{1+\vartheta}}{f\paren{0}}.
\end{equation}
Notice that $\paren{a}$ is exactly the same as inequality \eqref{HelDistnBound2}. Moreover, $\paren{b}$ follows from Corollary \ref{TauSigmaInvCor} and lastly, there is $n_0\in\bb{N}$, which depends on $f$ and $\vartheta$, for which inequality $\paren{c}$ holds. Using \eqref{KeyLowBndIneq}, one can easily verify inequality \eqref{HelDistnBound} and concluding the proof.
\end{proof}

\section{Auxiliary results}\label{AuxRes}
This section contains several technical results needed in Appendix \ref{Proofs}.

\begin{lem}\label{ConcIneqNonCentChi2}
Let $\sigma_0 \geq 1$ and $n\geq 2$. Let $\ff{Z}\in\bb{R}^{n}$ be a Gaussian random vector with $\bb{E}\ff{Z} = \mu$ and $\var \ff{Z}_k\leq \sigma^2_0$ for any $1\leq k\leq n$. Moreover, let $R_n = 1+2\paren{\log\paren{\frac{2n}{\delta}} + \sqrt{\log\paren{\frac{2n}{\delta}}}}$. For any $\delta\in\paren{0,1}$ and any $n\in\bb{N}$, the following results hold.
\begin{enumerate}
\item If $\mu = 0$, then $\bb{P}\brac{\max\limits_{1\leq j\leq n}\ff{Z}^2_j \geq \sigma^2_0R_n}\leq \frac{\delta}{2}$.
\item If $\max\limits_{1\leq j\leq n} \abs{\mu_j}\geq 4\sigma_0\sqrt{\log\paren{\frac{2n}{\delta}}}$, then $\bb{P}\brac{\max\limits_{1\leq j\leq n}\ff{Z}^2_j \leq \sigma^2_0R_n}\leq \frac{\delta}{2}$.
\end{enumerate}
\end{lem}
\begin{proof}
For brevity, let $\sigma_j = \var \ff{Z}_j,\;j=1,\ldots,n$. Notice that $\paren{\frac{\ff{Z}_j}{\sigma_j}}^2$ are standard $\chi^2_1$ random variables, for any $j=1,\ldots,n$. Lemma $8.1$ in \cite{L.Birge} implies that $\bb{P}\paren{\ff{Z}^2_j\geq \sigma^2_jR_n}\leq \frac{\delta}{2n}$. Thus, $\bb{P}\paren{\ff{Z}^2_j\geq \sigma^2_0R_n}\leq \frac{\delta}{2n}$ due to $\sigma_j\leq \sigma_0$. We conclude the proof of the first part by a union bound argument. Now, we turn to prove the second part. Define $k \coloneqq \arg\max\limits_{1\leq j\leq n} \abs{\mu_j}$. It is easy to verify that $R_n\leq 4\log\paren{\frac{2n}{\delta}}$. Observe that
\begin{equation*}
\bb{P}\brac{\max\limits_{1\leq j\leq n}\ff{Z}^2_j \leq \sigma^2_0R_n} \leq \bb{P}\brac{\frac{\ff{Z}^2_k}{\sigma^2_k} \leq \paren{\frac{\sigma_0}{\sigma_k}}^2R_n}\leq \bb{P}\brac{\frac{\ff{Z}^2_k}{\sigma^2_k} \leq 4\paren{\frac{\sigma_0}{\sigma_k}}^2\log\paren{\frac{2n}{\delta}}}.
\end{equation*}
Moreover, $\frac{\ff{Z}^2_k}{\sigma^2_k}$ is a non-central $\chi^2_1$ random variables with non-centrality parameter $B_k\coloneqq\abs{\frac{\mu_k}{\sigma_k}}$. The lower bound condition on $\abs{\mu_k}$ implies that $B_k\geq 4\frac{\sigma_0}{\sigma_k}\sqrt{\log\paren{\frac{2n}{\delta}}}$. We finish the proof by the following inequality,

\begin{equation*}
\bb{P}\brac{\frac{\ff{Z}^2_k}{\sigma^2_k} \leq 4\paren{\frac{\sigma_0}{\sigma_k}}^2\log\paren{\frac{2n}{\delta}}}\RelNum{\paren{a}}{\leq}\bb{P}\brac{\frac{\ff{Z}^2_k}{\sigma^2_k}\leq 1+B^2_k - 2\sqrt{\paren{1+2B^2_k}\log\paren{\frac{2}{\delta}}}}\RelNum{\paren{b}}{\leq} \frac{\delta}{2}.
\end{equation*}
In order to demonstrate inequality $\paren{a}$, we need to show that $1+B^2_k - 2\sqrt{\paren{1+2B^2_k}\log\paren{\frac{2}{\delta}}}\geq 4\paren{\frac{\sigma_0}{\sigma_k}}^2\log\paren{\frac{2n}{\delta}}$ which can be shown by obvious inequality $\sigma_0/\sigma_k\geq 1$ and a few lines of algebra. Inequality $\paren{b}$ can be inferred from Lemma $8.1$ of \cite{L.Birge}.
\end{proof}

\begin{prop}[\emph{Kantorovich} inequality, (p. $452$, \cite{RA.Horn})]\label{KanIneq}
Let $\Sigma\in\bb{R}^{n\times n}$ be a non-singular covariance matrix and let $V\in\bb{R}^n$ be a non-zero vector. Then, $V^\top \Sigma^{-1}V \geq \frac{\LpNorm{V}{2}^4}{V^\top \Sigma V}$.
\end{prop}

\begin{lem}\label{VarContIneq}
Let $K_0$ and $K_1$ be two covariance function with spectral densities $\hat{K}_0$ and $\hat{K}_1$, respectively. Define, $\Sigma_0 \coloneqq \brac{K_0\paren{\frac{r-s}{n}}}^n_{r,s=1}$ and $\Sigma_1 \coloneqq \brac{K_1\paren{\frac{r-s}{n}}}^n_{r,s=1}$. Suppose that there exists $\rho\in\paren{0,1}$ such that $\norm{\frac{K_1}{K_0}-1}{\infty}\leq \rho$. Then, 
\begin{enumerate}[label=(\alph*),leftmargin=*]
\item $\Sigma_0-\Sigma_1\preceq\frac{\rho}{1-\rho} \Sigma_1$.
\item $\Sigma^{-1}_1 \paren{\Sigma_0-\Sigma_1}\Sigma^{-1}_1\preceq\frac{\rho}{1-\rho} \Sigma^{-1}_1$.
\end{enumerate}
\end{lem}

\begin{proof}
Trivial calculations on $\norm{\frac{K_1}{K_0}-1}{\infty}\leq \rho$ shows that for any $\omega\in\bb{R}$,
\begin{equation}\label{UppLowBndDiffSpDens}
\frac{-\rho}{1+\rho} \hat{K}_1\paren{\omega}\leq \hat{K}_0\paren{\omega} - \hat{K}_1\paren{\omega}\leq\frac{\rho}{1-\rho} \hat{K}_1\paren{\omega}.
\end{equation}
Choose $v\in\bb{R}^n$ arbitrarily. The basic properties of spectral density and inequality \eqref{UppLowBndDiffSpDens} imply that

\begin{eqnarray*}
2\pi v^\top\paren{\Sigma_0-\Sigma_1}v &=& \int\limits_{\bb{R}}\paren{\hat{K}_0\paren{\omega} - \hat{K}_1\paren{\omega}}\abs{\sum\limits_{\ell=1}^{n}v_{\ell}e^{j\omega\ell/n}}^2d\omega \leq \frac{\rho}{\paren{1+\rho}}\int\limits_{\bb{R}}\hat{K}_1\paren{\omega}\abs{\sum\limits_{\ell=1}^{n}v_{\ell}e^{j\omega\ell/n}}^2d\omega\\
&=&\frac{2\pi\rho}{1-\rho}v^\top \Sigma_1 v.
\end{eqnarray*}
Thus, $\Sigma_0-\Sigma_1\preceq \frac{\rho}{1-\rho} \Sigma_1$. The second inequality is an obvious implication of the first inequality.
\end{proof} 

\begin{lem}\label{PowExpSpDenDecay}
Let $\delta\in\paren{0,2}$, $d\in\paren{0,\infty}$ and define $K:\bb{R}\mapsto\bb{R}$ by $K\paren{r} = \sigma^2\exp\paren{-\abs{\frac{r}{d}}^{\delta}}$. Then,
\begin{equation*}
\lim\limits_{\omega\rightarrow\infty} \hat{K}\paren{\omega}\abs{\omega}^{1+\delta} = C_{\delta}\paren{d,\sigma}\coloneqq\frac{\sigma^2\delta\Gamma\paren{\delta}\sin\paren{\frac{\pi\delta}{2}}}{\pi d^{\delta}}.
\end{equation*}
\end{lem}

\begin{proof}
It is obvious that $C_{\delta}\paren{d,\sigma} = \sigma^2 C_{\delta}\paren{d,1}$, so without loss of generality assume that $\sigma = 1$. It is trivial that $K\paren{r}$ is of index $\delta$ as $\abs{r}\rightarrow 0$, i.e. $\lim\limits_{\abs{r}\rightarrow 0} \frac{1-K\paren{r\lambda}}{1-K\paren{r}} = \lambda^{\delta}$ for any $\lambda>0$. Thus, applying the Tauberian Theorem (p. 35, \cite{ML.Stein}) leads to
\begin{equation}\label{TauberianThm}
\lim\limits_{\omega\rightarrow\infty} \brac{1-K\paren{1/\omega}}^{-1} \int\limits_{\omega}^{\infty} \hat{K}\paren{u}du=  \frac{\Gamma\paren{\delta}\sin\paren{\frac{\pi\delta}{2}}}{\pi} = \frac{C_{\delta}\paren{d,1} d^{\delta}}{\delta}.
\end{equation}
Moreover, the first order Taylor expansion of $e^{-x}$ is at $0$, implies that $\brac{1-K\paren{1/\omega}}^{-1}\paren{d\omega}^{-\delta}\rightarrow 1$ as $\omega\rightarrow 0$. Thus, \eqref{TauberianThm} can be rewritten by last limiting identity and applying L'Hospital's rule.
\begin{equation*}
C_{\delta}\paren{d,1} = \lim\limits_{\omega\rightarrow\infty}\delta d^{-\delta} \paren{d\omega}^{\delta} \delta\omega^{\delta}\int\limits_{\omega}^{\infty} \hat{K}\paren{u}du = \lim\limits_{\omega\rightarrow\infty}\delta\omega^{\delta} \int\limits_{\omega}^{\infty} \hat{K}\paren{u}du =  \lim\limits_{\omega\rightarrow\infty} \hat{K}\paren{\omega}\abs{\omega}^{1+\delta}.
\end{equation*}
\end{proof}

The following Lemma is probably well-known in the literature of Gaussian processes (e.g. the identity $2$ of \cite{ML.Stein} (p.112) is analogous but not exactly same as the part $\paren{a}$ of Lemma \ref{EqGaussMeasLemma}). Because of the absence of direct references, we include and prove the following result in this section.

\begin{lem}\label{EqGaussMeasLemma}
Let $G_i,\;i=1,2$ be two zero mean stationary Gaussian process in $\brac{0,1}$ associated to covariance functions $K_i,\;i=1,2$, respectively. For any $n\in\bb{N}$, define two positive definite covariance matrices by $\Sigma_n \coloneqq \brac{K_1\paren{\frac{r-s}{n}}}$ and $\Psi_n \coloneqq \brac{K_2\paren{\frac{r-s}{n}}}$. If $G_1$ and $G_2$ induce equivalent measures on the Hilbert space of $\bb{L}^2\paren{\brac{0,1}}$, then there exists an scalar $B\in\brapar{1,\infty}$ for which
\begin{enumerate}
\item $\frac{1}{B}\leq \lim\limits_{n\rightarrow\infty} \inf_{v\ne \zero_n} \frac{v^\top \Sigma_n v}{v^\top\Psi_n v}\leq \lim\limits_{n\rightarrow\infty} \sup_{v\ne \zero_n} \frac{v^\top \Sigma_n v}{v^\top\Psi_n v}\leq B$.
\item $\frac{1}{B}\leq \lim\limits_{n\rightarrow\infty} \inf_{v\ne \zero_n} \frac{v^\top \Sigma^{-1}_n v}{v^\top\Psi^{-1}_n v}\leq \lim\limits_{n\rightarrow\infty} \sup_{v\ne \zero_n} \frac{v^\top \Sigma^{-1}_n v}{v^\top\Psi^{-1}_n v}\leq B$.
\end{enumerate}
\end{lem}

\begin{proof}
We use $\bb{P}_i,i=1,2$ to denote the probability measures with respect to $G_i,\;i=1,2$, respectively. Abusing the notation, $\boldsymbol{X}\in\bb{R}^n$ represents the random vector generated by sampling Gaussian process at $\set{k/n}^n_{k=1}$ for any $n\in\bb{N}$. We prove the existence of a finite scalar $B_1$ for which $\lim\limits_{n\rightarrow\infty} \sup_{v\ne \zero_n} \frac{v^\top \Sigma_n v}{v^\top\Psi_n v}\leq B_1$. Assume toward contradiction that $\lim\limits_{n\rightarrow\infty} \sup_{v\ne \zero_n} \frac{v^\top \Sigma_n v}{v^\top\Psi_n v}$ tends to infinity. So, there is a sequence of non-zero vectors $\set{v_n\in\bb{R}^n}^{\infty}_{n=1}$ such that 
\begin{equation}\label{ContradictionPremise}
\limsup\limits_{n\rightarrow\infty} \frac{v_n^\top \Sigma_n v_n}{v_n^\top\Psi_n v_n} = \infty.
\end{equation} 
Consider the measurable event ${\rm E}_n = \brac{\abs{\InnerProd{v_n}{\boldsymbol{X}}}\geq \sqrt{v_n^\top \Sigma_n v_n}}$. Simple calculations shows that
\begin{equation}\label{PofAn}
\bb{P}_1\paren{{\rm E}_n} = Q\paren{1},\quad \bb{P}_2\paren{{\rm E}_n} = Q\paren{\sqrt{\frac{v_n^\top \Sigma_n v_n}{v_n^\top\Psi_n v_n}}},
\end{equation} 
in which $Q\paren{\cdot}$ stands for the $Q$-function, i.e. $Q\paren{r} = \int\limits \frac{1}{\sqrt{2\pi}} \exp\paren{-x^2/2}\bbM{1}\paren{\abs{x}\geq r}dx$. Combining \eqref{ContradictionPremise} and \eqref{PofAn} leads to $\limsup\limits_{n\rightarrow\infty} \frac{\bb{P}_1\paren{{\rm E}_n}}{\bb{P}_2\paren{{\rm E}_n}} = \infty$ which contradicts the absolute continuity of $\bb{P}_1$ with respect to $\bb{P}_2$. One cam show using the same technique that there is $B_2\in\paren{1,\infty}$ such that 
\begin{equation*}
\frac{1}{B_2}\leq \lim\limits_{n\rightarrow\infty} \inf_{v\ne \zero_n} \frac{v^\top \Sigma_n v}{v^\top\Psi_n v}.
\end{equation*}
We conclude the proof by choosing $B = B_1 \vee B_2$. Now, we turn to substantiate the second claim. Pick a non-zero vector $v\in\bb{R}^n$. According to Lemma \ref{ConvertInvCovToCovLem}, there is an suitably chosen $n-$dimensional vector $u$ (The inner product of $u$ and $v$ is necessarily $1$) such that
\begin{equation*}
\frac{v^\top\Sigma^{-1}_n v}{v^\top\Psi^{-1}_n v} = v^\top\Sigma^{-1}_n v u^\top\Psi_n u = \frac{u^\top\Psi_n u}{\max_{\InnerProd{\omega}{v} = 1} \omega^\top\Sigma_n\omega}\leq \frac{u^\top\Psi_n u}{u^\top\Sigma_n u}\RelNum{\paren{a}}{\leq} B.
\end{equation*}
Note that the inequality $\paren{a}$ is obtained from the first part of this Lemma. Taking supremum over all non-zero $v\in\bb{R}^n$ and $n\in\bb{N}$ terminates the proof. 
\end{proof}

\begin{lem}\label{ConvertInvCovToCovLem}
Let $\Sigma\in\bb{R}^{n\times n}$ be a non-singular covariance matrix and let $\omega\in\bb{R}^n$ be a non-zero vector. Then,
\begin{equation}\label{MinofConvProg}
\paren{\omega^\top \Sigma^{-1}\omega}^{-1} = \min _{\InnerProd{v}{\omega}{} = 1} v^\top \Sigma v.
\end{equation}
\end{lem}

\begin{proof}
Since the optimization problem in \eqref{MinofConvProg} is a convex program with continuously differentiable objective function and constraint, so its minimal value can be obtained solving the KKT equations. That is, there are $\hat{\lambda}\geq 0$ and $\hat{v}$ such that
\begin{equation*}
2\Sigma \hat{v} - \hat{\lambda}\omega = 0,\quad\hat{\lambda}\paren{\InnerProd{\hat{v}}{\omega}{}-1} = 0.
\end{equation*}
Solving the above set of equations yields, $\hat{v} = \frac{\Sigma^{-1}\omega}{\omega^\top \Sigma^{-1}\omega}$. The desired result will be established by replacing $\hat{v}$ into the right hand side of \eqref{MinofConvProg}.
\end{proof}

\begin{lem}\label{VerifofReqCondsAGLRT}
Let $\Omega = \set{\eta=\paren{d,\sigma}:\eta\in\Omega}\subset\paren{0,\infty}^2$ be a compact set such that 
\begin{equation}\label{DistOmegaXYAxes}
\dist\paren{\Omega, \set{\paren{x,y}:\;x=0\;\mbox{or}\; y=0}} > 0.
\end{equation}
The conditions of Proposition \ref{PropFxdPluginGLRT} are satisfied for the following scenarios.
\begin{enumerate}
\item $K\paren{\cdot,\eta}$ is the powered exponential covariance function with known $\beta\in\paren{0,2}$ and $\eta\in\Omega$. 
\item $\hat{K}\paren{\cdot,\eta}$ is the Matern spectral density with known $\nu$, given by \eqref{MaternSpDen}.
\end{enumerate}
\end{lem}

\begin{proof}
We first substantiate part $1$. It is easy to verify the first condition in Proposition \ref{PropFxdPluginGLRT} for powered exponential covariance. For proving the continuity condition, let $\eta_m$ be a convergent sequence in $\Omega$ to $\eta$. So
\begin{equation*}
\norm{\hat{K}\paren{\cdot,\eta_m}-\hat{K}\paren{\cdot,\eta}}{\infty}\leq\int\limits_{\bb{R}}\abs{K\paren{r,\eta_m}-K\paren{r,\eta}} dr \rightarrow 0,\quad\mbox{as}\;m\rightarrow\infty
\end{equation*}
due to the dominated convergence Theorem. Finally, Lemma \ref{PowExpSpDenDecay} and a few lines of algebra imply that
\begin{equation*}
\sup_{\eta\in\Omega}\LpNorm{\nabla\log \ff{C}_{K,\Omega}\paren{\eta}}{2} = 2\sup_{\paren{\sigma,d}\in\Omega} \sqrt{\sigma^{-2}+\paren{\beta/d}^2} < \infty.
\end{equation*} 
Notice that the last inequality is a consequence of \eqref{DistOmegaXYAxes}. Now we turn to the proof of second part. The verification of the first two conditions in Proposition \ref{PropFxdPluginGLRT} is analogous to the powered exponential covariance function. Note that $\ff{C}_{K,\Omega}\paren{\eta} = \frac{\sigma^2\sqrt{4\pi}\Gamma\paren{\nu+1/2}}{\Gamma\paren{\nu}}\rho^{-2\nu}$ for Matern covariance. Thus, 
\begin{equation*}
\sup_{\eta\in\Omega}\LpNorm{\nabla\log \ff{C}_{K,\Omega}\paren{\eta}}{2} = 2\sup_{\paren{\sigma,d}\in\Omega} \sqrt{\sigma^{-2}+\paren{\nu/d}^2} < \infty.
\end{equation*}
\end{proof}

\begin{lem}\label{ContL2NormConv}
Suppose that $\hat{K}\paren{\omega,\eta}$ satisfies the conditions in the part $2$ of Lemma \ref{VerifofReqCondsAGLRT} and the taper function $\hat{f}_{\tap}\paren{\omega}$ admits Assumption \ref{PlgCPDAssu2}. Let $\hat{K}_{\tap}\paren{\cdot,\eta}$ be the convolution of $\hat{K}\paren{\omega,\eta}$ and $\hat{f}_{\tap}$. Then,
\begin{equation*}
g\paren{\eta} = \norm{\frac{\hat{K}_{\tap}\paren{\cdot,\eta}}{\hat{K}\paren{\cdot,\eta}}-1}{2},
\end{equation*}
is a continuous function of $\eta$ in $\Omega$.
\end{lem}

\begin{proof}
Let $d_{\max} = \sup\set{d:\;\paren{\sigma,d}\in\Omega}$ and $d_{\min} = \inf\set{d:\;\paren{\sigma,d}\in\Omega}$. Note that $0<d_{\min}\leq d_{\max} < \infty$ as $\Omega$ is a compact set in the interior of the upper right half-plane. define $h:\bb{R}\times\Omega\mapsto\paren{0,\infty}$ by 
\begin{equation*}
h\paren{\omega,\eta} \coloneqq \frac{\hat{K}_{\tap}\paren{\omega,\eta}}{\hat{K}\paren{\omega,\eta}}-1 = \int\limits_{-\infty}^{\infty} \hat{f}_{\tap}\paren{\omega-u}\frac{\hat{K}\paren{u,\eta}}{\hat{K}\paren{\omega,\eta}}du - 1.
\end{equation*}
We first show that for any fixed $\omega\in\bb{R}$, $h\paren{\omega,\eta}$ is a continuous function of $\eta$. Define two real valued function $h_{d,\omega}\paren{u} = \paren{\frac{d^2_{\max}+\omega^2}{d^2_{\min}+u^2}}^{\nu+1/2}$ and $p\paren{u,\eta}\coloneqq\hat{f}_{\tap}\paren{\omega-u}\frac{\hat{K}\paren{u,\eta}}{\hat{K}\paren{\omega,\eta}}$. Note that $h_{d,\omega}$ is integrable and for any $u\in\bb{R}$, $p\paren{u,\cdot}$ is a continuous function of $\eta$. Moreover,
\begin{equation*}
\sup_{\eta=\paren{\sigma,d}\in\Omega} p\paren{u,\eta} = \hat{f}_{\tap}\paren{\omega-u} \sup_{\eta\in\Omega}\paren{\frac{d^2+\omega^2}{d^2+u^2}}^{\nu+1/2}\leq \sup_{\eta\in\Omega}\paren{\frac{d^2+\omega^2}{d^2+u^2}}^{\nu+1/2}\leq h_{d,\omega}\paren{u}.
\end{equation*}
Choose $\eta\in\Omega$ and a convergent sequence $\set{\eta_n}^{\infty}_{n=1}\subset\Omega
$ to $\eta$, in an arbitrary fashion. The dominated convergence theorem shows that
\begin{equation*}
\lim\limits_{n\rightarrow\infty}h\paren{\omega,\eta_n} = \lim\limits_{n\rightarrow\infty}\int\limits_{-\infty}^{\infty} p\paren{u,\eta_n} du = h\paren{\omega,\eta},
\end{equation*}
which confirms the continuity of $h$ with respect to $\eta$. Now, we demonstrate the continuity of $g$. Notice that $g\paren{\eta} = \norm{h\paren{\cdot,\eta}}{2}$. One can show using the exact same technique as Lemma $4$ of \cite{J.Du} (with a minor algebraic difference) that there is a function $q:\bb{R}\mapsto\bb{R}$ given by $q\paren{\omega} = c_1 \vee c_2\abs{\omega}^{-\paren{1+r}}$, where $c_1$, $c_2$ and $r$ are appropriately chosen finite positive scalars, uniformly dominating $h\paren{\omega,\eta}$ over $\Omega$. Namely,
\begin{equation*}
\sup_{\eta\in\Omega} h\paren{\omega,\eta}\leq q\paren{\omega}.
\end{equation*}
Thus,
\begin{equation*}
\abs{\lim\limits_{n\rightarrow\infty} g\paren{\eta_n}-g\paren{\eta}}\leq \lim\limits_{n\rightarrow\infty} \norm{h\paren{\cdot,\eta}-h\paren{\cdot,\eta_n}}{2} =  \paren{\int\limits_{-\infty}^{\infty}\brac{h^2\paren{\omega,\eta}-h^2\paren{\omega,\eta_n}}d\omega}^{1/2} \RelNum{\paren{a}}{=} 0,
\end{equation*} 
terminating the proof. Notice that identity $\paren{a}$ follows from the continuity of $h$ with respect to $\eta$ and the dominated convergence theorem.
\end{proof}

\begin{lem}\label{LowBndGBeta}
Let $K$ be a covariance function such that $\norm{\hat{K}'}{\infty}<\infty$ and define $G_{\beta}:\bb{R}\mapsto\brac{0,1}$ by \eqref{GBetaFunc}. Then, there is a universal constant $c>0$ such that
\begin{equation*}
\inf_{\beta\in\paren{0,1}} \int\limits_{-\infty}^{\infty}\hat{K}\paren{\omega}G_{\beta}\paren{\omega}d\omega\geq c.
\end{equation*}
\end{lem}

\begin{proof}
Observe that for any $\omega\in\bb{R}$, $G_{\beta}\paren{\omega}$ is a quadratic function of $\beta$ in the compact interval $\brac{0,1}$ and $\lim\limits_{n\rightarrow\infty}\norm{G_{\beta_n}-G_{\beta}}{\infty}=0$ for any convergent sequence $\beta_n\rightarrow\beta$. This property implies that 
\begin{equation}\label{ModalIneq1}
\inf_{\beta\in\paren{0,1}} \int\limits_{-\infty}^{\infty}\hat{K}\paren{\omega}G_{\beta}\paren{\omega}d\omega\geq \frac{1}{2}\brac{\inf_{\beta\in\paren{0,1},\; \abs{\beta-1/2}\geq r} \int\limits_{-\infty}^{\infty}\hat{K}\paren{\omega}G_{\beta}\paren{\omega}d\omega \wedge \int\limits_{-\infty}^{\infty}\hat{K}\paren{\omega}G_{0.5}\paren{\omega}d\omega}.
\end{equation}
for some sufficiently small $r>0$. Observe that, $G_{\beta}\paren{0} = \paren{1-2\beta}^2>0$ for $\beta\ne 1/2$. The differentiability of $G_{\beta}$ and $\hat{K}\paren{\omega}$ implies the existence of a non-degenerate open interval $\cc{I}_{\beta}$ centered at $0$ such that,
\begin{equation*}
\inf_{\omega\in \cc{I}_{\beta}} \hat{K}\paren{\omega}G_{\beta}\paren{\omega}\geq \frac{\paren{1-2\beta}^2\hat{K}\paren{0}}{2}\;\;\Rightarrow\;\;\int\limits_{-\infty}^{\infty}\frac{\hat{K}\paren{\omega}G_{\beta}\paren{\omega}}{2\pi}d\omega\geq\frac{\paren{1-2\beta}^2\hat{K}\paren{0}}{4\pi}\abs{\cc{I}_{\beta}}.
\end{equation*}
Notice that $\inf_{\abs{\beta-1/2}\geq r}\paren{1-2\beta}^2\abs{\cc{I}_{\beta}}>0$. So, we just need to show that the corresponding term to $\beta=1/2$ in the right hand side of \eqref{ModalIneq1} is strictly positive. For $\beta = 1/2$, $G_{\beta}\paren{\omega} = \brac{\sinc\paren{\omega/4}\sin\paren{\omega/2}}^2$ and so
\begin{equation*}
\int\limits_{-\infty}^{\infty}\frac{\hat{K}\paren{\omega}G_{\beta}\paren{\omega}}{2\pi}d\omega \geq \int\limits_{-2\pi}^{2\pi}\frac{\hat{K}\paren{\omega}\brac{\sinc\paren{\omega/4}\sin\paren{\omega/2}}^2}{2\pi}d\omega \RelNum{\paren{b}}{\geq} \frac{2}{\pi^3}\int\limits_{-2\pi}^{2\pi}\hat{K}\paren{\omega}\sin^2\paren{\omega/2}d\omega\RelNum{\paren{c}}{>}0.
\end{equation*}
Note that $\paren{b}$ is a consequence of monotonicity of $\sinc\paren{\cdot}$ in the interval $\paren{0,\pi/2}$ and inequality $\paren{c}$ follows from the combination of $\abs{\hat{K}'\paren{0}}<\infty$ and $\hat{K}\paren{0}>0$.
\end{proof}

\section{Non-asymptotic behaviour of the inverse of Toeplitz matrices}\label{NonAsympInvLargToep}

In this section, we investigate some non-asymptotic properties of the inverse of Toeplitz matrices with polynomially decaying off-diagonal entries. The developed results plays a crucial role in the analysis of GLRT in increasing domain. We first introduced some simplifying notation. For any symmetric and periodic function $f\in\bb{L}^{\infty}\paren{\bb{R}}$ with period $2\pi$, define 
\begin{equation*}
\cc{T}_{\bb{N}}\paren{f} \coloneqq \paren{\begin{array}{cccc}
f_0 & f_1 & f_2 & \cdots \\
f_1 & f_0 & f_1 & \ddots\\
f_2 & f_1 & f_0 & \ddots \\
\vdots & \ddots & \ddots & \ddots \\
\end{array}}.
\end{equation*}
in which $\set{f_m}_{m\in\bb{Z}}$ denotes the set of Fourier coefficients of $f$ and $f_{-m} = f_m$ for any $m\in\bb{Z}$. Moreover, let $\cc{T}_n\paren{f} = \brac{\cc{T}_{\bb{N}}\paren{f}}^n_{r,s=1}$. Finally, for any $S_n,S'_n\subseteq\set{1,\ldots,n}$ and any $L_n\in\bb{R}^{n\times n}$ define
\begin{equation*}
\cc{\tau}\paren{L_n, S_n, S'_n} \coloneqq \frac{\bbM{1}^\top_{S_n} L_n \bbM{1}_{S'_n}}{\sqrt{\abs{S_n}\;\abs{S'_n}}}, 
\end{equation*}
in which $\bbM{1}_{S_n}\in\bb{R}^n$ represents the indicator vector of $S_n$ in $\set{1,\ldots,n}$. For the sake of brevity, we use the shorthand notation $\cc{\tau}\paren{L_n, S_n}$ if $S_n = S'_n$. Throughout this section we assume that there exists some $\alpha\in\paren{0,1/2}$ such that $\cc{C}_{n,\alpha} = \brac{\alpha n, \paren{1-\alpha}n}$ and $L_n$ is a Toeplitz covariance matrix where its generator satisfies Assumption \ref{AssumpToepCov} for some $c,\lambda,m_f$ and $M_f$. We also use $L_{\bb{N}}$ to indicate the infinite sized version of $L_n$.

\begin{lem}\label{UppBndTauSnScn}
Let $S_n = \set{1,\ldots,t}$ for some $t\in\cc{C}_{n,\alpha}$ and $S^c_n = \set{1,\ldots,n}\setminus S_n$. There is a bounded constant $C\geq 0$ depending on $f$ and $\alpha$ such that
\begin{equation*}
\abs{\cc{\tau}\paren{L^{-1}_n, S_n, S^c_n}} \;\leq \frac{C}{\sqrt{\alpha\paren{1-\alpha}}}n^{-\paren{\lambda\wedge 1}} \paren{\bbM{1}_{\set{\lambda\ne 1}}+ \bbM{1}_{\set{\lambda=1}}\log n}.
\end{equation*}
\end{lem}

\begin{proof}
Without loss of generality assume that $t\leq n-t$. The proof is based upon the well known fact that for a positive definite matrix with polynomially decaying off-diagonal entries, the corresponding elements of its inverse shrink with same rate (See e.g. \cite{S.Jaffard}). For any $d\in\set{1,\ldots,n-1}$, let 
\begin{equation*}
\zeta\paren{d,t} \coloneqq \set{\paren{r,s}:\; r\in S_n,\; j\notin S_n,\; \abs{r-s} = d}.
\end{equation*}
A simple counting argument leads to 
\begin{equation}\label{eq1Lem1A1}
\abs{\zeta\paren{d,t}} = \set{\begin{array}{cl}
d, &\; d\in \set{1,\ldots,t}\\
t,&\;  d\in\set{t,\ldots,n-t}\\
n-d,&\;  d\in\set{n-t+1,\ldots,n-1}\\
\end{array}}.
\end{equation}
Using triangle inequality and rearranging the different components in $\cc{\tau}\paren{L^{-1}_n,S_n, S^c_n}$, we get
\begin{equation}\label{eq2Lem1A1}
\abs{\cc{\tau}\paren{L^{-1}_n, S_n, S^c_n}}\leq \frac{1}{\sqrt{t\paren{n-t}}}\sum\limits_{d=1}^{n-1} \abs{\zeta\paren{d,t}} \max_{\paren{r,s}\in\zeta\paren{d,t}} \abs{L^{-1}_n\paren{r,s}}.
\end{equation}
According to Lemma A1 of \cite{P.Hall}, there is a $C\paren{f}$ for which the entries of $L^{-1}_n$ admits the following inequality,
\begin{equation}\label{eq3Lem1A1}
\abs{L^{-1}_n\paren{r,s}}\;\leq \frac{\tilde{C}\paren{f}}{\paren{1+\abs{r-s}}^{1+\lambda}}.
\end{equation}
We terminate the proof by substituting inequality \eqref{eq3Lem1A1} and identity \eqref{eq1Lem1A1} into inequality \eqref{eq2Lem1A1}. We will skip the algebraic details due to lack of space.
\begin{eqnarray*}
\abs{\cc{\tau}\paren{L^{-1}_n, S_n, S^c_n}} &\leq& \frac{\tilde{C}}{\sqrt{t\paren{n-t}}}\paren{\sum\limits_{d=1}^{t} \frac{d}{\paren{1+d}^{1+\lambda}}+\sum\limits_{d=t+1}^{n-t} \frac{t}{\paren{1+d}^{1+\lambda}}}+\sum\limits_{d=n-t+1}^{n-1} \frac{n-d}{\paren{1+d}^{1+\lambda}}\\
&\RelNum{\paren{a}}{\leq}&\frac{C}{\sqrt{\alpha\paren{1-\alpha}}}n^{-\paren{\lambda\wedge 1}} \paren{\bbM{1}_{\set{\lambda\ne 1}}+ \bbM{1}_{\set{\lambda=1}}\log n}.
\end{eqnarray*}
Note that inequality $\paren{a}$ follows from the fact that $\sqrt{t\paren{n-t}}\geq n\sqrt{\alpha\paren{1-\alpha}}$.
\end{proof}

The proof of the following result is omitted due to its analogy to the proof of Lemma \ref{UppBndTauSnScn}.

\begin{lem}\label{UppBndTauSnScn2}
With the same conditions and notation as Lemma \ref{UppBndTauSnScn}, there is a non-negative constant $C'\coloneqq C'\paren{f}$ such that
\begin{equation*}
\abs{\cc{\tau}\paren{L_n, S_n, S^c_n}} \;\leq \frac{C'}{\sqrt{\alpha\paren{1-\alpha}}}n^{-\paren{\lambda\wedge 1}} \paren{\bbM{1}_{\set{\lambda\ne 1}}+ \bbM{1}_{\set{\lambda=1}}\log n}.
\end{equation*} 
\end{lem}

\begin{prop}\label{UppBndTauSnInfMatrix}
Define the infinite column vector $\nu_n\in\bb{R}^{\bb{N}}$ by $\nu_n\paren{r} = \bbM{1}_{\set{r\leq n}}$. Then, there is a constant $C \coloneqq C\paren{f,\alpha}$ such that
\begin{equation*}
\abs{\frac{1}{n} \nu^\top_n L^{-1}_{\bb{N}}\nu_n-\frac{1}{f\left(0\right)}}\leq C \paren{n^{-\paren{\lambda-\epsilon}}\bbM{1}_{\set{\lambda\in\parbra{0,1}}}+n^{-1}\bbM{1}_{\set{\lambda>1}}},
\end{equation*}
in which $\epsilon\in\paren{0,\lambda}$ is chosen in an arbitrary way.
\end{prop}

\begin{proof}
We obtain an upper bound on $\paren{ \nu^\top_n L^{-1}_{\bb{N}}\nu_n/n-1/f\left(0\right)}$ and the lower bound can be achieved using akin techniques. For brevity, let $g = 1/f$. According to Proposition $1.12$ of \cite{A.Bottcher}, $L^{-1}_{\bb{N}}$ satisfies the following identity,

\begin{equation}\label{eq1Prop1App2}
L^{-1}_{\bb{N}} = \cc{T}_{\bb{N}}\paren{g} + L^{-1}_{\bb{N}} H_{\bb{N}}\paren{f}H_{\bb{N}}\paren{g},
\end{equation}
in which $H_{\bb{N}}\paren{f}$ is the generated \textit{Hankel} matrix by $f$ as
\begin{equation*}
H_{\bb{N}}\paren{f} = \paren{\begin{array}{cccc}
f_1 & f_2 & f_3 & \cdots \\
f_2 & f_3 & f_4 & \cdots \\
f_3 & f_4 & f_5 & \cdots \\
\vdots & \vdots & \vdots & \vdots
\end{array}}.
\end{equation*}
$H_{\bb{N}}\paren{g}$ can also be defined in a similar way. A simple algebraic manipulation on identity \ref{eq1Prop1App2} leads to
\begin{eqnarray}\label{eq2Prop1A2}
\frac{1}{n} \nu^\top_n L^{-1}_{\bb{N}}\nu_n &=& \frac{1}{n} \nu^\top_n \cc{T}_{\bb{N}}\paren{g}\nu_n+\frac{\nu^\top_n L^{-1}_{\bb{N}}H_{\bb{N}}\paren{f}H_{\bb{N}}\paren{g}\nu_n }{n} \nonumber \\
&=& \frac{1}{n} \nu^\top_n \cc{T}_{\bb{N}}\paren{g}\nu_n \;+\; \frac{1}{n}\sum\limits_{r=1}^{n}\sum\limits_{s=1}^{n} \InnerProd{L^{-1}_{\bb{N}}H_{\bb{N}}\paren{f}e_{\bb{N}}\paren{r}}{H_{\bb{N}}\paren{g}e_{\bb{N}}\paren{k}}\nonumber \\
&\RelNum{\paren{a}}{\leq}& \frac{1}{n} \nu^\top_n \cc{T}_{\bb{N}}\paren{g}\nu_n+ \frac{\OpNorm{L^{-1}_{\bb{N}}}{2}{2}}{n} \sum\limits_{r=1}^{n}\LpNorm{H_{\bb{N}}\paren{g}e_{\bb{N}}\paren{r}}{2}\sum\limits_{s=1}^{n}\LpNorm{H_{\bb{N}}\paren{f}e_{\bb{N}}\paren{s}}{2}\nonumber \\
&\leq& \frac{1}{n} \nu^\top_n \cc{T}_{\bb{N}}\paren{g}\nu_n+\frac{1}{nm_f} \paren{\sum\limits_{r=1}^{n}\LpNorm{H_{\bb{N}}\paren{g}e_{\bb{N}}\paren{r}}{2}\;\vee\;\sum\limits_{s=1}^{n}\LpNorm{H_{\bb{N}}\paren{f}e_{\bb{N}}\paren{s}}{2}}^2.\nonumber\\
\end{eqnarray}
Note that $\paren{a}$ is direct consequence of the combination of generalized Cauchy-Schwartz inequality and simple properties of operator norm; The last inequality of \ref{eq2Prop1A2} can be obtained by inequality $\paren{1.14}$ of \cite{A.Bottcher} on the operator norm of Toeplitz matrices. As next step, we control $\sum\limits_{r=1}^{n}\LpNorm{H_{\bb{N}}\paren{g}e_{\bb{N}}\paren{r}}{2}$ from above. It is known that if $m_f > 0$ and $f$ satisfies condition \ref{UppBndfk}, then $g$ does as well. So, there bounded are constants $c\paren{\lambda},c'\paren{\lambda},c''\paren{\lambda}>0$ such that
\begin{eqnarray}\label{eq3Prop1App2}
\sum\limits_{j=1}^{n}\LpNorm{H_{\bb{N}}\paren{g}e_{\bb{N}}\paren{r}}{2} &=& \sum\limits_{r=1}^{n}\sqrt{\sum\limits_{k=r}^{\infty} \abs{g_k}^2}\;\leq c\sum\limits_{r=1}^{n} \sqrt{\int\limits_{r}^{\infty}x^{-\paren{2+2\lambda}}dx} \;\leq\; c'\sum\limits_{r=1}^{n}j^{-\paren{1/2+\lambda}}\nonumber\\
&\leq&c''\paren{n^{1/2-\lambda}\bbM{1}_{\set{\lambda\in\paren{0,1/2}}} + \log n\bbM{1}_{\set{\lambda=1/2}}+ \bbM{1}_{\set{\lambda>1/2}}}.
\end{eqnarray}
Henceforth, there is some constant strictly positive $\tilde{c}$ (depending on $\lambda$) such that
\begin{equation}\label{UppBndIneq}
\frac{1}{n} \nu^\top_n L^{-1}_{\bb{N}}\nu_n \leq \frac{1}{n} \nu^\top_n \cc{T}_{\bb{N}}\paren{g}\nu_n+ \tilde{c}\paren{n^{-2\lambda}\bbM{1}_{\set{\lambda\in\paren{0,1/2}}} + n^{-1}\log^2 n\bbM{1}_{\set{\lambda=1/2}}+ n^{-1}\bbM{1}_{\set{\lambda>1/2}}}.
\end{equation}
Lastly, we obtain an upper bound on $\frac{1}{n} \nu^\top_n \cc{T}_{\bb{N}}\paren{g}\nu_n$. $\nu^\top_n \cc{T}_{\bb{N}}\paren{g}\nu_n$ can be viewed as variance of $S_n = \sum\limits_{k=1}^{n} X_k$ where $\set{X_k}_{k\in\bb{N}}$ is a stationary process with spectral density $g$. Thus, identity $2.1.2$ of \cite{L.Zhengyan} shows that

\begin{eqnarray}\label{UppBndIneq2}
\frac{1}{n} \nu^\top_n \cc{T}_{\bb{N}}\paren{g}\nu_n &=& \frac{1}{\pi n}\int\limits_{-\frac{\pi}{2}}^{\frac{\pi}{2}} \paren{\frac{\sin\paren{n\omega}}{\sin\paren{\omega}}}^2 g\paren{2\omega} d\omega =
\frac{g\paren{0}}{\pi n}\int\limits_{-\frac{\pi}{2}}^{\frac{\pi}{2}} \paren{\frac{\sin\paren{n\omega}}{\sin\paren{\omega}}}^2 d\omega\nonumber\\
&+& \frac{2}{\pi n}\int\limits_{-\frac{\pi}{2}}^{\frac{\pi}{2}} \paren{\frac{\sin\paren{n\omega}}{\sin\paren{\omega}}}^2 \brac{g\paren{2\omega}-g\paren{0}} d\omega \RelNum{\paren{a}}{=} \frac{1}{f\paren{0}}+\frac{2}{\pi n}\int\limits_{-\frac{\pi}{2}}^{\frac{\pi}{2}} \paren{\frac{\sin\paren{n\omega}}{\sin\paren{\omega}}}^2 \brac{g\paren{2\omega}-g\paren{0}} d\omega\nonumber\\
&\RelNum{\paren{b}}{\leq}& \frac{1}{f\paren{0}} + \frac{2\sum\limits_{m\in\bb{Z}}\abs{f_m} \abs{m}^{\lambda'}}{\pi nm_f f\paren{0}}\int\limits_{-\frac{\pi}{2}}^{\frac{\pi}{2}} \paren{\frac{\sin\paren{n\omega}}{\sin\paren{\omega}}}^2 \abs{\omega}^{\lambda'} d\omega.
\end{eqnarray}
Identity $\paren{a}$ follows from the following results which can be proved by applying \emph{Parseval's identity} on the triangular pulse centred at $0$. 
\begin{equation*}
\frac{1}{\pi n}\int\limits_{-\frac{\pi}{2}}^{\frac{\pi}{2}} \paren{\frac{\sin\paren{n\omega}}{\sin\paren{\omega}}}^2 d\omega = 1
\end{equation*}
Based upon the following inequality, $g$ inherits the Holder property from $f$.
\begin{equation*}
\abs{g\paren{2\omega}-g\paren{0}} = \abs{\frac{f\paren{2\omega}-f\paren{0}}{f\paren{0}f\paren{2\omega}}}\leq \frac{\abs{f\paren{2\omega}-f\paren{0}}}{f\paren{0}m_f}.
\end{equation*}
Hence, $\paren{b}$ is implied by this inheritance property and Proposition $3.2.12$ of \cite{L.Grafakos} regarding that $f$ is a $\paren{\lambda'}$-Holder function for any $0\leq\lambda'<\lambda$. 

For any $\lambda\in\parbra{0,1}$, choose $\epsilon = \epsilon\paren{\lambda}>0$ such that $\paren{\lambda-\epsilon} > 0$. Let $\lambda' = \paren{\lambda-\epsilon}\bbM{1}_{\set{\lambda\in\parbra{0,1}}} + \frac{1+\lambda}{2}\bbM{1}_{\set{\lambda > 1}}$. Obviously $\lambda' < \lambda$ and $\lambda'>1$ for any $\lambda>1$. Notice that $\sum\limits_{m\in\mathbb{Z}}\abs{f_m} \abs{m}^{\lambda'}$ is bounded, so,
\begin{eqnarray}\label{UppBndIneq3}
\frac{1}{n}\int\limits_{-\frac{\pi}{2}}^{\frac{\pi}{2}} \paren{\frac{\sin\paren{n\omega}}{\sin\paren{\omega}}}^2 \abs{\omega}^{\lambda'} d\omega &\RelNum{\paren{a_0}}{\asymp}& \int\limits_{0}^{\frac{\pi}{2}} \frac{\sin^2\paren{n\omega}}{\paren{n\omega}^2}\abs{\omega}^{\lambda'} dn\omega
= n^{-\lambda'} \int\limits_{0}^{\frac{n\pi}{2}} \frac{\sin^2\paren{u}}{u^{2-\lambda'}} du\nonumber\\
&\RelNum{\paren{a_1}}{\leq}&\hat{c} \paren{n^{-\lambda+\epsilon}\bbM{1}_{\set{\lambda\in\parbra{0,1}}}+n^{-1}\bbM{1}_{\set{\lambda>1}}}.
\end{eqnarray}
Inequality $\paren{a_0}$ follows from the fact that $\frac{2}{\pi}\leq\abs{\frac{\sin \omega}{\omega}} \leq 1$ for any $\omega\in\paren{0,\frac{\pi}{2}}$ and $\paren{a_1}$ is given by simple integration techniques. Combination of \eqref{UppBndIneq}-\eqref{UppBndIneq3} and  concludes the proof.
\end{proof}

\begin{prop}
Let $S_n = \set{1,\ldots,n}$. Under the same notation and assumptions as Proposition \ref{UppBndTauSnInfMatrix}, we have
\begin{equation*}
\abs{\frac{1}{n} \nu^\top_n L^{-1}_{\bb{N}}\nu_n-\cc{\tau}\paren{L^{-1}_{\bb{N}}, S_n}}\leq C \paren{n^{-\paren{\lambda-\epsilon}}\bbM{1}_{\set{\lambda\in\parbra{0,1}}}+n^{-1}\bbM{1}_{\set{\lambda>1}}}.
\end{equation*}
\end{prop}

\begin{proof}
Using \emph{Widom's theorem} (Theorem 2.14, \cite{A.Bottcher}) and one line of straightforward algebra, there is a matrix $D_n\in\mathbb{R}^{n\times n}$ such that
\begin{eqnarray*}
\abs{\frac{1}{n} \nu^\top_n L^{-1}_{\bb{N}}\nu_n-\cc{\tau}\paren{L^{-1}_{\bb{N}}, S_n}}&=& \abs{\frac{1}{n} \nu^\top_n V_n\paren{L^{-1}_{\bb{N}}-\cc{T}_{\bb{N}}\paren{g}}V_n\nu_n+ \cc{\tau}\paren{D_n, S_n}} \\
&\leq& \OpNorm{D_n}{2}{2} +\abs{\frac{1}{n}\nu^\top_n V_n\paren{L^{-1}_{\bb{N}}-\cc{T}_{\bb{N}}\paren{g}}V_n\nu_n}.
\end{eqnarray*}
$V_n$ is a bounded operator which is defined by $V_n\paren{\nu} = \paren{\nu_n,\cdots, \nu_1,0,0,\cdots}$ for any $\nu\in \ell^2$. Based upon Theorem $2.15$ of \cite{A.Bottcher} (first equation, p. 44), $\OpNorm{D_n}{2}{2} \;= o\paren{n^{-\lambda}}$. Moreover, the special form of $V_n$, gives

\begin{equation*}
\abs{\frac{1}{n}\paren{V_n\nu_n}^\top\paren{L^{-1}_{\bb{N}}-\cc{T}_{\bb{N}}\paren{g}}V_n\nu_n} = \abs{\frac{1}{n} \nu^\top_n \paren{L^{-1}_{\bb{N}}-\cc{T}_{\bb{N}}\paren{g}}\nu_n}.
\end{equation*}
Ultimately, combining the same tricks as inequality $\paren{a}$ in \eqref{eq2Prop1A2}, identity \eqref{eq1Prop1App2} and inequality \eqref{eq3Prop1App2} yields
\begin{equation*}
\abs{\frac{1}{n} \nu^\top_n V_n\paren{L^{-1}_{\bb{N}}-\cc{T}_{\bb{N}}\paren{g}}V_n\nu_n}\leq C'\paren{n^{-\paren{\lambda-\epsilon}}\bbM{1}_{\set{\lambda\in\parbra{0,1}}}+n^{-1}\bbM{1}_{\set{\lambda>1}}},
\end{equation*}
for some constant $C'\paren{\lambda}>0$. We end the proof by using triangle inequality.
\end{proof}

\begin{cor}\label{TauSigmaInvCor}
There is a constant $C>0$ depending on $f$ such that
\begin{equation}\label{eq5Cor1App2}
\abs{\cc{\tau}\paren{\Sigma^{-1}_n, \set{1,\ldots,n}} - \frac{1}{f\paren{0}}} \;\leq C'\paren{n^{-\paren{\lambda-\epsilon}}\bbM{1}_{\set{\lambda\in\parbra{0,1}}}+n^{-1}\bbM{1}_{\set{\lambda>1}}}.
\end{equation}
\end{cor}

\paragraph{Acknowledgment.}
This research is partially supported by NSF grant ACI-1047871. Additionally, CS is partially supported by NSF grants 1422157, 1217880, and 0953135, and LN
by NSF CAREER award DMS-1351362, NSF CNS-1409303, and NSF CCF-1115769.

\end{document}